\documentclass[12pt]{mythesis}
\usepackage{amsfonts}
\usepackage{amssymb}
\usepackage{amsmath}
\usepackage{graphicx}
\usepackage{epsfig}

\college{Jesus College}
\degree{Doctor of Philosophy}
\degreedate{Michaelmas 2005}
\newtheorem{theorem}{Theorem}

\newtheorem{corollary}[theorem]{Corollary}

\newtheorem{definition}[theorem]{Definition}
\newtheorem{example}[theorem]{Example}

\newtheorem{lemma}[theorem]{Lemma}

\newtheorem{proposition}[theorem]{Proposition}
\newtheorem{remark}[theorem]{Remark}

\newenvironment{proof}[1][Proof]{\textbf{#1.} }{\hfill $\Box$}

\begin{document}

\title{The Brownian Frame Process as a Rough Path}
\author{Ben Hoff}
\maketitle

\begin{abstract}
The Brownian frame process $\mathcal{T}^{B}$ is defined as 
\begin{equation*}
\mathcal{T}_{t}^{B}:=\left( B_{t-1+u}\right) _{0\leq u\leq 1},\text{ \ \ }%
t\in \left[ 0,1\right] ,
\end{equation*}%
where $B$ is a real-valued Brownian motion with parameter set $\left[ -1,1%
\right] $. This thesis investigates properties of the path-valued Brownian
frame process relevant to establishing an integration theory based on the
theory of rough paths (\cite{Lyons98}). The interest in studying this object
comes from its connection with Gaussian Volterra processes (e.g. \cite%
{Decreusefond05}) and stochastic delay differential equations (e.g. \cite%
{Mohammed84}). Chapter 2 establishes the existence of $\mathcal{T}^{B}$. We
then examine the convergence of dyadic polygonal approximations to $\mathcal{%
T}^{B}$ if the path-space $V$ where $\mathcal{T}^{B}$ takes its values is
equipped with first the $p$-variation norm ($p>2$) and second the $\sup $%
-norm. In the case of the $p$-variation norm, the Brownian frame process is
shown to have finite $\acute{p}$-variation for $\acute{p}>\frac{2p}{p-2}$.
In the case of the $\sup $-norm, it is shown to have finite $\acute{p}$%
-variation for $\acute{p}>2$. Chapter 3 provides a tail estimate for the
probability that two evaluations of the Brownian frame process are far apart
in the $p$-variation norm. Chapter 4 shows that $\mathcal{T}^{B}$ does not
have a L\'{e}vy area if $V\otimes V$ is equipped with the injective tensor
product norm (where $V=C\left( \left[ 0,1\right] \right) $).
\end{abstract}

\bigskip \begin{dedication}
To my parents.
\end{dedication}

\begin{acknowledgements}
My thanks go to my supervisor, Terry Lyons,  without whom this thesis would
never have happened. I am grateful to Ben Hambly and James Norris for their
careful reading and valuable comments. As a recipient of a doctoral training
grant I\ am endebted to \emph{EPSRC}. I gratefully acknowledge the
contribution of a \emph{Jesus College} \emph{Old Members Graduate Scholarship}. 
Most importantly, I want to thank my parents for their generosity in all respects.
\end{acknowledgements}

\begin{romanpages}
\tableofcontents
\end{romanpages}

\chapter{\protect\bigskip Introduction}

\section{Motivation}

Let $x_{.}:\left[ 0,t\right] \rightarrow V$, where $V~$is some linear space
and%
\begin{equation*}
T:\left[ 0,t\right] \rightarrow \left[ 0,t\right]
\end{equation*}%
be such that $T\left( h\right) \leq h$. By the \emph{historic frame path }$%
\mathcal{T}^{x}$ of $x$, we mean the path-valued path%
\begin{equation}
\mathcal{T}_{h}^{x}:=\left( x_{h-T\left( h\right) +u}\right) _{0\leq u\leq
T\left( h\right) },\text{ \ \ }h\in \left[ 0,t\right] .
\label{algebraic frame process definition}
\end{equation}%
$\mathcal{T}^{x}$ takes its values in the space of paths from $\left[
0,T\left( h\right) \right] $ to $V$, i.e. $\mathcal{T}_{h}^{x}\in V^{\left[
0,T\left( h\right) \right] }$ (where $V^{\left[ a,b\right] }$ denotes the
space of all paths from $\left[ a,b\right] $ to $V$). The map%
\begin{equation*}
T:\left[ 0,t\right] \rightarrow \left[ 0,t\right]
\end{equation*}%
determines the frame length of $\mathcal{T}^{x}$ at time $h$ so that the
evaluation of $\mathcal{T}^{x}$ at time $h$ is a $V$-valued path of length $%
T\left( h\right) $.

\begin{example}
\bigskip If $T\equiv 0$, then $\mathcal{T}^{x}\equiv x$.
\end{example}

\ We give two examples to illustrate the relevance of the \emph{historic
frame path:}

\begin{example}
From a control theory point of view, a differential equation is interpreted
as follows: The system we want to control has state $y_{t}$ at time $t$
(where $y_{t}$ lives on some manifold $W$). The initial state is given by $%
y_{0}$. The driving signal $x_{t}$ that determines the state of the system
at time $t$ lives on some Banach space $V$. It filters through a map 
\begin{equation*}
f\left( .\right) :V\rightarrow \left( W\rightarrow TW\right)
\end{equation*}%
to produce an effect in the state space $W$ (here $TW$ denotes the tangent
space of $W$). We call $f$ the \emph{vector field map} and write 
\begin{equation}
dy_{u}=f\left( y_{u}\right) dx_{u},\text{ \ \ }y_{0}=a.
\label{differential equation}
\end{equation}%
Subject to $f$ and $x$ satisfying certain conditions, \cite{Lyons98} gives
meaning to (\ref{differential equation}) (the exact result is quoted in
Theorem \ref{Universal Limit Theorem} below). In the classical theory of
dynamical systems with memory (\emph{delay differential equations}), the
system response at time $r$ may depend on the entire solution trajectory up
to time $r$, that is the vector field map $f$ is a function of the path
segment $\left( y_{u}:0\leq u\leq r\right) $, so that 
\begin{equation*}
f\left( \left( y_{u}\right) _{0\leq u\leq r}\right) :V\rightarrow TW.
\end{equation*}%
The \textquotedblleft delay\textquotedblright\ here occurs as an argument of
the vector field map. One could imagine a dynamical system with a different
type of memory where the delay occurs in the driving signal $x$: \emph{\ }%
Since each state $y_{r}$ is determined by the initial state $y_{0}$ and the
signal path segment $\left( x_{u}:0\leq u\leq r\right) $, an approach where $%
y$ is driven by the \emph{entire} historic trajectory $\left( x_{u}:0\leq
u\leq .\right) $ could prove to be interesting. In other words, we might
like to give meaning to (\ref{differential equation}) in the sense of \cite%
{Lyons98} if the driving signal $x$ evaluated at time $r$ equals $\left(
z_{u}\right) _{0\leq u\leq r}$ and $z$ is some $V$-valued path. In this
case, $x=\mathcal{T}^{z}$ is the frame process defined in (\ref{algebraic
frame process definition}) associated to $z$ with variable frame length $%
T\left( h\right) =h$.
\end{example}

\begin{example}
We consider a Volterra Gaussian Process, say Fractional Brownian Motion
(fBM) with Hurst parameter $H$. $B$ is a real-valued Brownian Motion. The
evaluation of fBM at time $t$ is obtained by applying the Wiener integral to
an appropriate deterministic $L^{2}$- kernel $K_{H}\left( t,.\right) $%
\footnote{%
The kernel $K_{H}$ is given by%
\begin{equation*}
K_{H}\left( t,s\right) =C_{H}\left[ \frac{2}{2H-1}\left( \frac{t\left(
t-s\right) }{s}\right) ^{H-\frac{1}{2}}-\int_{s}^{t}\left( \frac{u\left(
u-s\right) }{u}\right) ^{H-\frac{1}{2}}\frac{du}{u}\right] \mathbf{1}_{\left]
0,t\right[ }\left( s\right)
\end{equation*}%
where $C_{H}=\frac{\Gamma \left( 2-2H\right) \cos \left( \pi H\right) }{%
\Gamma \left( H-\frac{1}{2}\right) \pi H\left( 1-2H\right) }.$}, 
\begin{equation*}
fBM_{H}\left( r\right) :=\int_{0}^{r}K_{H}\left( t,s\right) dB_{s}.
\end{equation*}%
This is a functional of the path-segment $\mathcal{T}^{B}\left( r\right)
=\left( B_{u}:0\leq u\leq r\right) $ -- i.e. a functional of the Brownian
Frame Process of frame length $T\left( r\right) =r$.
\end{example}

\bigskip This thesis investigates some of the rough path properties (c.f. 
\cite{Lyons98}) of the \ \emph{historic frame path }on the Wiener space $%
\left( C_{0}\left( \left[ -1,1\right] \right) ,\sigma _{\left\Vert
.\right\Vert _{\infty }},\mathbb{P}\right) $\footnote{$C_{0}\left( \left[
-1,1\right] \right) $ denotes the space of continuous functions on $\left[
-1,1\right] $ that are $0$ at $-1$. $\sigma _{\left\Vert .\right\Vert
_{\infty }}$ denotes the completion of the $\sigma $-algebra generated by
the $\sup $-norm and $\mathbb{P}$ is the Wiener measure.}, defined for $f\in
C_{0}\left( \left[ -1,1\right] \right) $ at $h\in \left[ 0,1\right] $ as%
\begin{equation}
\mathcal{T}_{h}^{f}:=\left( f_{h-1+u}\right) _{0\leq u\leq 1}.
\label{Brownian frame process}
\end{equation}%
Here, we have a constant frame length $T\equiv 1$. When viewed as a (Borel)
random variable on Wiener space, $\mathcal{T}_{h}^{B}$ will be called the 
\emph{Brownian frame} random variable at $h$ (associated to the Brownian
Motion $B$ supported on $\left( C_{0}\left( \left[ -1,1\right] \right)
,\sigma _{\left\Vert .\right\Vert _{\infty }},\mathbb{P}\right) $).

\section{Rough Path Theory}

The theory of Rough Paths as developed in \cite{Lyons98} shows how to
construct solutions to differential equations driven by paths that are not
of bounded variation but have controlled roughness. The $p$-variation
(Definition\emph{\ }\ref{Defn:p-varn} below) is taken as a measure of a
Banach space valued path's roughness. The analysis is independent of the
dimension of the Banach space and so appears particularly well adapted to
the case where the Banach space in question is a path-space. We give a brief
overview:

Let $V$ be a Banach space with norm $\left\vert .\right\vert _{V}$. We
denote by $T^{n}\left( V\right) $ the truncated tensor algebra of $V$, that
is 
\begin{equation*}
T^{n}\left( V\right) :=\mathbb{R\oplus }V\oplus V^{\otimes 2}\oplus
...\oplus V^{\otimes n},
\end{equation*}%
where $V^{\otimes i}:=V\otimes V\otimes ...\otimes V$ ($i$ copies).

\begin{definition}[multiplicative functional]
A \emph{continuous} functional 
\begin{equation*}
\left\{ 
\begin{array}{c}
\mathbf{x}:\left\{ \left( s,t\right) :0\leq s\leq t\leq T\right\}
\rightarrow T^{n}\left( V\right) \\ 
\mathbf{x}_{s,t}=\left(
x_{s,t}^{0},x_{s,t}^{1},x_{s,t}^{2},...,x_{s,t}^{n}\right)%
\end{array}%
\right.
\end{equation*}%
is called \emph{multiplicative }if it satisfies \emph{Chen's identity}%
\begin{equation*}
\mathbf{x}_{s,t}\otimes \mathbf{x}_{t,u}=\mathbf{x}_{s,u}\text{ \ for }0\leq
s\leq t\leq u\leq T\text{ }
\end{equation*}%
and 
\begin{equation*}
x_{s,t}^{0}\equiv 0\text{ on }\left\{ \left( s,t\right) :0\leq s\leq t\leq
T\right\} .
\end{equation*}
\end{definition}

\begin{example}
\label{Ex: level 1 multiplicative lift}Every continuous $V$-valued path $x$
determines a multiplicative functional $\left( 1,x_{t}-x_{s}\right) $ in $%
T^{1}\left( V\right) $. The converse is not true: With every multiplicative
functional $y\in T^{1}\left( V\right) $ we may associate a collection of $V$%
-valued paths $\left\{ y_{0,t}^{1}+c:c\in V\right\} $.
\end{example}

A wide class of \emph{multiplicative functionals }is given by

\begin{theorem}[Chen's theorem]
\label{Theorem: Chen's theorem}Suppose $x$ is a bounded variation path in $V$%
. Then $x$ has a canonical \emph{multiplicative} lift $\mathbf{x}%
_{s,t}=\left( x_{s,t}^{0},x_{s,t}^{1},...,x_{s,t}^{n}\right) $ to $%
T^{n}\left( V\right) $ given by its sequence of iterated integrals:%
\begin{eqnarray*}
x_{s,t}^{0} &\equiv &1,\text{ \ }x_{s,t}^{1}\equiv x_{t}-x_{s},\text{ \ } \\
\text{\ }x_{s,t}^{k} &=&\int ...\int_{s\leq u_{1}\leq u_{2}\leq ...\leq
u_{k}\leq t}dx_{u_{1}}\otimes dx_{u_{2}}\otimes ...\otimes dx_{u_{k}}.
\end{eqnarray*}%
We call $\mathbf{x}$ the \emph{Chen lift of }$x$.
\end{theorem}

In general, if $\mathbf{y=}\left( y^{0},y^{1},y^{2},...,y^{m}\right) $ is a
multiplicative functional with values in $T^{m}\left( V\right) $ and for
some $\mathbf{z=}\left( z^{0},z^{1},z^{2},...,z^{n}\right) $ taking values
in $T^{n}\left( V\right) $ with $m\geq n$, we have that%
\begin{equation*}
\mathbf{z}=\left( y^{0},y^{1},y^{2},...,y^{n}\right) ,
\end{equation*}%
then we say that $\mathbf{z}$ is the projection of $\mathbf{y}$ onto $%
T^{n}\left( V\right) $. We also say that $\mathbf{y}$ is \emph{a}
multiplicative functional lying above $\mathbf{z}$ or that $\mathbf{y}$ is 
\emph{a lift} of $\mathbf{z}$. It is important to note that in general for a
given $\mathbf{z}$, neither existence nor uniqueness of $\mathbf{y}$ is
obvious (Theorem \ref{FirstTheoremLyons}\ below deals with the existence of
a unique lift for a particular class of $\mathbf{z}$).

For each $n$, $V^{\otimes n}$ is assumed to be equipped with a \emph{%
compatible tensor norm }$\left\Vert .\right\Vert _{V^{\otimes n}}$: If $v\in
V^{\otimes i}$ and $w\in V^{\otimes j}$, then $\left\Vert .\right\Vert
_{V^{\otimes n}}$ is said to be \emph{compatible }if for any $\left(
i,j\right) $ with $i+j\leq n$, we have 
\begin{equation}
\left\Vert v\otimes w\right\Vert _{V^{\otimes \left( i+j\right) }}\leq
\left\Vert v\right\Vert _{V^{\otimes i}}\left\Vert w\right\Vert _{V^{\otimes
j}},  \label{Defn: compatible norm}
\end{equation}%
and 
\begin{equation*}
\left\Vert v\right\Vert _{V^{\otimes 1}}\equiv \left\vert v\right\vert _{V}.
\end{equation*}

\begin{definition}
\label{Defn:p-varn}\bigskip\ Let $\mathcal{D}\left( \left[ 0,T\right]
\right) $ denote the set of all finite dissections of $\left[ 0,1\right] $,
that is 
\begin{equation*}
\mathcal{D}\left( \left[ 0,T\right] \right) :=\left\{ \left\{
t_{0},t_{1,},...,t_{n}\right\} :t_{0}=0<t_{1}<t_{2}<...<t_{n}=T\text{ and }n%
\text{ is finite}\right\} .
\end{equation*}%
For $p\geq n$, the $p$\emph{-variation functional of level }$i$ of a
function 
\begin{equation*}
\mathbf{x}:\left\{ \left( s,t\right) :0\leq s\leq t\leq T\right\}
\rightarrow T^{n}\left( V\right)
\end{equation*}%
with%
\begin{equation*}
\mathbf{x}_{s,t}=\left(
x_{s,t}^{0},x_{s,t}^{1},x_{s,t}^{2},...,x_{s,t}^{n}\right)
\end{equation*}%
is defined as 
\begin{equation*}
\mathcal{V}_{p}\left( x^{i}\right) :=\sup \left\{ \sum_{D}\left\Vert
x_{t_{j-1},t_{j}}^{i}\right\Vert _{V^{\otimes i}}^{\frac{p}{i}}:D\in 
\mathcal{D}\left( \left[ 0,T\right] \right) \right\} ^{\frac{i}{p}}.
\end{equation*}%
$\mathbf{x}$ is said to have finite $p$\emph{-variation} if%
\begin{equation*}
\mathcal{V}_{p}\left( \mathbf{x}\right) :=\max_{1\leq i\leq n}\mathcal{V}%
_{p}\left( x^{i}\right)
\end{equation*}%
is finite.\ 
\end{definition}

\begin{definition}[rough path]
A multiplicative functional of finite $p$-variation with values in $%
T^{\left\lfloor p\right\rfloor }\left( V\right) $ (where $\left\lfloor
p\right\rfloor :=\max \left( n\in \mathbb{N}:n\leq p\right) $) is called a
rough path of roughness $p$. The set of all $p$-rough paths is denoted as $%
\Omega _{p}\left( V\right) .$ For any $\mathbf{x,y}\in $ $\Omega _{p}\left(
V\right) $, the $p$-variation distance $d_{p}$ is defined as%
\begin{equation*}
d_{p}\left( \mathbf{x,y}\right) =\mathcal{V}_{p}\left( \mathbf{x-y}\right) 
\text{.}
\end{equation*}
\end{definition}

\begin{remark}
$\left( \Omega _{p}\left( V\right) ,d_{p}\right) $ is a complete metric
space as are all projections of $\Omega _{p}\left( V\right) $ onto $%
T^{n}\left( V\right) $ for $n\leq \left\lfloor p\right\rfloor $ (Lemma 3.3.3
in \cite{Lyons02}). However, $\Omega _{p}\left( V\right) $ is not a linear
space -- in general, the sum of two multiplicative functionals fails to be
multiplicative.
\end{remark}

\begin{lemma}
If $\acute{p}\geq p$, then $d_{\acute{p}}\leq d_{p}$.
\end{lemma}

\begin{proof}
For $\left( a_{j}\right) _{1\leq j\leq n}\in \mathbb{R}^{n}$, 
\begin{eqnarray*}
\left[ \sum_{j=1}^{n}\left\vert a_{j}\right\vert \right] ^{\frac{\acute{p}}{p%
}} &=&\left[ \sum_{j=1}^{n}\left\vert a_{j}\right\vert \right] \left[
\sum_{j=1}^{n}\left\vert a_{j}\right\vert \right] ^{\frac{\acute{p}}{p}-1} \\
&=&\sum_{j=1}^{n}\left\vert a_{j}\right\vert \left[ \sum_{i=1}^{n}\left\vert
a_{i}\right\vert \right] ^{\frac{\acute{p}}{p}-1}\geq
\sum_{j=1}^{n}\left\vert a_{j}\right\vert ^{\frac{\acute{p}}{p}}.
\end{eqnarray*}%
We now fix $\mathbf{x}\in \Omega _{p}\left( V\right) $ and a dissection $D$
in $\mathcal{D}\left( \left[ 0,T\right] \right) $ (Definition \ref%
{Defn:p-varn} above). If we take $\left\vert a_{j}\right\vert =\left\vert
x_{t_{j},t_{j+1}}^{i}\right\vert ^{\frac{p}{i}}$ for $1\leq i\leq p$ and
then take the $\sup $ over $\mathcal{D}\left( \left[ 0,T\right] \right) $,
the result follows (we set $x_{s,t}^{n}\equiv 0$ for $p<n\leq \acute{p}$).
\end{proof}

\begin{definition}[smooth rough path]
\label{Defn: smooth rough path}Suppose $x$ is a bounded variation path in $V$%
. We will denote by $\mathbf{x}$ the Chen lift of $x$ to $T^{n}\left(
V\right) $ and call $\mathbf{x}$ a \emph{smooth rough path}. The collection
of smooth rough paths is contained in $\cap _{p\geq 1}\Omega _{p}\left(
V\right) $.
\end{definition}

\begin{definition}[geometric rough path]
\label{Defn: geometric rough path}The closure of the smooth rough paths
under $d_{p}$ is called the space of \emph{geometric rough paths\ of
roughness }$p$ and denoted $\Omega G\left( V\right) _{p}$.
\end{definition}

The \emph{\textquotedblleft First Theorem\textquotedblright\ }from \cite%
{Lyons98} gives sufficient conditions for the existence of a unique lift:

\begin{theorem}[Theorem 2.2.1 in \protect\cite{Lyons98}]
\label{FirstTheoremLyons}Let $\mathbf{X}_{s,t}^{\left( n\right) }$ be a
continuous multiplicative functional in $T^{\left( n\right) }\left( V\right) 
$ of finite $p$-variation where $n=\left\lfloor p\right\rfloor $. There
exists a multiplicative extension $\mathbf{X}_{s,t}^{\left( m\right) }$ to $%
T^{\left( m\right) }\left( V\right) $, $m>n$ which is of finite $p$%
-variation. The extension is unique in this class.
\end{theorem}

The \emph{Universal Limit Theorem }below establishes the connection between
multiplicative functionals and solutions to differential equations:

\begin{theorem}[ Theorem 4.1.1 in \protect\cite{Lyons98}]
\label{Universal Limit Theorem}Let $f$ be a linear map from $V$ to the space
of Lipschitz-$\gamma $ vector fields on $W$, that is $f\in L\left(
V,Lip\left( \gamma ,W,W\right) \right) $ \footnote{$Lip\left( n+\varepsilon
,W,W\right) $ denotes the space of $n$-times differentiable functions from $%
W $ to $W$ whose $n^{\text{th}}$ derivative is $\varepsilon $-Lipschitz}%
where $1\leq p<\gamma $. Then consider the It\^{o} map $I:X\rightarrow
\left( X,Y\right) $ defined for bounded variation paths by 
\begin{equation}
dY_{t}=f\left( Y_{t}\right) dX_{t},\text{ \ \ \ }Y_{0}=a.  \label{Ch1diffeqn}
\end{equation}%
Define the one form $h$ by%
\begin{equation*}
h\left( x,y\right) \left( dX,dY\right) =h\left( y\right) \left( dX,dY\right)
=\left( dX,f\left( y\right) dX\right) .
\end{equation*}%
For any geometric multiplicative functional $\mathbf{X}\in \Omega G\left(
V\right) _{p}$ there is exactly one geometric multiplicative functional
extension $\mathbf{Z=}\left( \mathbf{X,Y}\right) \in \Omega G\left( V\oplus
W\right) _{p}$ such that if $Y_{t}=$ $\mathbf{Y}_{0,t}^{1}+a,$ then $\mathbf{%
Z}$ satisfies the rough\footnote{%
By a \emph{rough }differential equation we mean a differential equation on
the truncated tensor algebra $T^{n}\left( V\right) $. To emphasise the
difference between a differential equation and a \emph{rough} differential
equation we write \textquotedblleft $\delta \mathbf{Z}$\textquotedblright\
instead of \textquotedblleft $d\mathbf{Z}$\textquotedblright (c.f Definition
4.1.1 in \cite{Lyons98})\textbf{\ }} differential equation%
\begin{equation*}
\delta \mathbf{Z=}h\left( \mathbf{Y}_{t}\right) \delta \mathbf{Z}.
\end{equation*}%
Such a solution exists on a small interval $\left[ 0,T\right] $ whose length
can be controlled entirely in terms of the control of the roughness of $X$
and of $f$ . The It\^{o} map is uniformly continuous and the map $\mathbf{%
X\rightarrow Z}$ is the unique continuous extension of the It\^{o} map from $%
\Omega G\left( V\right) ^{p}$ to $\Omega G\left( V\oplus W\right) ^{p}$.
\end{theorem}

\section{Results}

In Chapter $2$ we establish the existence of a continuous Brownian frame
process $\mathcal{T}^{B}$ into $\left( C\left( \left[ 0,1\right] \right)
,\left\Vert .\right\Vert _{\infty }\right) $. $\mathcal{T}^{B}$
algebraically agrees with the object defined in (\ref{Brownian frame process}%
), i.e. for any $f\in C_{0}\left( \left[ -1,1\right] \right) $, $\mathcal{T}%
^{B}\left( f\right) :=\mathcal{T}^{f}$. We call $\mathcal{T}^{B}$ the
Brownian $\sup $\emph{-frame process}.\emph{\ }For\ $\acute{p}>2,$ $\mathcal{%
T}_{.}^{B}$ is shown to have finite $\acute{p}$-variation. Furthermore, the
dyadic polygonal approximations of $\mathcal{T}^{B}$ converge to $\mathcal{T}%
^{B}$ in $\acute{p}$-variation.

If $p\geq 1$, $\left( C\left( \left[ 0,1\right] \right) _{p},\left\Vert
.\right\Vert _{p}\right) $ denotes the Banach space of continuous functions
of finite $p$-variation equipped with the $p$-variation norm 
\begin{equation*}
\left\Vert .\right\Vert _{p}:=\left\Vert .\right\Vert _{\infty }+\mathcal{V}%
_{p}\left( .\right) .
\end{equation*}%
For $p>2$, we establish the existence of a continuous Brownian frame process 
$\mathcal{S}_{.}^{B,p}$ into $\left( C\left( \left[ 0,1\right] \right)
_{p},\left\Vert .\right\Vert _{p}\right) $ that algebraically agrees with
the object defined in (\ref{Brownian frame process}) off a nullset. That is, 
$\mathcal{S}^{B,p}\left( f\right) $ is equal to $\mathcal{T}^{f}$ for $%
\mathbb{P}$-a.e. $f\in C_{0}\left( \left[ -1,1\right] \right) .$ \ $\mathcal{%
S}^{B,p}$ is called the $p$\emph{-variation frame process. }For $\acute{p}>%
\frac{2p}{p-2}$, we show that the dyadic polygonal approximations of $%
\mathcal{S}^{B,p}$ converge to $\mathcal{S}_{.}^{B,p}$ in $\acute{p}$%
-variation.

In Chapter $3$, we examine further properties of $\mathcal{S}^{B,p}$ ($p>2$%
): We find constants, $d_{1}\left( \alpha ,p\right) $ \ ($\alpha $ is a
constant strictly greater than $1-\frac{1}{p}$) and $d_{2}\left( p\right) $,
so that the random variable 
\begin{equation*}
\frac{\left\Vert \mathcal{S}_{h_{2}}^{B,p}-\mathcal{S}_{h_{1}}^{B,p}\right%
\Vert _{p}}{d_{2}\left( p\right) \left( h_{2}-h_{1}\right) ^{\frac{1}{2}-%
\frac{1}{p}}}-d_{1}\left( \alpha ,p\right)
\end{equation*}%
has Gaussian tails. More precisely, we find constants $d_{1}\left( \alpha
,p\right) $ and $d_{2}\left( p\right) $ such that 
\begin{equation*}
\mathbb{P}\left( \frac{\left\Vert \mathcal{S}_{h_{2}}^{B,p}-\mathcal{S}%
_{h_{1}}^{B,p}\right\Vert }{d_{2}\left( p\right) \left( h_{2}-h_{1}\right) ^{%
\frac{1}{2}-\frac{1}{p}}}-d_{1}\left( \alpha ,p\right) >r\right) \leq \frac{1%
}{\sqrt{2\pi }r}\exp \left( -\frac{r^{2}}{2}\right) .
\end{equation*}

In Chapter $4$ we are concerned with the $\sup $-frame process $\mathcal{T}%
^{B}$. The \emph{Universal Limit Theorem} (Theorem \ref{Universal Limit
Theorem} above) tells us how to solve differential equations driven by rough
paths --\emph{\ provided the driving signal is the level }$1$ \emph{%
projection of some geometric }$\acute{p}$\emph{-rough path. }Since 
\begin{equation*}
\mathcal{T}^{B}:\left[ 0,1\right] \rightarrow C\left( \left[ 0,1\right]
,\left\Vert .\right\Vert _{\infty }\right)
\end{equation*}%
has finite $\acute{p}$-variation for $\acute{p}>2$, by Theorem \ref%
{FirstTheoremLyons}, any lift of $\mathcal{T}^{B}$ to $T^{2}\left( V\right) $
that preserves finite $\acute{p}$-variation establishes an integration
theory for $\mathcal{T}^{B}$ (here, $V=\left( C\left( \left[ 0,1\right]
\right) ,\left\Vert .\right\Vert _{\infty }\right) $). From Chapter 2 we
know that the (smooth) dyadic polygonal approximations $\mathcal{T}%
^{B}\left( m\right) $ of $\mathcal{T}^{B}$ converge to $\mathcal{T}^{B}$ in $%
\acute{p}$-variation norm. So it is natural to ask whether the lifts of $%
\mathcal{T}^{B}\left( m\right) $ converge in $T^{2}\left( V\right) $,
thereby providing a natural lift of $\mathcal{T}^{B}$ to $T^{2}\left(
V\right) $. We show that while the canonical lifts of the dyadic polygonal
approximations to $\mathcal{T}^{B}$ converge, the convergence is to an
object that lives \emph{outside} the injective tensor product $V\otimes
_{\vee }V$: After proving that $V\otimes _{\vee }V$ is isomorphic to $%
C\left( \left[ 0,1\right] \times \left[ 0,1\right] \right) $, we show that
the L\'{e}vy Area $\mathcal{A}_{0,1}\left( \mathcal{T}^{B}\right) $ of $%
\mathcal{T}^{B}$ -- which is the limit of the antisymmetric component of $%
\int \int_{0\leq u\leq v\leq 1}d\mathcal{T}^{B}\left( m\right) _{u}\otimes d%
\mathcal{T}^{B}\left( m\right) _{v}$ as $m\rightarrow \infty $ -- is
continuous \emph{off} the diagonal of the unit square. However, it is shown
to have a jump-discontinuity \emph{on} the diagonal and so does not exist in 
$V\otimes _{\vee }V$.

We believe that a further examination of the \textquotedblleft
obstruction\textquotedblright\ -- which is intimately linked to the
quadratic variation of Brownian motion -- will prove interesting and
important.

\chapter{\protect\bigskip The Brownian frame process}

\section{Notation}

We work on the classical Wiener Space $\left( C_{0}\left( \left[ -1,1\right]
,\mathbb{R}\right) ,\sigma _{\left\Vert .\right\Vert _{\infty }},\mathbb{P}%
\right) $. Here, 
\begin{equation*}
C_{0}\left( \left[ -1,1\right] \right) :=\left\{ f\in C\left( \left[ -1,1%
\right] \right) :f\left( -1\right) =0\right\}
\end{equation*}%
and $\sigma _{\left\Vert .\right\Vert _{\infty }}$ denotes the completion of
the Borel $\sigma $-algebra generated by the sup-norm 
\begin{equation*}
\left\Vert f\right\Vert _{\infty }:=\sup \left\{ \left\vert f\left( t\right)
\right\vert :t\in \left[ -1,1\right] \right\} ,
\end{equation*}%
with respect to the Wiener measure $\mathbb{P}$.

$B$ denotes the coordinate process 
\begin{equation*}
B_{t}:C_{0}\left( \left[ -1,1\right] \right) \rightarrow \mathbb{R}%
:B_{t}\left( f\right) :=f\left( t\right) \text{ \ \ },t\in \left[ -1,1\right]
\end{equation*}
(on $C_{0}\left( \left[ -1,1\right] \right) $) which (under Wiener measure $%
\mathbb{P}$) is a Brownian Motion.

$C\left( \left[ a,b\right] \right) _{p}$ denotes the space of continuous
real-valued paths on $\left[ a,b\right] $ of finite $p$-variation, i.e.%
\begin{equation*}
C\left( \left[ a,b\right] \right) _{p}=\left\{ f\in C\left( \left[ a,b\right]
\right) :\mathcal{V}_{p}\left( f^{1}\right) <\infty \right\} ,
\end{equation*}%
where $f^{1}\left( s,t\right) :=f\left( t\right) -f\left( s\right) $ (c.f.
Theorem \ref{Theorem: Chen's theorem} in Chapter 1) and the $p$-variation
functional of level $1$ -- $\mathcal{V}_{p}\left( .\right) $ -- is defined
in Definition \ref{Defn:p-varn} in Chapter 1.

Similarly, $C_{0}\left( \left[ a,b\right] \right) _{p}$ denotes the space of
continuous functions of finite $p$-variation that are $0$ at $a$, i.e. 
\begin{equation*}
C_{0}\left( \left[ a,b\right] \right) _{p}=\left\{ f\in C_{0}\left( \left[
a,b\right] \right) :\mathcal{V}_{p}\left( f^{1}\right) <\infty \right\} .
\end{equation*}

When $f\in C\left( \left[ a,b\right] \right) _{p}$, we say that $f$ has
finite $p$-variation. We write $\mathcal{V}_{p}\left( f\right) $ in place of 
$\mathcal{V}_{p}\left( f^{1}\right) $ and in this way view $\mathcal{V}%
_{p}\left( .\right) $ as a functional on $C\left( \left[ a,b\right] \right)
_{p}$ .

From Chapter 1 we recall that $\mathcal{V}_{p}\left( .\right) +\left\Vert
.\right\Vert _{\infty }$ is the $p$-variation norm on $C\left( \left[ 0,1%
\right] \right) _{p}$ for which we write%
\begin{equation*}
\left\Vert .\right\Vert _{p}:=\mathcal{V}_{p}\left( .\right) +\left\Vert
.\right\Vert _{\infty }.
\end{equation*}%
For the Banach space $C\left( \left[ 0,1\right] \right) _{p}$ equipped with $%
\left\Vert .\right\Vert _{p}$ we write $\left( C\left( \left[ 0,1\right]
\right) _{p},\left\Vert .\right\Vert _{p}\right) $. When the norm is not
explicitly mentioned, $C\left( \left[ 0,1\right] \right) _{p}$ is viewed as
a subspace of $\left( C\left( \left[ 0,1\right] \right) ,\left\Vert
.\right\Vert _{\infty }\right) $.

For $\alpha >1-\frac{1}{p}$ and $p>1$, we define the function $c\left(
.,.\right) $ as%
\begin{equation}
c\left( \alpha ,p\right) :=\left( \sum_{n=1}^{\infty }n^{-\alpha \frac{p}{p-1%
}}\right) ^{\frac{p-1}{p}}.  \label{c(alpha,p)-definition}
\end{equation}

The evaluation at $h\in \left[ 0,1\right] $ of the (deterministic) \emph{%
frame operator }$\mathcal{T}_{.}\mathcal{\ }$on the set of functions from $%
\left[ -1,1\right] $ to $\mathbb{R}$ -- denoted as $\mathbb{R}^{\left[ -1,1%
\right] }$ -- is defined as%
\begin{equation*}
\left\{ 
\begin{array}{c}
\mathcal{T}_{h}^{.}:\mathbb{R}^{\left[ -1,1\right] }\rightarrow \mathbb{R}^{%
\left[ 0,1\right] } \\ 
\mathcal{T}_{h}^{f}=\left( f\left( h-1+u\right) \right) _{u\in \left[ 0,1%
\right] }%
\end{array}%
\right. .
\end{equation*}%
We will be interested in $\mathcal{T}_{h}$ as a function on Wiener space and
so will be considering the restriction of $\mathcal{T}_{h}$ to $C_{0}\left( %
\left[ -1,1\right] \right) $.

If $M=\left\{ x_{1},x_{2},...,x_{n}\right\} $, we write $\left\Vert
M\right\Vert $ to denote the number of elements in $M$, i.e. $\left\Vert
M\right\Vert =n$.

The $\Gamma $-function is defined as $\Gamma \left( t\right)
:=\int_{0}^{\infty }x^{t-1}e^{-x}dx$.

$D\left( \left[ 0,1\right] ^{d}\right) $ is the $d$-dimensional unit cube of
dyadic rationals, i.e. 
\begin{equation}
D\left( \left[ 0,1\right] ^{d}\right) :=\left( \bigcup_{n\in \mathbb{N}%
}\left\{ \frac{k}{2^{n}}:0\leq k\leq 2^{n}\right\} \right) ^{d},
\label{Eqn: dyadic n-cube}
\end{equation}%
and for $u,t\in \left[ 0,1\right] ^{d}$ , \textquotedblleft $u\geq t$%
\textquotedblright\ means that 
\begin{equation*}
u_{i}\geq t_{i}\text{ \ \ \ for \ }1\leq i\leq d\text{.}
\end{equation*}

Finally, for any $h\in \left[ 0,1\right] $, $n\left( h\right) $ denotes the
unique integer such that $2^{-n\left( h\right) }\leq h<2^{-n\left( h\right)
+1}$.

\section{Main Results}

In the following section, we show that for every fixed $h\in \left[ 0,1%
\right] $, $\left. \mathcal{T}_{h}\right\vert _{C_{0}\left( \left[ -1,1%
\right] \right) }$ is a Borel random variable on Wiener space, mapping into $%
\left( C\left( \left[ 0,1\right] \right) ,\left\Vert .\right\Vert _{\infty
}\right) $. In order to distinguish between the algebraic frame operator $%
\mathcal{T}_{h}$ (which we recall is into $\mathbb{R}^{\left[ 0,1\right] }$
without a topology) and the analytic object that is the Borel random
variable into $\left( C\left( \left[ 0,1\right] \right) ,\left\Vert
.\right\Vert _{\infty }\right) $, we write 
\begin{equation*}
\mathcal{T}_{h}^{B}:=\left. \mathcal{T}_{h}^{.}\right\vert _{C_{0}\left( %
\left[ -1,1\right] \right) }.
\end{equation*}%
We show that $\mathcal{T}^{B}$ has continuous sample paths. As mentioned in
Chapter $1$, $\mathcal{T}^{B}$ is called the $\sup $-frame-process.

For $p>2,$ we show that for any fixed $h\in \left[ 0,1\right] $ the function 
\begin{equation}
\mathcal{T}_{h}^{B,p}:=\left\{ 
\begin{array}{ll}
\mathcal{T}_{h}^{B} & \text{on }C_{0}\left( \left[ -1,1\right] \right) _{p}
\\ 
0 & \text{otherwise}%
\end{array}%
\right.  \label{p-variation frame process}
\end{equation}%
is a Borel random variable on Wiener space.

We produce H\"{o}lder-type moment bounds for $\left\Vert \mathcal{T}%
_{h_{2}}^{B}-\mathcal{T}_{h_{1}}^{B}\right\Vert _{\infty }$ and $\left\Vert 
\mathcal{T}_{h_{2}}^{B,p}-\mathcal{T}_{h_{1}}^{B,p}\right\Vert _{p}$. In the
case of $\mathcal{T}^{B,p}$, we use Kolmogorov's lemma to show that there
exists a modification with continuous sample paths. We denote this
continuous modification as $\mathcal{S}^{B,p}$ and call it the $p$-variation
frame process. We then show that if $\acute{p}>\frac{2p}{p-2}$, the dyadic
polygonal approximations of $\mathcal{S}^{B,p}$ converge to $\mathcal{S}%
^{B,p}$ in $\acute{p}$-variation,\ $\mathbb{P}$-a.s. In the case of $%
\mathcal{T}^{B}$ we deduce that $\mathcal{T}^{B}$ has finite $\acute{p}$%
-variation for $\acute{p}>2$, $\mathbb{P}$-a.s, and that the polygonal
dyadic approximations of $\mathcal{T}^{B}$ converge to $\mathcal{T}^{B}$ in $%
\acute{p}$-variation, $\mathbb{P}$-a.s.

\section{\protect\bigskip Existence of the Brownian Frame random variable}

\subsection{\protect\bigskip The $\sup $-norm frame process $\mathcal{T}^{B}$%
}

\begin{lemma}
\label{Lemma: sup-norm process is Borel}For any $h\in \left[ 0,1\right] $, $%
\mathcal{T}_{h}^{B}$ is continuous and hence Borel measurable.
\end{lemma}

\begin{proof}
Fix $\varepsilon >0$. Suppose $f,g\in C_{0}\left( \left[ -1,1\right] \right) 
$ and that $\left\Vert f-g\right\Vert _{\infty }<\delta =\varepsilon .$ Then 
\begin{equation*}
\sup_{u\in \left[ 0,1\right] }\left\vert f\left( h-1+u\right) -g\left(
h-1+u\right) \right\vert \leq \left\Vert f-g\right\Vert _{\infty
}<\varepsilon \text{.}
\end{equation*}
\end{proof}

\begin{lemma}
The $\sup $-norm frame process $\mathcal{T}^{B}$ has continuous sample paths.
\end{lemma}

\begin{proof}
Fix $f\in C_{0}\left( \left[ -1,1\right] \right) $ and $\varepsilon >0$.
Since $\left[ -1,1\right] $ is compact, $f$ is uniformly continuous, so that 
$\exists \delta _{\varepsilon }>0$ such that $\left\vert x-y\right\vert
<\delta _{\varepsilon }\Longrightarrow \left\vert f\left( x\right) -f\left(
y\right) \right\vert <\varepsilon $, i.e.%
\begin{equation*}
\sup_{x,y\in \left[ -1,1\right] :\left\vert x-y\right\vert <\delta
_{\varepsilon }}\left\vert f\left( y\right) -f\left( x\right) \right\vert
\leq \varepsilon \text{,}
\end{equation*}%
and so for $\left\vert h_{2}-h_{1}\right\vert <\delta _{\varepsilon }$, 
\begin{eqnarray*}
&&\sup_{u\in \left[ 0,1\right] }\left\vert \mathcal{T}_{h_{2}}^{f}\left(
u\right) -\mathcal{T}_{h_{1}}^{f}\left( u\right) \right\vert \\
&=&\sup_{u\in \left[ 0,1\right] }\left\vert f\left( h_{2}-1+u\right)
-f\left( h_{1}-1+u\right) \right\vert \\
&\leq &\sup_{x,y\in \left[ -1,1\right] :\left\vert x-y\right\vert <\delta
_{\varepsilon }}\left\vert f\left( y\right) -f\left( x\right) \right\vert
\leq \varepsilon .
\end{eqnarray*}
\end{proof}

\subsection{The $p$-variation-norm frame process $\mathcal{S}^{B,p}$}

We regard the $p$-variation functional $\mathcal{V}_{p}$ as a functional on $%
\left( C\left( \left[ 0,1\right] \right) _{p},\left\Vert .\right\Vert
_{\infty }\right) $ and prove that it is lower semi-continuous as a map%
\begin{equation*}
\mathcal{V}_{p}:\left( C\left( \left[ 0,1\right] \right) _{p},\left\Vert
.\right\Vert _{\infty }\right) \rightarrow \left( \mathbb{R}^{+},\left\vert
.\right\vert \right) .
\end{equation*}

\begin{definition}
A real valued function $f$ defined on a topological space $\left( X,\tau
_{X}\right) $ is \emph{lower semi-continuous }if for any $\alpha \in \mathbb{%
R}$ the set $\left\{ x\in X:f\left( x\right) \leq \alpha \right\} $ is
closed.
\end{definition}

So any lower-semicontinuous real-valued function is Borel-measurable.\emph{\ 
}

\begin{lemma}
\label{lsc_lemma}A real-valued function $f$ is \emph{lower semi-continuous }%
iff \ $x_{n}\rightarrow x$, then $\lim \inf f\left( x_{n}\right) \geq
f\left( x\right) $.
\end{lemma}

\begin{proof}
$\Longrightarrow $: Fix any $\alpha \in \mathbb{R}$ and consider $M_{\alpha
}=\left\{ f\left( x\right) \leq \alpha \right\} $. Suppose $\left\{
x_{n}:n\in \mathbb{N}\right\} \subset M_{\alpha }$ and $x_{n}\rightarrow x$.
Since $f\left( x_{n}\right) \leq \alpha $ for all $n\in \mathbb{N}$, it
follows that $\lim \inf f\left( x_{n}\right) \leq \alpha $. But by
assumption $f\left( x\right) \leq \lim \inf f\left( x_{n}\right) \leq \alpha 
$, so that $x\in M_{\alpha }$. Thus, $M_{\alpha }$ is closed.

$\Longleftarrow $: Let $m:=\lim \inf f\left( x_{n}\right) $ and fix $%
\varepsilon >0$. An infinite subsequence of $\left( x_{n}\right) $ is
contained in $A_{\varepsilon }:=\left\{ x\in X:f\left( x\right) \leq
m+\varepsilon \right\} $. Hence, $x\in A_{\varepsilon }$ for all $%
\varepsilon >0$. Hence, $f\left( x\right) \leq m$.
\end{proof}

\begin{proposition}
\label{Propn: Semicontinuity of p-variation}We recall the definition of the $%
p$-variation functional (c.f. Definition \ref{Defn:p-varn}): Let $\mathcal{D}%
\left( \left[ 0,1\right] \right) $ denote the set of all finite dissections
of $\left[ 0,1\right] $. If $p\geq 1$, the $p$-variation functional $%
\mathcal{V}_{p}$ on $C\left( \left[ 0,1\right] \right) _{p}$, defined as%
\begin{equation*}
\mathcal{V}_{p}:C\left( \left[ 0,1\right] \right) _{p}\rightarrow \mathbb{R}%
^{+}:\mathcal{V}_{p}\left( f\right) :=\sup \left\{ \sum_{D}\left\vert
f\left( t_{i+1}\right) -f\left( t_{i}\right) \right\vert ^{p}:D\in \mathcal{D%
}\left( \left[ 0,1\right] \right) \right\} ^{\frac{1}{p}}
\end{equation*}%
is lower-semicontinuous.
\end{proposition}

\begin{proof}
Fix $f\in C\left( \left[ 0,1\right] \right) _{p}$. Suppose that $f_{n}\in
C\left( \left[ 0,1\right] \right) $ and that $f_{n}\rightarrow f$. \ Since $%
\left[ 0,1\right] $ is compact and $f$ and $\left\{ f_{n}\,:n\in \mathbb{N}%
\right\} $ are continuous, we may assume w.l.o.g. that the sequence $\left(
f_{n}\right) $ is uniformly bounded by $K<\infty $, i.e. that $\sup_{n\in 
\mathbb{N}
}\left\Vert f_{n}\right\Vert _{\infty }<K$.\ Since $f$ \ has finite $p$%
-variation, $\mathcal{V}_{p}^{p}\left( f\right) =c<\infty $, there exists a
sequence of dissections $\left( D_{m}\right) _{m\in 
\mathbb{N}
}$ with 
\begin{equation*}
D_{m}=\left\{ t_{0}^{\left( m\right) }=0;t_{1}^{\left( m\right)
};t_{2}^{\left( m\right) };...;t_{\left\Vert D_{m}\right\Vert }^{\left(
m\right) }=1\right\} \in \mathcal{D}\left( \left[ 0,1\right] \right)
\end{equation*}%
such that%
\begin{equation*}
\mathcal{V}_{p}^{p}\left( \sum_{D_{m}}f\left( t_{j}^{\left( m\right)
}\right) \mathbf{1}_{\left[ t_{j}^{\left( m\right) },t_{j+1}^{\left(
m\right) }\right) }\right) \rightarrow c.
\end{equation*}%
In other words: the sequence of step-functions defined as 
\begin{equation*}
\hat{f}_{m}=\sum_{j=0}^{\left\Vert D_{m}\right\Vert -1}f\left( t_{j}^{\left(
m\right) }\right) \mathbf{1}_{\left[ t_{j}^{\left( m\right)
},t_{j+1}^{\left( m\right) }\right) },
\end{equation*}%
converges to $f$ in $p$-variation. We start by fixing $\varepsilon >0$.
W.l.o.g. we may assume that $\left\Vert D_{m}\right\Vert $ is non-decreasing
(otherwise we eventually end up with the trivial dissection $\left\{
0,1\right\} $) so that we may choose 
\begin{equation*}
\delta _{m}:=\frac{\varepsilon }{\left\Vert D_{m}\right\Vert 2^{p}pK^{p-1}}.
\end{equation*}%
By definition of $\hat{f}_{m}$, $\exists N_{\varepsilon }^{1}$ such that 
\begin{equation*}
\inf_{m>N_{\varepsilon }^{1}}\mathcal{V}_{p}^{p}\left( \hat{f}_{m}\right)
>c-\varepsilon ,
\end{equation*}%
and since $f_{n}\rightarrow f$ in $C\left( \left[ 0,1\right] \right) $ and $%
\left[ 0,1\right] $ is compact, $\exists N_{\varepsilon ,m}^{2}$ such that 
\begin{equation*}
\sup_{n>N_{\varepsilon ,m}^{2}}\left\Vert f_{n}-f\right\Vert _{\infty
}<\delta _{m}\text{.}
\end{equation*}%
Fix $m>N_{\varepsilon }^{1}$ and $n>N_{\varepsilon ,m}^{2}$ and consider any 
$t_{i}^{\left( m\right) },t_{i+1}^{\left( m\right) }\in D_{m}$. Then

\begin{eqnarray}
&&\left\vert f_{n}\left( t_{i+1}^{\left( m\right) }\right) -f_{n}\left(
t_{i}^{\left( m\right) }\right) \right\vert  \label{some bounds} \\
&\geq &\left\vert f\left( t_{i+1}^{\left( m\right) }\right) -f\left(
t_{i}^{\left( m\right) }\right) \right\vert -2\delta _{m}.  \notag
\end{eqnarray}%
By the mean value theorem, for any differentiable function $g$ on $\left[
x-\lambda ,x\right] $,%
\begin{equation}
g\left( x\right) \leq g\left( x-\lambda \right) +\lambda \sup_{t\in \left[
x-\lambda ,x\right] }\left\vert g^{\prime }\left( t\right) \right\vert .
\label{mean value theorem}
\end{equation}%
In particular, we choose $g\left( x\right) =x^{p}$ and $\lambda =2\delta
_{m} $ and apply (\ref{mean value theorem}) to the right hand side of (\ref%
{some bounds}), noting that 
\begin{equation*}
\min_{t_{j}^{\left( m\right) },t_{j+1}^{\left( m\right) }\in
D_{m}}\left\vert f\left( t_{j+1}^{\left( m\right) }\right) -f\left(
t_{j}^{\left( m\right) }\right) \right\vert \leq \left\vert f\left(
t_{i+1}^{\left( m\right) }\right) -f\left( t_{i}^{\left( m\right) }\right)
\right\vert \leq 2K.
\end{equation*}%
We then get%
\begin{equation*}
\left( \left\vert f\left( t_{i+1}^{\left( m\right) }\right) -f\left(
t_{i}^{\left( m\right) }\right) \right\vert \right) ^{p}\leq \left\vert
f_{n}\left( t_{i+1}^{\left( m\right) }\right) -f_{n}\left( t_{i}^{\left(
m\right) }\right) \right\vert ^{p}+2^{p}\delta _{m}pK^{p-1}.
\end{equation*}%
Hence, 
\begin{eqnarray*}
c-\varepsilon &\leq &\inf_{m>N_{\varepsilon }^{1}}\mathcal{V}_{p}^{p}\left( 
\hat{f}_{m}\right) \\
&\leq &\inf_{m>N_{\varepsilon }^{1}}\inf_{n>N_{\varepsilon ,m}^{2}}\left( 
\mathcal{V}_{p}^{p}\left( \sum_{D_{m}}f_{n}\left( t_{j}^{\left( m\right)
}\right) \mathbf{1}_{\left[ t_{j}^{\left( m\right) },t_{j+1}^{\left(
m\right) }\right] }\right) \right. \\
&&\left. +2^{p}\delta _{m}pK^{p-1}\left\Vert D_{m}\right\Vert \right) \\
&=&\inf_{m>N_{\varepsilon }^{1}}\inf_{n>N_{\varepsilon ,m}^{2}}\mathcal{V}%
_{p}^{p}\left( \sum_{D_{m}}f_{n}\left( t_{j}^{\left( m\right) }\right) 
\mathbf{1}_{\left[ t_{j}^{\left( m\right) },t_{j+1}^{\left( m\right) }\right]
}\right) +\varepsilon \\
&\leq &\inf_{n>N_{\varepsilon ,m}^{2}}\mathcal{V}_{p}^{p}\left( f_{n}\right)
+\varepsilon .
\end{eqnarray*}%
Letting $\varepsilon \rightarrow 0$, we get%
\begin{equation*}
\mathcal{V}_{p}^{p}\left( f\right) \leq \liminf \mathcal{V}_{p}^{p}\left(
f_{n}\right) .
\end{equation*}%
The result now follows from Lemma \ref{lsc_lemma}.
\end{proof}

Hence,

\begin{proposition}
Fix $p\geq 1$. Let $\mathbf{i}$ denote the imbedding map%
\begin{equation*}
\mathbf{i:}\left( C\left( \left[ 0,1\right] \right) _{p},\left\Vert
.\right\Vert _{\infty }\right) \rightarrow \left( C\left( \left[ 0,1\right]
\right) _{p},\left\Vert .\right\Vert _{\infty }+\mathcal{V}_{p}\left(
.\right) \right)
\end{equation*}%
For every fixed $h\in \left[ 0,1\right] $, 
\begin{equation*}
\mathbf{i\circ }\mathcal{T}_{h}^{B}:\left( C_{0}\left( \left[ -1,1\right]
\right) _{p},\left\Vert .\right\Vert _{\infty }\right) \rightarrow \left(
C\left( \left[ 0,1\right] \right) _{p},\left\Vert .\right\Vert _{p}\right)
\end{equation*}%
is $\sigma _{\left\Vert .\right\Vert _{\infty }}$-measurable. For $p>2$, 
\begin{equation*}
\mathcal{T}_{h}^{B,p}:\left( C_{0}\left( \left[ -1,1\right] \right)
,\left\Vert .\right\Vert _{\infty }\right) \rightarrow \left( C\left( \left[
0,1\right] \right) _{p},\left\Vert .\right\Vert _{p}\right)
\end{equation*}%
as defined in (\ref{p-variation frame process}) is a Borel random variable.
\end{proposition}

\begin{proof}
Proposition \ref{Propn: Semicontinuity of p-variation} proves that the
imbedding map $\mathbf{i}$ is \emph{lower-semicontinuous }and hence $\sigma
_{\left\Vert .\right\Vert _{\infty }}$-measurable\emph{\ }for any $p\geq 1$. 
\emph{\ }In addition, by Lemma \ref{Lemma: sup-norm process is Borel}, for
any fixed $h\in \left[ 0,1\right] $, $\mathcal{T}_{h}:\left( C_{0}\left( %
\left[ -1,1\right] \right) ,\left\Vert .\right\Vert _{\infty }\right)
\rightarrow \left( C\left( \left[ 0,1\right] \right) ,\left\Vert
.\right\Vert _{\infty }\right) $ is continuous, so that the restriction of $%
\mathcal{T}_{h}$ to $\left( C_{0}\left( \left[ -1,1\right] \right)
_{p},\left\Vert .\right\Vert _{\infty }\right) $ is also continuous. Hence, $%
\mathbf{i\circ }\mathcal{T}_{h}$ is $\sigma \left( \left\Vert .\right\Vert
_{\infty }\right) $-measurable.

A result by P. L\'{e}vy (\emph{Theorem 9} in \cite{Levy40}) states that $%
\mathbb{P}$-a.e. Brownian sample path has infinite $p$-variation if $p\leq 2$%
, but finite $p$-variation if $p>2$, so that 
\begin{equation*}
\mathbb{P}\left( C_{0}\left( \left[ -1,1\right] \right) _{p}\right) =1\text{
if }p>2\text{.}
\end{equation*}%
Since $\sigma _{\left\Vert .\right\Vert _{\infty }}$ is complete w.r.t. $%
\mathbb{P}$, both $C_{0}\left( \left[ -1,1\right] \right) _{p}$ and $%
C_{0}\left( \left[ -1,1\right] \right) \setminus C_{0}\left( \left[ -1,1%
\right] \right) _{p}$ are in $\sigma _{\left\Vert .\right\Vert _{\infty }}$.
Hence,%
\begin{equation*}
\mathcal{T}_{h}^{B,p}=\left\{ 
\begin{array}{ll}
\mathcal{T}_{h}^{B} & \text{on }C_{0}\left( \left[ -1,1\right] \right) _{p}
\\ 
0 & \text{otherwise}%
\end{array}%
\right.
\end{equation*}%
is Borel-measurable.
\end{proof}

\section{The H\"{o}lder condition and Corollaries}

The following inequality (\emph{Lemma 3} in \cite{Lyons99}) will be used
frequently:

\begin{lemma}
\label{Lemma: l_p bound}Suppose that $a_{n}\in \mathbb{R}$ for all $n\in 
\mathbb{N}$. Then 
\begin{equation}
\left( \sum_{i=1}^{\infty }\left\vert a_{i}\right\vert \right) ^{p}\leq
c\left( \alpha ,p\right) ^{p}\sum_{i=1}^{\infty }i^{\alpha p}\left\vert
a_{i}\right\vert ^{p}.  \label{Hoelder's inequality}
\end{equation}%
$c\left( \alpha ,p\right) $ is as defined in Section 2.1, equation (\ref%
{c(alpha,p)-definition}), i.e. 
\begin{equation*}
c\left( \alpha ,p\right) =\left( \sum_{n=1}^{\infty }n^{-\alpha \frac{p}{p-1}%
}\right) ^{\frac{p-1}{p}}
\end{equation*}%
where $\alpha >1-\frac{1}{p}$ and $p>1$.
\end{lemma}

\begin{proof}
By H\"{o}lder's inequality%
\begin{eqnarray*}
\sum_{i=1}^{\infty }\left\vert a_{i}\right\vert &=&\sum_{i=1}^{\infty
}i^{-\alpha }i^{\alpha }\left\vert a_{i}\right\vert \\
&\leq &\left( \sum_{i=1}i^{-\alpha \frac{p}{p-1}}\right) ^{\frac{p-1}{p}%
}\left( \sum_{i=1}^{\infty }i^{\alpha p}\left\vert a_{i}\right\vert
^{p}\right) ^{\frac{1}{p}}.
\end{eqnarray*}
\end{proof}

\begin{proposition}
\label{Propn: pathwise p-varn bound}Suppose that\emph{\ }$h\in \left[ 0,1%
\right] $ and that $f\in C\left( \left[ -1,1\right] \right) $. If $p>1$ and $%
\alpha >1-\frac{1}{p}$, then there is a finite constant $c\left( \alpha
,p\right) $ such that%
\begin{eqnarray*}
&&c\left( \alpha ,p\right) \left( \sum_{n=0}^{\infty }\left( n+1\right)
^{\alpha p}\right. \\
&&\left. \sum_{k=0}^{2^{n+n\left( h\right) +1}-1}\left\vert f\left( \frac{k+1%
}{2^{n+n\left( h\right) }}-1\right) -f\left( \frac{k}{2^{n+n\left( h\right) }%
}-1\right) \right\vert ^{p}\right) ^{\frac{1}{p}}
\end{eqnarray*}%
bounds both 
\begin{equation*}
\frac{1}{4}\mathcal{V}_{p}\left( \mathcal{T}_{h}^{f}-\mathcal{T}%
_{0}^{f}\right)
\end{equation*}%
and 
\begin{equation*}
2^{-\frac{p-1}{p}}\left\Vert \mathcal{T}_{h}^{f}-\mathcal{T}%
_{0}^{f}\right\Vert _{\infty }
\end{equation*}%
and hence 
\begin{equation*}
\left( 4+2^{\frac{p-1}{p}}\right) ^{-1}\left\Vert \mathcal{T}_{h}^{f}-%
\mathcal{T}_{0}^{f}\right\Vert _{p}.
\end{equation*}%
Here, $c\left( .,.\right) $ is defined in Section 2.1, equation (\ref%
{c(alpha,p)-definition}).
\end{proposition}

\begin{proof}
If $x\in \mathbb{R}$, let $\left\lceil x\right\rceil $ ($\left\lfloor
x\right\rfloor $) denote the smallest (largest) integer that is greater
(less)\ than $x$. Any interval $\left[ s,t\right) \subset \left[ 0,2\right] $
with $t-s\leq 1$ is a countable union of disjoint dyadic intervals contained
in $\left[ s,t\right) $ of the form%
\begin{equation*}
\left[ \text{\b{t}}_{j-1},\text{\b{t}}_{j}\right) :=\left[ \left\lfloor
2^{j-1}t\right\rfloor 2^{-\left( j-1\right) },\left\lfloor
2^{j}t\right\rfloor 2^{-j}\right)
\end{equation*}%
and 
\begin{equation*}
\left[ \bar{s}_{j},\bar{s}_{j-1}\right) :=\left[ \left\lceil
2^{j}s\right\rceil 2^{-j},\left\lceil 2^{j-1}s\right\rceil 2^{-\left(
j-1\right) }\right) ,
\end{equation*}%
where $j>n\left( t-s\right) $ ($n\left( h\right) $ is defined in \emph{%
Section 2.1}). Suppose that $f\in C\left( \left[ -1,1\right] \right) $. In
order to simplify notation, we \textquotedblleft shift\textquotedblright\ $f$
to the right by $1$ and work with 
\begin{equation*}
\hat{f}\equiv f\left( .-1\right) \in C\left( \left[ 0,2\right] \right) .
\end{equation*}%
With $\hat{f}$ we associate the function $\hat{f}^{1}\left( s,t\right) :=%
\hat{f}\left( t\right) -\hat{f}\left( s\right) $ $\in C\left( \left[ 0,2%
\right] \times \left[ 0,2\right] \right) $. Since $\bar{s}_{j}\searrow s$
and \b{t}$_{j}\nearrow t$ and $\hat{f}$ $^{1}\ $is continuous, we have the
following representation for $\hat{f}^{1}\left( s,t\right) $: 
\begin{equation}
\hat{f}^{1}\left( s,t\right) =\sum_{j=n\left( t-s\right) }^{\infty }\left( 
\hat{f}^{1}\left( \bar{s}_{j+1},\bar{s}_{j}\right) +\hat{f}^{1}\left( \text{%
\b{t}}_{j},\text{\b{t}}_{j+1}\right) \right) .  \label{representation}
\end{equation}

Using (\ref{representation}) we may now apply (\ref{Hoelder's inequality})
from Lemma \ref{Lemma: l_p bound} to bound the $p^{\text{th}}$ power of $%
\hat{f}^{1}\left( s,t\right) $ by the $p^{\text{th}}$ powers of the
increments over dyadics contained in $\left[ s,t\right) :$ 
\begin{eqnarray}
&&\left\vert \hat{f}^{1}\left( s,t\right) \right\vert ^{p}  \notag \\
&\leq &c\left( \alpha ,p\right) ^{p}\sum_{i=0}^{\infty }\left( i+1\right)
^{\alpha p}\left\vert \hat{f}^{1}\left( \bar{s}_{i+1+n\left( t-s\right) },%
\bar{s}_{i+n\left( t-s\right) }\right) +\hat{f}^{1}\left( \text{\b{t}}%
_{i+n\left( t-s\right) },\text{\b{t}}_{i+1+n\left( t-s\right) }\right)
\right\vert ^{p},  \notag
\end{eqnarray}%
which by Jensen's inequality is bounded by

\begin{eqnarray}
&&2^{p-1}c\left( \alpha ,p\right) ^{p}\sum_{i=0}^{\infty }\left( i+1\right)
^{\alpha p}\left( \left\vert \hat{f}^{1}\left( \bar{s}_{i+1+n\left(
t-s\right) },\bar{s}_{i+n\left( t-s\right) }\right) \right\vert ^{p}\right. 
\notag \\
&&\left. +\left\vert \hat{f}^{1}\left( \text{\b{t}}_{i+n\left( t-s\right) },%
\text{\b{t}}_{i+1+n\left( t-s\right) }\right) \right\vert ^{p}\right)
\label{1st p-bound}
\end{eqnarray}%
But 
\begin{equation*}
\left\{ \bar{s}_{j},\text{\b{t}}_{j}:j\geq n\left( t-s\right) \right\}
\end{equation*}%
is contained in 
\begin{equation*}
\left[ s,t\right] \cap \left\{ \frac{k}{2^{n}}:0\leq k\leq 2^{n+1};n\geq
n\left( t-s\right) \right\} .
\end{equation*}%
Hence, (\ref{1st p-bound}) and in turn $\left\vert \hat{f}^{1}\left(
s,t\right) \right\vert ^{p}$ is bounded by 
\begin{equation}
2^{p-1}c\left( \alpha ,p\right) ^{p}\sum_{i=0}^{\infty }\left( i+1\right)
^{\alpha p}\sum_{k:s\leq \frac{k}{2^{n\left( t-s\right) +i}}<\frac{k+1}{%
2^{n\left( t-s\right) +i}}\leq t}\left\vert \hat{f}^{1}\left( \frac{k}{%
2^{n\left( t-s\right) +i}},\frac{k+1}{2^{n\left( t-s\right) +i}}\right)
\right\vert ^{p}.  \label{p-bound for continuous functions}
\end{equation}%
We now use (\ref{p-bound for continuous functions}) to bound

\begin{enumerate}
\item $\sup_{t\in \left[ 0,1\right] }\left\vert \hat{f}^{1}\left(
t,t+h\right) \right\vert $ where $h\in \left[ 0,1\right] .$

For each fixed $t$, (\ref{p-bound for continuous functions}) applies. Note
that (\ref{p-bound for continuous functions}) is itself bounded by 
\begin{equation*}
2^{p-1}c\left( \alpha ,p\right) ^{p}\sum_{i=0}^{\infty }\left( i+1\right)
^{\alpha p}\sum_{k=0}^{2^{i+n\left( h\right) +1}-1}\left\vert \hat{f}%
^{1}\left( \frac{k}{2^{n\left( h\right) +i}},\frac{k+1}{2^{n\left( h\right)
+i}}\right) \right\vert ^{p},
\end{equation*}%
which only depends on $h$ -- not on $\ t$ -- so that 
\begin{eqnarray}
&&\sup_{t\in \left[ 0,1\right] }\left\vert \hat{f}^{1}\left( t,t+h\right)
\right\vert  \notag \\
&\leq &2^{\frac{p-1}{p}}c\left( \alpha ,p\right) \left[ \sum_{i=0}^{\infty
}\left( i+1\right) ^{\alpha p}\sum_{k=0}^{2^{i+n\left( h\right)
+1}-1}\left\vert \hat{f}^{1}\left( \frac{k}{2^{n\left( h\right) +i}},\frac{%
k+1}{2^{n\left( h\right) +i}}\right) \right\vert ^{p}\right] ^{\frac{1}{p}}.
\label{sup-bound}
\end{eqnarray}

\item $\mathcal{V}_{p}\left( \hat{f}\left( .+h\right) -\hat{f}\right) $
where $h\in \left[ 0,1\right] $.

Let $D$ be a dissection of $\left[ 0,1\right] $. We may then rewrite $D$ as 
\begin{equation*}
D=D^{\geq h}\cup D^{<h},
\end{equation*}%
where%
\begin{equation*}
D^{\geq h}:=\left\{ t_{i}\in D:t_{i+1}-t_{i}\geq h\right\}
\end{equation*}%
and%
\begin{equation*}
D^{<h}:=\left\{ t_{i}\in D:t_{i+1}-t_{i}<h\right\} .
\end{equation*}%
Since 
\begin{equation*}
\hat{f}^{1}\left( t_{i}+h,t_{i+1}+h\right) -\hat{f}^{1}\left(
t_{i},t_{i+1}\right) =\hat{f}^{1}\left( t_{i+1},t_{i+1}+h\right) -\hat{f}%
^{1}\left( t_{i},t_{i}+h\right) ,
\end{equation*}%
and by Jensen's inequality, we have 
\begin{eqnarray*}
&&\sum_{D}\left\vert \hat{f}^{1}\left( t_{i}+h,t_{i+1}+h\right) -\hat{f}%
^{1}\left( t_{i},t_{i+1}\right) \right\vert ^{p} \\
&=&\sum_{D^{<h}}\left\vert \hat{f}^{1}\left( t_{i}+h,t_{i+1}+h\right) -\hat{f%
}^{1}\left( t_{i},t_{i+1}\right) \right\vert ^{p} \\
&&+\sum_{D^{\geq h}}\left\vert \hat{f}^{1}\left( t_{i+1},t_{i+1}+h\right) -%
\hat{f}^{1}\left( t_{i},t_{i}+h\right) \right\vert ^{p} \\
&\leq &2^{p-1}\underset{\text{Sum }A}{\underbrace{\sum_{D^{<h}}\left\vert 
\hat{f}^{1}\left( t_{i},t_{i+1}\right) \right\vert ^{p}}}+2^{p-1}\underset{%
\text{Sum }B}{\underbrace{\sum_{D^{<h}}\left\vert \hat{f}^{1}\left(
t_{i}+h,t_{i+1}+h\right) \right\vert ^{p}}} \\
&&+2^{p-1}\underset{\text{Sum }C}{\underbrace{\sum_{D^{\geq h}}\left\vert 
\hat{f}^{1}\left( t_{i},t_{i}+h\right) \right\vert ^{p}}}+2^{p-1}\underset{%
\text{Sum }D}{\underbrace{\sum_{D^{\geq h}}\left\vert \hat{f}^{1}\left(
t_{i+1},t_{i+1}+h\right) \right\vert ^{p}}}.
\end{eqnarray*}%
Note that 
\begin{equation*}
\left\{ \left[ t_{i},t_{i+1}\right) :t_{i},t_{i+1}\in D^{<h}\right\}
\end{equation*}%
contains pairwise disjoint intervals as do 
\begin{eqnarray*}
&&\left\{ \left[ t_{i}+h,t_{i+1}+h\right) :t_{i},t_{i+1}\in D^{<h}\right\} ,
\\
&&\left\{ \left[ t_{i},t_{i}+h\right) :t_{i}\in D^{\geq h}\right\}
\end{eqnarray*}%
and 
\begin{equation*}
\left\{ \left[ t_{i+1},t_{i+1}+h\right) :t_{i+1}\in D^{\geq h}\right\} .
\end{equation*}%
Hence, applying (\ref{p-bound for continuous functions}) to every term in
Sum $A$ and then summing, there is no double counting of increments over
dyadics so that%
\begin{equation*}
\text{Sum }A\leq 2^{p-1}c\left( \alpha ,p\right) ^{p}\sum_{i=0}^{\infty
}\left( i+1\right) ^{\alpha p}\sum_{k=0}^{2^{i+n\left( h\right)
+1}-1}\left\vert \hat{f}^{1}\left( \frac{k}{2^{n\left( h\right) +i}},\frac{%
k+1}{2^{n\left( h\right) +i}}\right) \right\vert ^{p}\text{.}
\end{equation*}%
The same argument applies to Sum $B$, Sum $C$ and Sum $D$. Hence,%
\begin{eqnarray*}
&&\sum_{D}\left\vert \hat{f}^{1}\left( t_{i}+h,t_{i+1}+h\right) -\hat{f}%
\left( t_{i},t_{i+1}\right) \right\vert ^{p} \\
&\leq &2^{p-1}\left( \text{Sum }A+\text{Sum }B+\text{Sum }C+\text{Sum }%
D\right) \\
&\leq &2^{2p}c\left( \alpha ,p\right) ^{p}\sum_{i=0}^{\infty }\left(
i+1\right) ^{\alpha p}\sum_{k=0}^{2^{i+n\left( h\right) +1}-1}\left\vert 
\hat{f}^{1}\left( \frac{k}{2^{n\left( h\right) +i}},\frac{k+1}{2^{n\left(
h\right) +i}}\right) \right\vert ^{p}\text{.}
\end{eqnarray*}%
But this bound does not depend on the particular dissection $D$ chosen, so
that 
\begin{equation*}
\mathcal{V}_{p}\left( \hat{f}\left( .+h\right) -\hat{f}\right) ^{p}\leq
2^{2p}c\left( \alpha ,p\right) ^{p}\sum_{i=0}^{\infty }\left( i+1\right)
^{\alpha p}\sum_{k=0}^{2^{i+n\left( h\right) +1}-1}\left\vert \hat{f}%
^{1}\left( \frac{k}{2^{n\left( h\right) +i}},\frac{k+1}{2^{n\left( h\right)
+i}}\right) \right\vert ^{p}\text{.}
\end{equation*}
\end{enumerate}

If we now \textquotedblleft shift\textquotedblright\ $\hat{f}\in C\left( %
\left[ 0,2\right] \right) $ back into $C\left( \left[ -1,1\right] \right) $
by replacing $t$ by $t-1,$ we retrieve $f$, and so taken together, $1.$ and $%
2.$ prove the Proposition.
\end{proof}

\begin{corollary}
\label{Cor: Corollary to Propn: pathwise p-varn bound}Suppose that $0\leq
h_{1}\leq h_{2}\leq 1$ and that $f\in C\left( \left[ -1,1\right] \right) $.
If $p>1$ and $\alpha >1-\frac{1}{p}$, then 
\begin{eqnarray*}
&&c\left( \alpha ,p\right) \left( \sum_{n=0}^{\infty }\left( n+1\right)
^{\alpha p}\right. \\
&&\left. \sum_{k=0}^{2^{n+n\left( h_{2}-h_{1}\right) +1}-1}\left\vert
f\left( \frac{k+1}{2^{n+n\left( h_{2}-h_{1}\right) }}-1\right) -f\left( 
\frac{k}{2^{n+n\left( h_{2}-h_{1}\right) }}-1\right) \right\vert ^{p}\right)
^{\frac{1}{p}}.
\end{eqnarray*}%
bounds both 
\begin{equation*}
\frac{1}{4}\mathcal{V}_{p}\left( \mathcal{T}_{h_{2}}^{f}-\mathcal{T}%
_{h_{1}}^{f}\right)
\end{equation*}%
and 
\begin{equation*}
2^{-\frac{p-1}{p}}\left\Vert \mathcal{T}_{h_{2}}^{f}-\mathcal{T}%
_{h_{1}}^{f}\right\Vert _{\infty }
\end{equation*}%
and hence 
\begin{equation*}
\left( 4+2^{\frac{p-1}{p}}\right) ^{-1}\left\Vert \mathcal{T}_{h_{2}}^{f}-%
\mathcal{T}_{h_{1}}^{f}\right\Vert _{p}.
\end{equation*}%
$c\left( .,.\right) $ is defined in Section 2.1, equation (\ref%
{c(alpha,p)-definition}).
\end{corollary}

\begin{proof}
Let $u=t+h_{1}$ and $h=h_{2}-h_{1}$, so that 
\begin{equation*}
\left( \mathcal{T}_{h_{2}}^{f}-\mathcal{T}_{h_{1}}^{f}\right) \left(
t\right) =\left( \mathcal{T}_{h}^{f}-\mathcal{T}_{0}^{f}\right) \left(
u\right) .
\end{equation*}%
Then apply Proposition \ref{Propn: pathwise p-varn bound}.
\end{proof}

\subsection{\protect\bigskip The $p$-variation-norm frame process $\mathcal{S%
}^{B,p}$}

We now prove the following H\"{o}lder-type inequality:

\begin{proposition}
\label{Hoelder's condition}Suppose that $0\leq h_{1}\leq h_{2}\leq 1$. If $%
\acute{p}>p>2$, $\alpha >1-\frac{1}{p}$ and $\beta >1-\frac{1}{\acute{p}}$,
then there is a finite constant $d\left( \alpha ,\beta ,p,\acute{p}\right) $
such that 
\begin{equation*}
\mathbb{E}\left[ \left\Vert \mathcal{T}_{h_{2}}^{B,p}-\mathcal{T}%
_{h_{1}}^{B,p}\right\Vert _{p}^{\acute{p}}\right] \leq d\left( \alpha ,\beta
,p,\acute{p}\right) \left( h_{2}-h_{1}\right) ^{\left( \frac{1}{2}-\frac{1}{p%
}\right) \acute{p}}.
\end{equation*}%
Here, 
\begin{equation*}
d\left( \alpha ,\beta ,p,\acute{p}\right) :=2^{\frac{\acute{p}}{2}}\left[
2\left( 2+2^{-\frac{1}{p}}\right) c\left( \alpha ,p\right) \right] ^{\acute{p%
}}c\left( \beta ,\frac{\acute{p}}{p}\right) ^{\frac{\acute{p}}{p}}\sqrt{%
\frac{2^{\acute{p}}}{\pi }}\Gamma \left( \frac{\acute{p}+1}{2}\right)
\sum_{n=0}^{\infty }\left( n+1\right) ^{\alpha \acute{p}+\beta \frac{\acute{p%
}}{p}}2^{n\left( \frac{1}{p}-\frac{1}{2}\right) \acute{p}}.
\end{equation*}%
If $\acute{p}=p>2$, 
\begin{equation*}
\mathbb{E}\left[ \left\Vert \mathcal{T}_{h_{2}}^{B,p}-\mathcal{T}%
_{h_{1}}^{B,p}\right\Vert _{p}^{p}\right] \leq d\left( \alpha ,p\right)
\left( h_{2}-h_{1}\right) ^{\left( \frac{1}{2}-\frac{1}{p}\right) p}.
\end{equation*}%
If $p>\acute{p}>2,$ 
\begin{equation*}
\mathbb{E}\left[ \left\Vert \mathcal{T}_{h_{2}}^{B,p}-\mathcal{T}%
_{h_{1}}^{B,p}\right\Vert _{p}^{\acute{p}}\right] \leq d\left( \alpha
,p\right) ^{\frac{\acute{p}}{p}}\left( h_{2}-h_{1}\right) ^{\left( \frac{1}{2%
}-\frac{1}{p}\right) \acute{p}}
\end{equation*}%
Here,%
\begin{equation*}
d\left( \alpha ,p\right) :=2^{\frac{p}{2}}\left( 4+2^{\frac{p-1}{p}}\right)
^{p}c\left( \alpha ,p\right) ^{p}\sqrt{\frac{2^{p}}{\pi }}\Gamma \left( 
\frac{p+1}{2}\right) \sum_{n=1}^{\infty }\left( n+1\right) ^{\alpha
p}2^{n\left( 1-\frac{p}{2}\right) }.
\end{equation*}%
$c\left( .,.\right) $ is defined in Section 2.1, equation (\ref%
{c(alpha,p)-definition}).
\end{proposition}

\begin{proof}
We are interested in the $\acute{p}^{\text{th }}$moment of 
\begin{equation*}
\left\Vert \mathcal{T}_{h_{2}}^{B,p}-\mathcal{T}_{h_{1}}^{B,p}\right\Vert
_{p}=\mathcal{V}_{p}\left( \mathcal{T}_{h_{2}}^{B,p}-\mathcal{T}%
_{h_{1}}^{B,p}\right) +\left\Vert \mathcal{T}_{h_{2}}^{B,p}-\mathcal{T}%
_{h_{1}}^{B,p}\right\Vert _{\infty }.
\end{equation*}%
Since $p>2$, \ the set $C_{0}\left( \left[ -1,1\right] \right) _{p}$ has
full $\mathbb{P}$-measure (Theorem 9 in \cite{Levy40}) , so that we may work
on $C_{0}\left( \left[ -1,1\right] \right) _{p}$ instead of $C_{0}\left( %
\left[ -1,1\right] \right) $. We fix $f\in C_{0}\left( \left[ -1,1\right]
\right) _{p}$ and consider the sample path $\mathcal{T}^{f}$ of $\mathcal{T}$
: For any $h\in \left[ 0,1\right] $, 
\begin{equation*}
\mathcal{T}_{h}^{f}\in \left( C\left( \left[ 0,1\right] \right)
_{p},\left\Vert .\right\Vert _{p}\right) .
\end{equation*}%
Let $h:=h_{2}-h_{1}$. We first consider the case where\ $\acute{p}>p>2$. By
Corollary \ref{Cor: Corollary to Propn: pathwise p-varn bound}, 
\begin{eqnarray*}
&&\left\Vert \mathcal{T}_{h_{2}}^{f}-\mathcal{T}_{h_{1}}^{f}\right\Vert
_{p}^{\acute{p}} \\
&\leq &\left[ 2\left( 2+2^{-\frac{1}{p}}\right) c\left( \alpha ,p\right) %
\right] ^{\acute{p}} \\
&&\left( \sum_{n=0}^{\infty }\underset{:=b_{n}}{\underbrace{\left(
n+1\right) ^{\alpha p}\sum_{k=0}^{2^{n+n\left( h\right) +1}-1}\left\vert
f\left( \frac{k+1}{2^{n+n\left( h\right) }}-1\right) -f\left( \frac{k}{%
2^{n+n\left( h\right) }}-1\right) \right\vert ^{p}}}\right) ^{\frac{\acute{p}%
}{p}}.
\end{eqnarray*}%
By Lemma \ref{Lemma: l_p bound},%
\begin{equation*}
\left( \sum_{n=0}^{\infty }b_{n}\right) ^{\frac{\acute{p}}{p}}\leq c\left(
\beta ,\frac{\acute{p}}{p}\right) ^{\frac{\acute{p}}{p}}\sum_{n=0}^{\infty
}\left( n+1\right) ^{\beta \frac{\acute{p}}{p}}b_{n}^{\frac{\acute{p}}{p}}.
\end{equation*}%
Furthermore, by Jensen's inequality,%
\begin{eqnarray*}
b_{n}^{\frac{\acute{p}}{p}} &\leq &\left( n+1\right) ^{\alpha \acute{p}%
}\left( 2^{n+n\left( h\right) +1}\right) ^{\frac{\acute{p}}{p}-1} \\
&&\sum_{k=0}^{2^{n+n\left( h\right) }-1}\left\vert f\left( \frac{k+1}{%
2^{n+n\left( h\right) }}-1\right) -f\left( \frac{k}{2^{n+n\left( h\right) }}%
-1\right) \right\vert ^{\acute{p}}.
\end{eqnarray*}%
Hence,%
\begin{eqnarray*}
&&\left\Vert \mathcal{T}_{h_{2}}^{f}-\mathcal{T}_{h_{1}}^{f}\right\Vert
_{p}^{\acute{p}} \\
&\leq &\left[ 2\left( 2+2^{-\frac{1}{p}}\right) c\left( \alpha ,p\right) %
\right] ^{\acute{p}}c\left( \beta ,\frac{\acute{p}}{p}\right) ^{\frac{\acute{%
p}}{p}} \\
&&\sum_{n=0}^{\infty }\left( n+1\right) ^{\alpha \acute{p}+\beta \frac{%
\acute{p}}{p}}\left( 2^{n+n\left( h\right) +1}\right) ^{\frac{\acute{p}}{p}%
-1} \\
&&\sum_{k=0}^{2^{n+n\left( h\right) +1}-1}\left\vert f\left( \frac{k+1}{%
2^{n+n\left( h\right) }}-1\right) -f\left( \frac{k}{2^{n+n\left( h\right) }}%
-1\right) \right\vert ^{\acute{p}}.
\end{eqnarray*}

In order to calculate the moment bound, we recall that 
\begin{equation*}
\sqrt{\frac{2}{\pi \sigma ^{2}}}\int_{0}^{\infty }x^{\acute{p}}\exp \left( -%
\frac{x^{2}}{2\sigma ^{2}}\right) dx=\sqrt{\frac{2^{\acute{p}}}{\pi }}\sigma
^{p}\underset{:=\Gamma \left( \frac{\acute{p}+1}{2}\right) }{\underbrace{%
\int_{0}^{\infty }x^{\frac{\acute{p}-1}{2}}e^{-x}dx}}.
\end{equation*}%
Hence, we get%
\begin{equation}
\mathbb{E}\left[ \left\vert B_{t}-B_{s}\right\vert ^{\acute{p}}\right] =%
\sqrt{\frac{2^{\acute{p}}}{\pi }}\Gamma \left( \frac{\acute{p}+1}{2}\right)
\left\vert t-s\right\vert ^{\frac{\acute{p}}{2}}.  \label{Gaussian moments}
\end{equation}%
So by (\ref{Gaussian moments}), 
\begin{eqnarray*}
&&\mathbb{E}\left[ \left( \left\Vert \mathcal{T}_{h_{2}}-\mathcal{T}%
_{h_{1}}\right\Vert _{p}^{p}\right) ^{\frac{\acute{p}}{p}}\right] \\
&\leq &2^{\frac{\acute{p}}{2}}\left[ 2\left( 2+2^{-\frac{1}{p}}\right)
c\left( \alpha ,p\right) \right] ^{\acute{p}}c\left( \beta ,\frac{\acute{p}}{%
p}\right) ^{\frac{\acute{p}}{p}}\sqrt{\frac{2^{\acute{p}}}{\pi }}\Gamma
\left( \frac{\acute{p}+1}{2}\right) \\
&&\sum_{n=0}^{\infty }\left( n+1\right) ^{\alpha \acute{p}+\beta \frac{%
\acute{p}}{p}}\left( 2^{n+n\left( h\right) +1}\right) ^{\frac{\acute{p}}{p}-%
\frac{\acute{p}}{2}} \\
&\leq &d\left( \alpha ,\beta ,p,\acute{p}\right) h^{\left( \frac{1}{2}-\frac{%
1}{p}\right) \acute{p}}
\end{eqnarray*}%
where%
\begin{eqnarray*}
&&d\left( \alpha ,\beta ,p,\acute{p}\right) \\
&=&2^{\frac{\acute{p}}{2}}\left[ 2\left( 2+2^{-\frac{1}{p}}\right) c\left(
\alpha ,p\right) \right] ^{\acute{p}}c\left( \beta ,\frac{\acute{p}}{p}%
\right) ^{\frac{\acute{p}}{p}}\sqrt{\frac{2^{\acute{p}}}{\pi }}\Gamma \left( 
\frac{\acute{p}+1}{2}\right) \\
&&\times \sum_{n=0}^{\infty }\left( n+1\right) ^{\alpha \acute{p}+\beta 
\frac{\acute{p}}{p}}2^{n\left( \frac{\acute{p}}{p}-\frac{\acute{p}}{2}%
\right) }
\end{eqnarray*}%
which is finite if $p>2$, $\alpha >1-\frac{1}{p}$ and $\beta >1-\frac{p}{%
\acute{p}}.$

Next, we consider the case where $\acute{p}=p>2$: \ By Corollary \ref{Cor:
Corollary to Propn: pathwise p-varn bound},%
\begin{eqnarray*}
&&\left\Vert \mathcal{T}_{h_{2}}^{f}-\mathcal{T}_{h_{1}}^{f}\right\Vert
_{p}^{p} \\
&\leq &\left( 4+2^{\frac{p-1}{p}}\right) ^{p}c\left( \alpha ,p\right)
^{p}\sum_{n=0}^{\infty }\left( n+1\right) ^{\alpha p} \\
&&\sum_{k=0}^{2^{n+n\left( h\right) +1}-1}\left\vert f\left( \frac{k+1}{%
2^{n+n\left( h_{2}-h_{1}\right) }}-1\right) -f\left( \frac{k}{2^{n+n\left(
h_{2}-h_{1}\right) }}-1\right) \right\vert ^{p}.
\end{eqnarray*}%
Using (\ref{Gaussian moments}) we find that 
\begin{eqnarray*}
&&\mathbb{E}\left[ \left\Vert \mathcal{T}_{h_{2}}^{B,p}-\mathcal{T}%
_{h_{1}}^{B,p}\right\Vert _{p}^{p}\right] \\
&\leq &2^{\frac{p}{2}}\left( 4+2^{\frac{p-1}{p}}\right) ^{p}c\left( \alpha
,p\right) ^{p}\sqrt{\frac{2^{p}}{\pi }}\Gamma \left( \frac{p+1}{2}\right)
\sum_{n=0}^{\infty }\left( n+1\right) ^{\alpha p}2^{\left( n+n\left(
h\right) +1\right) \left( 1-\frac{p}{2}\right) } \\
&\leq &\underset{=d\left( \alpha ,p\right) }{\underbrace{\left( 2^{\frac{p}{2%
}}\left( 4+2^{\frac{p-1}{p}}\right) ^{p}c\left( \alpha ,p\right) ^{p}\sqrt{%
\frac{2^{p}}{\pi }}\Gamma \left( \frac{p+1}{2}\right) \right)
\sum_{n=0}^{\infty }\left( n+1\right) ^{\alpha p}2^{n\left( 1-\frac{p}{2}%
\right) }}}h^{\frac{p}{2}-1}.
\end{eqnarray*}%
Finally, if $2<\acute{p}<p$, by Lyapunov's inequality we have that 
\begin{eqnarray*}
&&\mathbb{E}\left[ \left\Vert \mathcal{T}_{h_{2}}^{B,p}-\mathcal{T}%
_{h_{1}}^{B,p}\right\Vert _{p}^{\acute{p}}\right] \\
&=&\mathbb{E}\left[ \left\Vert \mathcal{T}_{h_{2}}^{B,p}-\mathcal{T}%
_{h_{1}}^{B,p}\right\Vert _{p}^{p\frac{\acute{p}}{p}}\right] \\
&\leq &\mathbb{E}\left[ \left\Vert \mathcal{T}_{h_{2}}^{B,p}-\mathcal{T}%
_{h_{1}}^{B,p}\right\Vert _{p}^{p}\right] ^{\frac{\acute{p}}{p}}\leq d\left(
\alpha ,p\right) ^{\frac{\acute{p}}{p}}\left( h_{2}-h_{2}\right) ^{\left( 
\frac{1}{2}-\frac{1}{p}\right) \acute{p}}.
\end{eqnarray*}
\end{proof}

We recall \emph{Kolmogorov's lemma }(e.g. \cite{Mckean69}):

\begin{theorem}[Kolmogorov's Lemma ]
\label{Thm: Kolmogorov's Lemma}Let $X_{t}$, $t\in \left[ 0,1\right] ^{d}$, \
be a Banach space $V$-valued process for which there exist three strictly
positive constants $\gamma $,$c$,$\varepsilon $ such that 
\begin{equation*}
\mathbb{E}\left[ \left\Vert X_{t}-X_{s}\right\Vert _{V}^{\gamma }\right]
\leq c\left\vert t-s\right\vert ^{d+\varepsilon };
\end{equation*}%
then the process%
\begin{equation*}
X_{t}^{\ast }\equiv \lim \inf_{u\in D\left( \left[ 0,1\right] ^{d}\right)
:u\geq t}X_{u}
\end{equation*}%
is a continuous modification\footnote{$X_{.}$ is called a modification of $%
Y_{.}$ if for any $t\in \left[ 0,1\right] ^{d}$ we have $\mathbb{P}\left(
X_{t}=Y_{t}\right) =1$.} of $X_{.}$.
\end{theorem}

\begin{corollary}
\label{continuous version}For $p>2$, the process $\mathcal{T}_{.}^{B,p}:%
\left[ 0,1\right] \rightarrow \left( C\left( \left[ 0,1\right] \right)
_{p},\left\Vert .\right\Vert _{p}\right) $ has a continuous modification --
the Brownian frame process whose evaluation at time $t\in \left[ 0,1\right] $
is 
\begin{equation}
\mathcal{S}_{t}^{B,p}:=\lim \inf_{u\in D\left( \left[ 0,1\right] \right)
:u\geq t}\mathcal{T}_{u}^{B,p}.  \label{frame process definition}
\end{equation}
\end{corollary}

\begin{proof}
Applying Proposition \ref{Hoelder's condition} with $\acute{p}>\frac{2p}{p-2}
$ and 
\begin{equation*}
C=\left\{ 
\begin{array}{ll}
d\left( \alpha ,\beta ,p,\acute{p}\right) & \text{if }\acute{p}>p>2 \\ 
d\left( \alpha ,p\right) & \text{if }\acute{p}=p>2 \\ 
d\left( \alpha ,p\right) ^{\frac{\acute{p}}{p}} & \text{if }p>\acute{p}>2%
\end{array}%
\right. ,
\end{equation*}%
we see that $\mathcal{T}_{.}^{B,p}$ satisfies the conditions of Kolmogorov's
lemma so that (\ref{frame process definition}) gives a continuous
modification of $\mathcal{T}_{.}^{B,p}$ .
\end{proof}

\begin{definition}[dyadic polygonal approximation]
\label{Defn: dyadic polygonal approximation}The dyadic polygonal
approximations $\left( X_{.}\left( m\right) \right) _{m\in \mathbb{N}}$ to a
Banach space $V$ valued path $X_{.}$ are defined as 
\begin{equation}
X\left( m\right) _{t}:=X_{\frac{k-1}{2^{m}}}+2^{m}\left( t-\frac{k-1}{2^{m}}%
\right) \Delta _{k}^{m}X_{.}\text{ \ \ \ if }\frac{k-1}{2^{m}}\leq t\leq 
\frac{k}{2^{m}}.  \label{dyadic polygonal approximations}
\end{equation}%
Here, 
\begin{equation*}
\Delta _{k}^{m}X:=X_{\frac{k}{2^{m}}}-X_{\frac{k-1}{2^{m}}}.
\end{equation*}
\end{definition}

\begin{proposition}[Proposition $4.3.1$ and $4.3.2$ in \protect\cite{Lyons02}%
]
\label{Lyons' Holder Proposition}Suppose $\left( X_{t}\right) $ is a
continuous Banach space $V$ valued stochastic process on a completed
probability space $\left( \Omega ,\mathcal{F},\mathbb{P}\right) $ for which
there are constants $\acute{p}>1$, $\kappa \in \left( 0,1\right) $ such that 
$\kappa \acute{p}>1$ as well as a third constant $C$, \ such that 
\begin{equation}
\mathbb{E}\left[ \left\vert X_{t}-X_{s}\right\vert _{V}^{\acute{p}}\right]
\leq C\left\vert t-s\right\vert ^{\kappa \acute{p}}\text{ }\forall s,t\in %
\left[ 0,1\right] .  \label{Lyons' Holder condition}
\end{equation}%
Then the dyadic polygonal approximations $X\left( m\right) $ have finite $%
\acute{p}$-variation uniformly in $m$, $\mathbb{P}$-a.s. Furthermore, $%
\left( X\left( m\right) \right) $ converges to $X$ in $\acute{p}$-variation $%
\mathbb{P}$-a.s.
\end{proposition}

Hence,

\begin{corollary}
For $p,\acute{p}>2$ the $p$-variation norm frame process%
\begin{equation*}
\mathcal{S}^{B,p}:\left[ 0,1\right] \rightarrow \left( C\left( \left[ 0,1%
\right] \right) _{p},\left\Vert .\right\Vert _{p}\right)
\end{equation*}%
\bigskip has finite $\acute{p}$-variation for $\acute{p}>\frac{2p}{p-2},$ $%
\mathbb{P}$-a.s. Furthermore, the dyadic polygonal approximations $\mathcal{S%
}^{B,p}\left( m\right) $ as defined in (\ref{dyadic polygonal approximations}%
) converge to $\mathcal{S}^{B,p}$ in $\acute{p}$-variation $\mathbb{P}$-a.s.
\end{corollary}

\begin{proof}
By Proposition \ref{Hoelder's condition}, for appropriate choices of $\alpha 
$ and $\beta $ and $0\leq h_{1}\leq h_{2}\leq 1$, we have that 
\begin{eqnarray*}
&&\mathbb{E}\left[ \left\Vert \mathcal{S}_{h_{2}}^{B,p}-\mathcal{S}%
_{h_{1}}^{B,p}\right\Vert _{p}^{\acute{p}}\right] \\
&\leq &\left\{ 
\begin{array}{ll}
d\left( \alpha ,\beta ,p,\acute{p}\right) \left( h_{2}-h_{1}\right) ^{\left( 
\frac{1}{2}-\frac{1}{p}\right) \acute{p}} & \text{if }\acute{p}>p>2 \\ 
d\left( \alpha ,p\right) \left( h_{2}-h_{1}\right) ^{\left( \frac{1}{2}-%
\frac{1}{p}\right) p} & \text{if }\acute{p}=p>2 \\ 
d\left( \alpha ,p\right) ^{\frac{\acute{p}}{p}}\left( h_{2}-h_{1}\right)
^{\left( \frac{1}{2}-\frac{1}{p}\right) p} & \text{if }p>\acute{p}>2%
\end{array}%
\right. .
\end{eqnarray*}%
Since $\mathcal{S}_{.}^{B,p}$ is by definition continuous and takes its
values in the separable Banach space $\left( C\left( \left[ 0,1\right]
\right) _{p},\left\Vert .\right\Vert _{p}\right) $, by choosing $\acute{p}>%
\frac{2p}{p-2}$, we may apply Proposition \ref{Lyons' Holder Proposition} to 
$\mathcal{S}^{B,p}$ with $\kappa =\frac{1}{2}-\frac{1}{p}$ and 
\begin{equation*}
C=\left\{ 
\begin{array}{ll}
d\left( \alpha ,\beta ,p,\acute{p}\right) & \text{if }\acute{p}>p>2 \\ 
d\left( \alpha ,p\right) & \text{if }\acute{p}=p>2 \\ 
d\left( \alpha ,p\right) ^{\frac{\acute{p}}{p}} & \text{if }p>\acute{p}>2%
\end{array}%
\right. .
\end{equation*}
\end{proof}

\subsection{\protect\bigskip The $\sup $-norm frame process $\mathcal{T}^{B}$%
}

\begin{proposition}
\label{Propn: Finite p-varn and convergence of dyad poly} $\mathcal{T}^{B}$
has finite $\acute{p}$-variation for $\acute{p}>2$, $\mathbb{P}$-a.s.
Furthermore, the dyadic polygonal approxmiations $\mathcal{T}^{B}\left(
m\right) $ as defined in (\ref{dyadic polygonal approximations}) converge to 
$\mathcal{T}^{B}$ in $\acute{p}$-variation $\mathbb{P}$-a.s.
\end{proposition}

\begin{proof}
Fix $n\in \mathbb{N}$ and $\acute{p}>2+\frac{1}{n}$. Then fix $p>$ $\acute{p}
$ so large that 
\begin{equation*}
\frac{2p}{p-2}<2+\frac{1}{n}.
\end{equation*}%
Since $\left\Vert .\right\Vert _{p}\geq \left\Vert .\right\Vert _{\infty }$,
using Proposition \ref{Hoelder's condition}, we have that 
\begin{eqnarray*}
&&\mathbb{E}\left[ \left\Vert \mathcal{T}_{h_{2}}^{B}-\mathcal{T}%
_{h_{1}}^{B}\right\Vert _{\infty }^{\acute{p}}\right] \\
&\leq &\mathbb{E}\left[ \left\Vert \mathcal{T}_{h_{2}}^{B}-\mathcal{T}%
_{h_{1}}^{B}\right\Vert _{p}^{\acute{p}}\right] \\
&=&d\left( \alpha ,p\right) ^{\frac{\acute{p}}{p}}\left( h_{2}-h_{1}\right)
^{\left( \frac{1}{2}-\frac{1}{p}\right) \acute{p}}.
\end{eqnarray*}%
By choice 
\begin{equation*}
\acute{p}>2+\frac{1}{n}>\frac{2p}{p-2}=\frac{1}{\frac{1}{2}-\frac{1}{p}},
\end{equation*}%
so that the statement of Proposition \ref{Lyons' Holder Proposition} applies
with $\kappa =\frac{1}{2}-\frac{1}{p}$ on the set $\Omega _{2+\frac{1}{n}}$
of full $\mathbb{P}$-measure. The statement of this Proposition then holds
on the set $\cap _{n\in \mathbb{N}}\Omega _{2+\frac{1}{n}}$ of full $\mathbb{%
P}$-measure.
\end{proof}

\chapter{\protect\bigskip A Tail Estimate for $\left\Vert \mathcal{S}%
_{h_{2}}^{B,p}-\mathcal{S}_{h_{1}}^{B,p}\right\Vert _{p}$}

\section{Notation}

As before, we work on the Wiener space $\left( C_{0}\left( \left[ -1,1\right]
\right) ,\sigma _{\left\Vert .\right\Vert _{\infty }},\mathbb{P}\right) $.

$W^{1,2}\left( \left[ -1,1\right] \right) $ denotes the Sobolev space of
differentiable functions with derivative in $L^{2}\left( \left[ -1,1\right]
\right) $, i.e.%
\begin{equation*}
W^{1,2}\left( \left[ -1,1\right] \right) :=\left\{ F\in C_{0}\left( \left[
-1,1\right] \right) :\frac{dF}{dx}\in L^{2}\left( \left[ -1,1\right] \right)
\right\} \text{,}
\end{equation*}%
equipped with its Hilbert space norm%
\begin{equation*}
\left\Vert F\right\Vert _{W^{1,2}\left( \left[ -1,1\right] \right)
}:=\left\Vert \frac{dF}{dx}\right\Vert _{L^{2}\left( \left[ -1,1\right]
\right) }\text{.}
\end{equation*}%
In the Wiener space setting, $\left( W^{1,2}\left( \left[ -1,1\right]
\right) ,\left\Vert .\right\Vert _{W^{1,2}\left( \left[ -1,1\right] \right)
}\right) $ is also known as the \emph{Cameron-Martin space}.

$\mathcal{O}_{W^{1,2}\left( \left[ -1,1\right] \right) }=\mathcal{O}$
denotes the closed unit ball of the \emph{Cameron-Martin space}, i.e.%
\begin{equation*}
\mathcal{O}:=\left\{ F\in W^{1,2}\left( \left[ -1,1\right] \right)
:\left\Vert \frac{dF}{dx}\right\Vert _{L^{2}\left( \left[ -1,1\right]
\right) }\leq 1\right\} .
\end{equation*}%
$\ j$ denotes the continuous imbedding map 
\begin{equation*}
j:\left( W^{1,2}\left( \left[ -1,1\right] \right) ,\left\Vert .\right\Vert
_{W^{1,2}\left( \left[ -1,1\right] \right) }\right) \hookrightarrow \left(
C_{0}\left( \left[ -1,1\right] \right) ,\left\Vert .\right\Vert _{\infty
}\right) .
\end{equation*}%
To see that $j$ is continous, fix $\varepsilon >0$ and choose $f,g\in
W^{1,2}\left( \left[ -1,1\right] \right) $ such that $\left\Vert
f-g\right\Vert _{W^{1,2}\left( \left[ -1,1\right] \right) }<\frac{%
\varepsilon }{\sqrt{2}}$. Then%
\begin{eqnarray*}
&&\sup_{t\in \left[ -1,1\right] }\left\vert j\left( f\right) \left( t\right)
-j\left( g\right) \left( t\right) \right\vert \\
&=&\sup_{t\in \left[ -1,1\right] }\left\vert \int_{-1}^{t}\left( \frac{df}{du%
}-\frac{dg}{du}\right) du\right\vert \\
&\leq &\sup_{t\in \left[ -1,1\right] }\left\vert t+1\right\vert ^{\frac{1}{2}%
}\left[ \int_{-1}^{1}\left( \frac{df}{du}-\frac{dg}{du}\right) ^{2}du\right]
^{\frac{1}{2}} \\
&=&\sup_{t\in \left[ -1,1\right] }\left\vert t+1\right\vert ^{\frac{1}{2}%
}\left\Vert f-g\right\Vert _{W^{1,2}\left( \left[ -1,1\right] \right)
}<\varepsilon .
\end{eqnarray*}%
Hence, $\left\Vert .\right\Vert _{\infty }$ is $\sigma _{\left\Vert
.\right\Vert _{W^{1,2}\left( \left[ -1,1\right] \right) }}$-measurable where 
$\sigma _{\left\Vert .\right\Vert _{W^{1,2}\left( \left[ -1,1\right] \right)
}}$ denotes the Borel-$\sigma $-algebra on $W^{1,2}\left( \left[ -1,1\right]
\right) .$

We say that a Borel random variable $\mathcal{F}$ on Wiener space into some
Banach space $\left( V,\left\Vert .\right\Vert _{V}\right) $ is Lipschitz in
the direction of the Cameron-Martin space (or $W^{1,2}\left( \left[ -1,1%
\right] \right) $\emph{-Lipschitz})\emph{\ }if 
\begin{equation*}
Lip_{W^{1,2}\left( \left[ -1,1\right] \right) }\left( \mathcal{F}\right)
:=\sup_{x\in C_{0}\left( \left[ -1,1\right] \right) }\sup_{r>0}\sup_{y\in
j\left( \mathcal{O}\right) \setminus \left\{ 0\right\} }\frac{\left\Vert 
\mathcal{F}\left( x+ry\right) -\mathcal{F}\left( x\right) \right\Vert _{V}}{r%
}
\end{equation*}%
is finite. $Lip_{W^{1,2}\left( \left[ -1,1\right] \right) }\left( \mathcal{F}%
\right) $ is called the Lipschitz-norm of $\mathcal{F}$ (in the direction of 
$W^{1,2}\left( \left[ -1,1\right] \right) $).

$\Phi $ denotes the cumulative distribution function of the Gaussian measure
on $\mathbb{R}$, i.e.%
\begin{equation*}
\Phi \left( x\right) :=\frac{1}{\sqrt{2\pi }}\int_{-\infty }^{x}\exp \left( -%
\frac{u^{2}}{2}\right) du.
\end{equation*}%
If $A$ and $B$ are subsets of a vector space $V$, we define 
\begin{equation*}
A+B:=\left\{ a+b:a\in A,b\in B\right\}
\end{equation*}%
and for any $r\in \mathbb{R}$,%
\begin{equation*}
rA:=\left\{ ra:a\in A\right\} .
\end{equation*}%
A median $m_{X}$ of a real-valued random variable $X$ is defined as a value
in $\mathbb{R}$ with the property that%
\begin{equation*}
\mathbb{P}\left( X\leq m_{X}\right) =\mathbb{P}\left( X>m_{X}\right) =\frac{1%
}{2}.
\end{equation*}

\section{Main Results}

In this chapter, we are concerned with the $p$-variation frame process $%
\mathcal{S}^{B,p}$ ($p>2$) as defined in (\ref{frame process definition}).
We find two constants $d_{1}\left( \alpha ,p\right) $ (where $\alpha >1-%
\frac{1}{p}$) and $d_{2}\left( p\right) $ so that the random variable 
\begin{equation*}
\frac{\left\Vert \mathcal{S}_{h_{2}}^{B,p}-\mathcal{S}_{h_{1}}^{B,p}\right%
\Vert _{p}}{d_{2}\left( p\right) \left( h_{2}-h_{1}\right) ^{\frac{1}{2}-%
\frac{1}{p}}}-d_{1}\left( \alpha ,p\right)
\end{equation*}%
has Gaussian tails, i.e.%
\begin{equation*}
\mathbb{P}\left( \frac{\left\Vert \mathcal{S}_{h_{2}}^{B,p}-\mathcal{S}%
_{h_{1}}^{B,p}\right\Vert _{p}}{d_{2}\left( p\right) \left(
h_{2}-h_{1}\right) ^{\frac{1}{2}-\frac{1}{p}}}-d_{1}\left( \alpha ,p\right)
\geq r\right) \leq \frac{1}{\sqrt{2\pi }r}\exp \left( -\frac{r^{2}}{2}%
\right) .
\end{equation*}%
A concentration of measure result for Wiener space (due to Borell \ and
independently roughly at the same time to Sudakov\emph{\ }\ and T'sirelson)
is applied to the functional $\left\Vert \mathcal{S}_{h_{2}}^{B,p}-\mathcal{S%
}_{h_{1}}^{B,p}\right\Vert _{p}$ on Wiener space in order to derive the
above bound.

\section{Borell's inequality}

\subsection{Borell's inequality on $\left( C_{0}\left( \left[ -1,1\right]
\right) ,\protect\sigma _{\left\Vert .\right\Vert _{\infty }},\mathbb{P}%
\right) $}

The following inequality due to \emph{Borell }(Theorem 3.1 in \cite{Borell75}%
\emph{\ }and \newline \emph{\ }\cite{Sudakov74})\emph{\ }is at the heart of our
argument -- we refer to it as \emph{Borell's inequality:}

\begin{theorem}[Borell's inequality]
If $A$ $\in \sigma _{\left\Vert .\right\Vert _{\infty }}$ and $r>0$, then we
have the following lower bound for the enlargement of $A$ in the direction
of the unit ball in the Cameron-Martin space $\mathcal{O}$:%
\begin{equation}
\mathbb{P}\left( A+rj\left( \mathcal{O}\right) \right) \geq \Phi \left( \Phi
^{-1}\left( \mathbb{P}\left( A\right) \right) +r\right) .
\label{Borell's inequality}
\end{equation}
\end{theorem}

\begin{remark}
It is a fact that $A+rj\left( \mathcal{O}\right) $ is in $\sigma
_{\left\Vert .\right\Vert _{\infty }}$ (\cite{Ledoux94}).
\end{remark}

\subsection{An Application of Borell's inequality to non-negative
functionals on Wiener space}

We are interested in applying Borell's inequality to real-valued
non-negative random variables on Wiener space.

\begin{lemma}
\label{Lemma: Deviation inequality} Suppose that $\mathcal{F}$ is a
non-negative random variable that has finite Lipschitz norm in the direction
of $W^{1,2}\left( \left[ -1,1\right] \right) $, i.e. 
\begin{equation*}
Lip_{W^{1,2}\left( \left[ -1,1\right] \right) }\left( \mathcal{F}\right)
<\infty .
\end{equation*}%
Then%
\begin{equation}
\mathbb{P}\left( \mathcal{F}>2\mathbb{E}\left[ \mathcal{F}\right]
+rLip_{W^{1,2}\left( \left[ -1,1\right] \right) }\left( \mathcal{F}\right)
\right) \leq \frac{1}{\sqrt{2\pi }r}\exp \left( -\frac{r^{2}}{2}\right) .
\label{Deviation inequality}
\end{equation}
\end{lemma}

\begin{proof}
Let $m_{\mathcal{F}}$ denote the median of $\mathcal{F}$. By \emph{Borell's
inequality, }%
\begin{equation*}
\mathbb{P}\left( \left\{ \mathcal{F}\leq m_{\mathcal{F}}\right\} +rj\left( 
\mathcal{O}\right) \right) \geq \Phi \left( r\right) .
\end{equation*}%
But 
\begin{eqnarray*}
&&\left\{ \mathcal{F}\leq m_{\mathcal{F}}\right\} +rj\left( \mathcal{O}%
\right) \\
&\subset &\left\{ y\in C_{0}\left( \left[ -1,1\right] \right) :\mathcal{F}%
\left( y\right) \leq m_{\mathcal{F}}+rLip_{W^{1,2}\left( \left[ -1,1\right]
\right) }\left( \mathcal{F}\right) \right\}
\end{eqnarray*}%
(which is a Borel set). Since $\mathcal{F}$ is \emph{non-negative}, we have
that 
\begin{equation*}
\mathbb{E}\left[ \mathcal{F}\right] \geq \int_{\mathcal{F}\geq m_{F}}%
\mathcal{F}dP\geq \frac{1}{2}m_{\mathcal{F}}.
\end{equation*}%
Therefore, 
\begin{equation*}
\mathbb{P}\left( \mathcal{F}\leq 2\mathbb{E}\left[ \mathcal{F}\right]
+rLip_{W^{1,2}\left( \left[ -1,1\right] \right) }\left( \mathcal{F}\right)
\right) \geq \Phi \left( r\right) .
\end{equation*}%
Since we have that 
\begin{equation*}
1-\Phi \left( r\right) \leq \frac{1}{\sqrt{2\pi }}\int_{r}^{\infty }\frac{x}{%
r}\exp \left( -\frac{x^{2}}{2}\right) dx=\frac{1}{\sqrt{2\pi }r}\exp \left( -%
\frac{r^{2}}{2}\right) ,
\end{equation*}%
the result follows.
\end{proof}

\section{Tail Estimate for $\left\Vert \mathcal{S}_{h_{2}}^{B,p}-\mathcal{S}%
_{h_{1}}^{B,p}\right\Vert _{p}$}

\subsection{Lipschitz-norm bounds}

We start with an observation on how to bound the Lipschitz-norm of a
semi-norm:

\begin{lemma}
\label{Lemma: Lipschitz-norm of seminorm functional}If $\mathcal{F}$ is a
seminorm, then%
\begin{equation*}
Lip_{W^{1,2}\left( \left[ -1,1\right] \right) }\left( \mathcal{F}\right)
\leq \sup_{y\in \mathcal{O}}\left\vert \mathcal{F}\left( j\left( y\right)
\right) \right\vert .
\end{equation*}
\end{lemma}

\begin{proof}
This follows from the triangle inequality.
\end{proof}

Next, we find the exact Lipschitz-norm for the non-negative functional 
\begin{equation*}
Lip_{W^{1,2}\left( \left[ -1,1\right] \right) }\left( \left\Vert \mathcal{T}%
_{h_{2}}-\mathcal{T}_{h_{1}}\right\Vert _{\infty }\right) .
\end{equation*}

\begin{lemma}
\label{Lemma:Lipschitz bound frame process sup norm}For any fixed $0\leq
h_{1}\leq h_{2}\leq 1$,%
\begin{equation*}
Lip_{W^{1,2}\left( \left[ -1,1\right] \right) }\left( \left\Vert \mathcal{T}%
_{h_{2}}-\mathcal{T}_{h_{1}}\right\Vert _{\infty }\right) =\left(
h_{2}-h_{1}\right) ^{\frac{1}{2}}\text{.}
\end{equation*}
\end{lemma}

\begin{proof}
We fix $f\in L^{2}\left( \left[ -1,1\right] \right) $ and note that 
\begin{eqnarray}
&&\mathcal{T}_{h_{2}}^{\int_{-1}^{.}f\left( s\right) ds}-\mathcal{T}%
_{h_{1}}^{\int_{-1}^{.}f\left( s\right) ds}  \notag \\
&=&\left( \int_{h_{1}-1+u}^{h_{2}-1+u}f\left( s\right) ds\right) _{0\leq
u\leq 1}.  \label{reexpressing frame process eval on W1,2}
\end{eqnarray}%
Hence,%
\begin{eqnarray*}
&&\sup_{t\in \left[ 0,1\right] }\left\vert \left( \mathcal{T}%
_{h_{2}}^{\int_{-1}^{.}f\left( s\right) ds}-\mathcal{T}_{h_{1}}^{%
\int_{-1}^{.}f\left( s\right) ds}\right) \left( t\right) \right\vert \\
&=&\sup_{t\in \left[ 0,1\right] }\left\vert
\int_{h_{1}-1+t}^{h_{2}-1+t}f\left( s\right) ds\right\vert \\
&\leq &\left( h_{2}-h_{1}\right) ^{\frac{1}{2}}\left\Vert f\right\Vert
_{L^{2}\left( \left[ -1,1\right] \right) }.
\end{eqnarray*}%
Furthermore, $\left\Vert \mathcal{T}_{h_{2}}^{.}-\mathcal{T}%
_{h_{1}}^{.}\right\Vert _{\infty }$ is a seminorm on Wiener space so that
using Lemma \ref{Lemma: Lipschitz-norm of seminorm functional}, we have
shown that 
\begin{equation*}
Lip_{W^{1,2}\left( \left[ -1,1\right] \right) }\left( \left\Vert \mathcal{T}%
_{h_{2}}-\mathcal{T}_{h_{1}}\right\Vert _{\infty }\right) \leq \left(
h_{2}-h_{1}\right) ^{\frac{1}{2}}.
\end{equation*}

In order to show that $Lip_{W^{1,2}\left( \left[ -1,1\right] \right) }\left(
\left\Vert \mathcal{T}_{h_{2}}-\mathcal{T}_{h_{1}}\right\Vert _{\infty
}\right) $ is at least $\left( h_{2}-h_{1}\right) ^{\frac{1}{2}}$, we
consider 
\begin{equation*}
g\left( s\right) :=\frac{1}{\sqrt{h_{2}-h_{1}}}\mathbf{1}_{\left[
h_{1}-1\leq s<h_{2}-1\right] },
\end{equation*}%
and note that $\int_{0}^{.}g\left( s\right) ds\in \mathcal{O}$. Furthermore,
the supremum%
\begin{eqnarray*}
&&\sup_{t\in \left[ 0,1\right] }\left\vert \left( \mathcal{T}%
_{h_{2}}^{\int_{-1}^{.}g\left( s\right) ds}-\mathcal{T}_{h_{1}}^{%
\int_{-1}^{.}g\left( s\right) ds}\right) \left( t\right) \right\vert \\
&=&\frac{1}{\sqrt{h_{2}-h_{1}}}\sup_{t\in \left[ 0,1\right] }\left\vert
\int_{h_{1}-1+t}^{h_{2}-1+t}\mathbf{1}_{\left[ h_{1}-1\leq s<h_{2}-1\right]
}ds\right\vert
\end{eqnarray*}%
is attained at $t=0$ and is equal to $\sqrt{h_{2}-h_{1}}$ which completes
the proof.
\end{proof}

We find an upper bound for the Lipschitz-norm of $\mathcal{V}_{p}\left( 
\mathcal{T}_{h_{2}}^{.}-\mathcal{T}_{h_{1}}^{.}\right) $:

\begin{lemma}
\label{Lemma: Lipschitz bound frame process p-functional}For $p>2$ and any
fixed pair $h_{1}$ and $h_{2}$ with $0\leq h_{1}\leq h_{2}\leq 1$,%
\begin{equation*}
Lip_{W^{1,2}\left( \left[ -1,1\right] \right) }\left( \mathcal{V}_{p}\left( 
\mathcal{T}_{h_{2}}^{.}-\mathcal{T}_{h_{1}}^{.}\right) \right) \leq
d_{p}\left( h_{2}-h_{1}\right) ^{\frac{1}{2}-\frac{1}{p}},
\end{equation*}%
where%
\begin{equation*}
d_{p}:=2^{\frac{1}{p}+\frac{1}{2}}\left( 1+2^{\frac{p}{2}}\right) ^{\frac{1}{%
p}}.
\end{equation*}
\end{lemma}

\begin{proof}
We fix a dissection $D=\left\{ t_{0}=-1;t_{1};...;t_{n}=1\right\} $ of $%
\left[ -1,1\right] $ and $f\in L^{2}\left( \left[ -1,1\right] \right) $.
Using the representation (\ref{reexpressing frame process eval on W1,2}), we
find that%
\begin{eqnarray*}
&&\sum_{t_{i},t_{i+1}\in D}\left\vert \left( \mathcal{T}_{h_{2}}^{%
\int_{-1}^{.}f\left( s\right) ds}-\mathcal{T}_{h_{1}}^{\int_{-1}^{.}f\left(
s\right) ds}\right) \left( t_{i+1}\right) \right. \\
&&\left. -\left( \mathcal{T}_{h_{2}}^{\int_{-1}^{.}f\left( s\right) ds}-%
\mathcal{T}_{h_{1}}^{\int_{-1}^{.}f\left( s\right) ds}\right) \left(
t_{i}\right) \right\vert ^{p} \\
&=&\sum_{t_{i},t_{i+1}\in D}\left\vert \int_{\left( h_{2}-1+t_{i}\right)
\vee \left( h_{1}-1+t_{i+1}\right) }^{h_{2}-1+t_{i+1}}f\left( u\right)
du-\int_{h_{1}-1+t_{i}}^{\left( h_{2}-1+t_{i}\right) \wedge \left(
h_{1}-1+t_{i+1}\right) }f\left( u\right) du\right\vert ^{p}.
\end{eqnarray*}%
For any dissection piece $\left[ t_{i},t_{i+1}\right] $, by the
Cauchy-Schwarz inequality, we have that 
\begin{eqnarray*}
&&\left\vert \left( \int_{\left( h_{2}-1+t_{i}\right) \vee \left(
h_{1}-1+t_{i+1}\right) }^{h_{2}-1+t_{i+1}}-\int_{h_{1}-1+t_{i}}^{\left(
h_{2}-1+t_{i}\right) \wedge \left( h_{1}-1+t_{i+1}\right) }\right) f\left(
u\right) du\right\vert \\
&\leq &\sqrt{2}\left\Vert f\right\Vert _{L^{2}\left( \left[ -1,1\right]
\right) }\min \left( \left( t_{i+1}-t_{i}\right) ^{\frac{1}{2}},\left(
h_{2}-h_{1}\right) ^{\frac{1}{2}}\right) .
\end{eqnarray*}%
Hence,%
\begin{eqnarray*}
&&\sum_{t_{i},t_{i+1}\in D}\left\vert \left( \mathcal{T}_{h_{2}}^{%
\int_{-1}^{.}f\left( s\right) ds}-\mathcal{T}_{h_{1}}^{\int_{-1}^{.}f\left(
s\right) ds}\right) \left( t_{i+1}\right) \right. \\
&&\left. \left( \mathcal{T}_{h_{2}}^{\int_{-1}^{.}f\left( s\right) ds}-%
\mathcal{T}_{h_{1}}^{\int_{-1}^{.}f\left( s\right) ds}\right) \left(
t_{i}\right) \right\vert ^{p} \\
&\leq &2^{\frac{p}{2}}\left\Vert f\right\Vert _{L^{2}\left( \left[ -1,1%
\right] \right) }^{p}\sum_{t_{i},t_{i+1}\in D}\min \left( \left(
t_{i+1}-t_{i}\right) ^{\frac{p}{2}},\left( h_{2}-h_{1}\right) ^{\frac{p}{2}%
}\right) .
\end{eqnarray*}%
As in the proof of Proposition \ref{Propn: pathwise p-varn bound}\ in
Chapter 2, we write $D$ as a disjoint union of $D^{<h}$ and $D^{\geq h}$
where $h=h_{2}-h_{1}$ and 
\begin{equation*}
D^{<h}:=\left\{ t_{i}\in D:t_{i+1}-t_{i}<h\right\}
\end{equation*}%
and 
\begin{equation*}
D^{\geq h}:=\left\{ t_{i}\in D:t_{i+1}-t_{i}\geq h\right\} .
\end{equation*}%
Accordingly, we split the sum:%
\begin{eqnarray*}
&&\sum_{t_{i},t_{i+1}\in D}\min \left( \left( t_{i+1}-t_{i}\right) ^{\frac{p%
}{2}},\left( h_{2}-h_{1}\right) ^{\frac{p}{2}}\right) \\
&=&\underset{\text{Sum }A}{\underbrace{\sum_{t_{i},t_{i+1}\in D^{\geq h}}h^{%
\frac{p}{2}}}}+\underset{\text{Sum }B}{\underbrace{\sum_{t_{i},t_{i+1}\in
D^{<h}}\left( t_{i+1}-t_{i}\right) ^{\frac{p}{2}}}}.
\end{eqnarray*}%
Since there are at most $\frac{2}{h}$ terms in Sum $A$ (an interval of
length $2$ can contain at most $\frac{2}{h}$ disjoint intervals of length at
least $h$), Sum $A$ is bounded by $2\left( h_{2}-h_{1}\right) ^{\frac{p}{2}%
-1}$.

In order to bound Sum $B$, we start by listing the dissection pieces with
endpoints in $D^{<h}$: 
\begin{equation*}
\left\{ \left[ t_{i_{1}},t_{i_{1}+1}\right] ,\left[ t_{i_{2}},t_{i_{2}+1}%
\right] ,...,\left[ t_{i_{m}},t_{i_{m}+1}\right] \right\} .
\end{equation*}%
We now consider the following subsequence of $\left\{ t_{i_{j}}:1\leq j\leq
m\right\} $: Let $t_{i_{j_{1}}}$ be the first element in the sequence such
that the summed length of all dissection pieces $\left[ t_{i_{l}},t_{i_{l}+1}%
\right] $ with index $i_{l}$ up to and including $i_{j_{1}}$ is at least $h$
and at most $2h$ (if no such index exists, then \ the sum of all dissection
pieces in $D^{<h}$ adds up to strictly less than $h$, so that $%
\sum_{t_{i},t_{i+1}\in D^{<h}}\left( t_{i+1}-t_{i}\right) ^{\frac{p}{2}}<h^{%
\frac{p}{2}}$). In this way, recursively define $t_{i_{j_{k}}}$ such that
the overall length of all dissection pieces whose endpoints are included in $%
D^{<h}$ with index $i_{l}$ between $i_{j_{k-1}+1}$ and $i_{j_{k}}$ is at
least $h$ and at most $2h$. Since an interval of length $2$ can contain at
most $\frac{2}{h}$ disjoint intervals of length at least $h$, the
subsequence $\left\{ i_{j_{k}}\right\} $ thus defined has at most $\frac{2}{h%
}$ elements. Furthermore since $p>2$, for every $k$%
\begin{equation*}
\sum_{l=j_{k-1}+1}^{j_{k}}\left( t_{i_{l}+1}-t_{i_{l}}\right) ^{\frac{p}{2}%
}\leq \left( 2h\right) ^{\frac{p}{2}}.
\end{equation*}%
Hence,%
\begin{equation*}
\sum_{t_{i},t_{i+1}\in D^{<h}}\left( t_{i+1}-t_{i}\right) ^{\frac{p}{2}}\leq
2^{\frac{p}{2}+1}h^{\frac{p}{2}-1}.
\end{equation*}%
Therefore,%
\begin{equation*}
Lip_{W^{1,2}\left( \left[ -1,1\right] \right) }\left( \mathcal{V}_{p}\left( 
\mathcal{T}_{h_{2}}-\mathcal{T}_{h_{1}}\right) \right) \leq \underset{:=d_{p}%
}{\underbrace{2^{\frac{1}{2}+\frac{1}{p}}\left( 1+2^{\frac{p}{2}}\right) ^{%
\frac{1}{p}}}}\left( h_{2}-h_{1}\right) ^{\frac{1}{2}-\frac{1}{p}}.
\end{equation*}
\end{proof}

\subsection{Tail Estimate}

\begin{theorem}
If $p>2$ , $\ \alpha >1-\frac{1}{p}$ and $0\leq h_{1}<h_{2}\leq 1$, then 
\begin{eqnarray*}
&&\mathbb{P}\left( \frac{\left\Vert \mathcal{S}_{h_{2}}^{B,p}-\mathcal{S}%
_{h_{1}}^{B,p}\right\Vert _{p}}{d_{2}\left( p\right) \left(
h_{2}-h_{1}\right) ^{\frac{1}{2}-\frac{1}{p}}}-d_{1}\left( \alpha ,p\right)
\geq r\right) \\
&\leq &\frac{1}{\sqrt{2\pi }}\frac{1}{r}\exp \left( -\frac{r^{2}}{2}\right) ,
\end{eqnarray*}%
where 
\begin{equation*}
d_{1}\left( \alpha ,p\right) =\frac{d\left( \alpha ,p\right) ^{\frac{1}{p}}}{%
d_{p}}
\end{equation*}%
and 
\begin{equation*}
d_{2}\left( p\right) =2d_{p}.
\end{equation*}%
$d\left( \alpha ,p\right) $ is the constant given in Proposition \ref%
{Hoelder's condition}. $d_{p}$ is the constant given in Lemma \ref{Lemma:
Lipschitz bound frame process p-functional}.
\end{theorem}

\begin{proof}
We recall (\ref{Deviation inequality}) from Lemma \ref{Lemma: Deviation
inequality},%
\begin{equation*}
\mathbb{P}\left( \mathcal{F}>2\mathbb{E}\left[ \mathcal{F}\right]
+rLip_{W^{1,2}\left( \left[ -1,1\right] \right) }\left( \mathcal{F}\right)
\right) \leq \frac{1}{\sqrt{2\pi }r}\exp \left( -\frac{r^{2}}{2}\right) ,
\end{equation*}%
which we will apply to 
\begin{equation*}
\mathcal{F}=\left\Vert \mathcal{S}_{h_{2}}^{B,p}-\mathcal{S}%
_{h_{1}}^{B,p}\right\Vert _{p}
\end{equation*}%
-- noting in what follows that $\mathcal{S}^{B,p}$ is a continuous
modification of $\mathcal{T}^{B,p}$. From Chapter 2, Proposition \ref%
{Hoelder's condition} , we recall that 
\begin{equation*}
\mathbb{E}\left[ \left\Vert \mathcal{S}_{h_{2}}^{B,p}-\mathcal{S}%
_{h_{1}}^{B,p}\right\Vert _{p}^{p}\right] \leq d\left( \alpha ,p\right)
\left( h_{2}-h_{1}\right) ^{\left( \frac{1}{2}-\frac{1}{p}\right) p},
\end{equation*}%
so that by Lyapunov's inequality%
\begin{equation*}
\mathbb{E}\left[ \left\Vert \mathcal{S}_{h_{2}}^{B,p}-\mathcal{S}%
_{h_{1}}^{B,p}\right\Vert _{p}\right] \leq \mathbb{E}\left[ \left\Vert 
\mathcal{S}_{h_{2}}^{B,p}-\mathcal{S}_{h_{1}}^{B,p}\right\Vert _{p}^{p}%
\right] ^{\frac{1}{p}}\leq d\left( \alpha ,p\right) ^{\frac{1}{p}}\left(
h_{2}-h_{1}\right) ^{\frac{1}{2}-\frac{1}{p}}.
\end{equation*}%
From Lemma \ref{Lemma:Lipschitz bound frame process sup norm} and Lemma%
\emph{\ }\ref{Lemma: Lipschitz bound frame process p-functional}, we have
that on a set of full $\mathbb{P}$-measure 
\begin{eqnarray*}
Lip_{W^{1,2}\left( \left[ -1,1\right] \right) }\left( \left\Vert \mathcal{S}%
_{h_{2}}^{B,p}-\mathcal{S}_{h_{1}}^{B,p}\right\Vert _{p}\right)
&=&Lip_{W^{1,2}\left( \left[ -1,1\right] \right) }\left( \left\Vert \mathcal{%
T}_{h_{2}}-\mathcal{T}_{h_{1}}\right\Vert _{p}\right) \\
&\leq &2d_{p}\left( h_{2}-h_{1}\right) ^{\frac{1}{2}-\frac{1}{p}}.
\end{eqnarray*}%
Hence,%
\begin{eqnarray*}
&&\mathbb{P}\left( \mathcal{F}>2d\left( \alpha ,p\right) ^{\frac{1}{p}%
}\left( h_{2}-h_{1}\right) ^{\frac{1}{2}-\frac{1}{p}}+2d_{p}\left(
h_{2}-h_{1}\right) ^{\frac{1}{2}-\frac{1}{p}}r\right) \\
&\leq &\frac{1}{\sqrt{2\pi }r}\exp \left( -\frac{r^{2}}{2}\right) ,
\end{eqnarray*}%
so that 
\begin{equation*}
\mathbb{P}\left( \frac{\left\Vert \mathcal{S}_{h_{2}}^{B,p}-\mathcal{S}%
_{h_{1}}^{B,p}\right\Vert _{p}}{2d_{p}\left( h_{2}-h_{1}\right) ^{\frac{1}{2}%
-\frac{1}{p}}}>\frac{d\left( \alpha ,p\right) ^{\frac{1}{p}}}{d_{p}}%
+r\right) \leq \frac{1}{\sqrt{2\pi }r}\exp \left( -\frac{r^{2}}{2}\right) .
\end{equation*}
\end{proof}

\chapter{\protect\bigskip Non-existence of L\'{e}vy Area for the frame
process $\mathcal{T}^{B}$}

\section{Notation}

As in previous chapters, $C\left( \left[ 0,1\right] \right) $ denotes the
Banach space of real-valued continuous functions equipped with the $\sup $%
-norm. Similarly, $C\left( \left[ 0,1\right] \times \left[ 0,1\right]
\right) $ denotes the Banach space of real-valued continuous functions on
the unit square equipped with the $\sup $-norm.

If $V$ is a normed space with norm $\left\vert .\right\vert _{V}$, $\mathcal{%
O}_{V}$ denotes the closed unit ball of $V$, i.e.%
\begin{equation*}
\mathcal{O}_{V}:=\left\{ v\in V:\left\vert v\right\vert _{V}\leq 1\right\} .
\end{equation*}

Suppose $V$ and $W$ are Banach spaces.\ $\mathcal{B}\left( V,W\right) $
denotes the Banach space of bounded linear maps from $V$ to $W$ equipped
with the operator norm.

$V^{\ast }$ denotes the topological dual of $V$, equipped with the operator
norm on $\mathcal{B}\left( V,\mathbb{R}\right) $, i.e. 
\begin{equation*}
\left\vert v^{\ast }\right\vert _{V^{\ast }}:=\sup_{x\in \mathcal{O}%
_{V}}v^{\ast }\left( x\right) .
\end{equation*}

If $\ x\in \left[ 0,1\right] $, then $\delta _{x}$ denotes the evaluation
functional at $x$, that is for any continuous function $f:\left[ 0,1\right]
\rightarrow \mathbb{R}$, 
\begin{equation*}
\delta _{x}\left( f\right) :=f\left( x\right) .
\end{equation*}%
A \emph{partition }of $\left[ 0,1\right] $ is defined as a countable
collection $\left\{ t_{0}=0,t_{1},t_{2}\,,...:t_{0}<t_{1}<t_{2}<...\right\} $
such that $\left[ t_{i},t_{i+1}\right) \cap \left[ t_{j},t_{j+1}\right)
=\emptyset $ if $i\neq j$ and $\cup _{i=1}^{\infty }\left[
t_{i},t_{i+1}\right) =\left[ 0,1\right) $ . The \emph{total variation} of a
signed Borel measure $\mu $ on $\left[ 0,1\right] $ is defined as 
\begin{equation*}
\left\vert \mu \right\vert _{M^{1}\left( \left[ 0,1\right] \right) }:=\sup
\left\{ \sum_{i=1}^{\infty }\left\vert \mu \left( \left[ t_{i},t_{i+1}%
\right) \right) \right\vert :\text{all partitions of }\left[ 0,1\right]
\right\} .
\end{equation*}

$M^{1}\left( \left[ 0,1\right] \right) $ denotes the space of signed Borel
mesures on $\left[ 0,1\right] $ equipped with the total variation norm $%
\left\vert .\right\vert _{M^{1}\left( \left[ 0,1\right] \right) }$.

\section{Main Results}

From Proposition \ref{Propn: Finite p-varn and convergence of dyad poly} in
Chapter 2 we know that the $\sup $-norm frame process $\mathcal{T}^{B}$ has
finite $\acute{p}$-variation if $\acute{p}>2$. Theorem \ref%
{FirstTheoremLyons} in Chapter 1 tells us that in order to establish a rough
path integration theory for $\mathcal{T}^{B}$, we need to find a lift of $%
\mathcal{T}^{B}$ to $\Omega _{2}\left( V\right) $ (where $V=C\left( \left[
0,1\right] \right) $). Proposition \ref{Propn: Finite p-varn and convergence
of dyad poly} furthermore asserts that the dyadic polygonal approximations $%
\mathcal{T}^{B}\left( m\right) $ converge to $\mathcal{T}^{B}$ in $\acute{p}$%
-variation norm. Chen's Theorem (Theorem \ref{Theorem: Chen's theorem})\
gives us a multiplicative lift for every $m$. Thus, one may be led to ask
whether $\mathcal{T}^{B}$ has a lift that is the limit of the dyadic
polygonal smooth rough paths associated to the sequence $\left( \mathcal{T}%
^{B}\left( m\right) \right) _{m\in \mathbb{N}}$. In this chapter, we show
that this is \emph{not }the case: If we consider the injective tensor
product $V\otimes _{\vee }V$, while 
\begin{equation}
\int \int_{0\leq u\leq v\leq 1}d\mathcal{T}^{B}\left( m\right) _{u}\otimes d%
\mathcal{T}^{B}\left( m\right) _{v}  \label{Eqn: Chen lift 2-tensor}
\end{equation}%
takes its values in $V\otimes _{\vee }V$ for every fixed $m$, the limit as $%
m\rightarrow \infty $ does not. The reason is the following: While $V\otimes
_{\vee }V$ is shown to be isomorphic to $C\left( \left[ 0,1\right] \times %
\left[ 0,1\right] \right) $, we prove that in the limit as $m\rightarrow
\infty $, the antisymmetric component of (\ref{Eqn: Chen lift 2-tensor})
given by 
\begin{eqnarray*}
&&\frac{1}{2}\lim_{m\rightarrow \infty }\left( \int \int_{0\leq u\leq v\leq
1}d\mathcal{T}^{B}\left( m\right) _{u}\otimes d\mathcal{T}^{B}\left(
m\right) _{v}-\right. \\
&&\left. \int \int_{0\leq u\leq v\leq 1}d\mathcal{T}^{B}\left( m\right)
_{v}\otimes d\mathcal{T}^{B}\left( m\right) _{u}\right) ,
\end{eqnarray*}%
exists and is continuous off the diagonal of the unit square. However, on
the diagonal it fails to be continuous. As we will see in the proof of
Proposition \ref{Propn: Representation of Levy Area random variable} and
Proposition \ref{Propn: discontinuity of LA} below, this breakdown in
continuity is intimately linked to the fact that the L\'{e}vy area of the
frame process $\mathcal{T}^{B}$ picks up the quadratic variation of the
Brownian sample path $B$.

This result has the following geometric interpretation: From L\'{e}vy's
modulus of continuity (e.g. \cite{Mckean69}), we know that for $\mathbb{P}$%
-a.e. sample path $f\in C_{0}\left( \left[ -1,1\right] \right) $, 
\begin{equation*}
\overline{\lim }_{h\rightarrow 0^{+}}\frac{\left\Vert \mathcal{T}_{1-h}^{f}-%
\mathcal{T}_{1}^{f}\right\Vert _{\infty }}{\sqrt{2h\ln \frac{1}{h}}}=1.
\end{equation*}%
Hence, the planar path 
\begin{equation*}
\left( \mathcal{T}_{1-h}^{f},\mathcal{T}_{1}^{f}\right) =\left( f\left(
t-h\right) ,f\left( t\right) \right) _{0\leq t\leq 1}
\end{equation*}%
converges uniformly to the path $\left( f\left( t\right) ,f\left( t\right)
\right) _{0\leq t\leq 1}$ at a rate $\sqrt{2h\ln \frac{1}{h}}$. The L\'{e}vy
area (Example \ref{Ex: R^2} \ below) of $\left( f\left( t\right) ,f\left(
t\right) \right) _{0\leq t\leq 1}$ is $0$. However, we will see that for $%
h_{n}\searrow 0$, even though $\left( \mathcal{T}_{1-h_{n}}^{f},\mathcal{T}%
_{1}^{f}\right) \rightarrow \left( f\left( t\right) ,f\left( t\right)
\right) _{0\leq t\leq 1}$ uniformly, the L\'{e}vy area of the path $\left( 
\mathcal{T}_{1-h_{n}}^{f},\mathcal{T}_{1}^{f}\right) $ converges to $-1/2$.
So the random path\ sequence $\left( \mathcal{T}_{1-h_{n}}^{B},\mathcal{T}%
_{1}^{B}\right) $ exhibits a similar behaviour to the deterministic path
sequence $\left( \frac{\cos n^{2}t}{n},\frac{\sin n^{2}t}{n}\right) $ known
from Example 1.1.1 in \cite{Lyons98}: Although $\left( \frac{\cos n^{2}t}{n},%
\frac{\sin n^{2}t}{n}\right) $ converges to $\left( 0,0\right) $ uniformly,
the associated sequence of L\'{e}vy areas converges to $1/2$. Similarly, $%
\left( \mathcal{T}_{1-h_{n}}^{f},\mathcal{T}_{1}^{f}\right) $ converges
uniformly to a process\ taking values on the diagonal -- which does not
generate area. But the associated sequence of L\'{e}vy areas tends to $-1/2$.

\section{The injective tensor algebra $C\left( \left[ 0,1\right] \right)
\otimes _{\vee }C\left( \left[ 0,1\right] \right) $}

\subsection{Tensor Products}

Suppose that $V$ is a Banach space with topological dual $V^{\ast }$. The
algebraic tensor product of $V$ with itself -- denoted as $V\otimes V$ -- is
the set of all elements $\sum_{i=1}^{n}x_{i}\otimes y_{i}$ where $%
x_{i},y_{i} $ $\in V$ and $n$ is finite: For any $u_{1},u_{2},u_{3}\in V$,
we have that 
\begin{equation*}
u_{1}\otimes \left( u_{2}+u_{3}\right) =u_{1}\otimes u_{2}+u_{1}\otimes
u_{3}.
\end{equation*}
The following Proposition asserts that any element $w\in V\otimes V$
identifies a unique finite rank element in $\mathcal{B}\left( V^{\ast
},V\right) $\ whose action on an element $x^{\ast }\in V^{\ast }$ is given
by 
\begin{equation*}
\text{\ }\bar{w}\left( x^{\ast }\right) =\sum_{i=1}^{n}x^{\ast }\left(
x_{i}\right) y_{i}.
\end{equation*}

\begin{proposition}[Lemma 1.2 in \protect\cite{Schatten50}]
\bigskip If $w\in V\otimes V$ has two representations, say%
\begin{equation*}
w=\sum_{i=1}^{n}x_{i}^{1}\otimes y_{i}^{1}=\sum_{j=1}^{m}x_{j}^{2}\otimes
y_{j}^{2},
\end{equation*}%
then for any $x^{\ast }\in V^{\ast }$, we have that%
\begin{equation*}
\sum_{i=1}^{n}x^{\ast }\left( x_{i}^{1}\right)
y_{i}^{1}=\sum_{j=1}^{m}x^{\ast }\left( x_{j}^{2}\right) y_{j}^{2}.
\end{equation*}%
Hence, every element of $V\otimes V$ identifies a unique finite rank element
of $\mathcal{B}\left( V^{\ast },V\right) $.
\end{proposition}

Since $V^{\ast }$ is also a Banach space, by replacing $V$ by $V^{\ast }$ in
the above Proposition, every element of $V^{\ast }\otimes V^{\ast }$
identifies a unique finite rank element of $\mathcal{B}\left( V^{\ast \ast
},V^{\ast }\right) $. In this way, we define a natural duality pairing
between $V^{\ast }\otimes V^{\ast }$ and $V\otimes V$ given by 
\begin{equation}
\left\langle x\otimes y,f\otimes g\right\rangle :=f\left( x\right) g\left(
y\right) .  \label{duality pairing}
\end{equation}

\begin{definition}
\label{Defn: symmetric + antisymmetric two-tensor}A tensor $w\in $ $V\otimes
V$ is called \emph{symmetric}, if for any $f,g\in V^{\ast }$, we have 
\begin{equation*}
\left\langle x\otimes y,f\otimes g\right\rangle =\left\langle x\otimes
y,g\otimes f\right\rangle .
\end{equation*}%
It is called \emph{antisymmetric} if for any $f,g\in V^{\ast }$, we have 
\begin{equation*}
\left\langle x\otimes y,f\otimes g\right\rangle =-\left\langle x\otimes
y,g\otimes f\right\rangle .
\end{equation*}
\end{definition}

Every tensor $w=\sum_{i=1}^{n}x_{i}\otimes y_{i}\in V\otimes V$ can be
written as the sum of a symmetric two-tensor and an antisymmetric two-tensor
:%
\begin{eqnarray}
&&\sum_{i=1}^{n}x_{i}\otimes y_{i}  \notag \\
&=&\frac{1}{2}\underset{\text{symmetric component}}{\underbrace{%
\sum_{i=1}^{n}\left( x_{i}\otimes y_{i}+y_{i}\otimes x_{i}\right) }}  \notag
\\
&&+\frac{1}{2}\underset{\text{antisymmetric component}}{\underbrace{%
\sum_{i=1}^{n}\left( x_{i}\otimes y_{i}-y_{i}\otimes x_{i}\right) }}.
\label{symmetrisation}
\end{eqnarray}

\subsection{Tensor Product Norms}

We recall the definition of a \emph{compatible tensor norm }on $V\otimes V$
as defined in Chapter 1, equation (\ref{Defn: compatible norm}):

\begin{definition}[compatible tensor norm]
A tensor norm $\left\Vert .\right\Vert _{c}$ on $V\otimes V$ is said to be
compatible if for any $v_{1},v_{2}\in V$ 
\begin{equation*}
\left\Vert v_{1}\otimes v_{2}\right\Vert _{c}\leq \left\vert
v_{1}\right\vert _{V}\left\vert v_{2}\right\vert _{V}.
\end{equation*}%
$V\otimes _{c}V$ denotes the closure of $V\otimes V$ in $\mathcal{B}\left(
V^{\ast },V\right) $ with respect to $\left\Vert .\right\Vert _{c}$.
\end{definition}

We next consider a subclass of compatible tensor norms:

\begin{definition}
\label{Defn: cross-norm}A tensor norm $\left\Vert .\right\Vert _{\times }$
on $V\otimes V$ is called a \emph{cross-norm} if for any $v_{1},v_{2}\in V$ 
\begin{equation*}
\left\Vert v_{1}\otimes v_{2}\right\Vert _{\times }=\left\vert
v_{1}\right\vert _{V}\left\vert v_{2}\right\vert _{V}.
\end{equation*}
\end{definition}

The \emph{injective tensor product norm} is an example of a cross norm:

\begin{definition}
\label{Ex: injective tensor product norm}\bigskip Suppose $V$ is a Banach
space. The closure of $V\otimes V$ in $\mathcal{B}\left( V^{\ast },V\right) $
is called the \emph{injective tensor product} -- denoted as $V\otimes _{\vee
}V$. As a consequence of the Hahn-Banach theorem (Theorem 4.3 b) in \cite%
{Rudin91}), for any $v\in V$, 
\begin{equation}
\left\vert x\right\vert _{V}=\sup_{x^{\ast }\in \mathcal{O}_{V^{\ast
}}}x^{\ast }\left( v\right) .  \label{eqn: Hahn-Banach}
\end{equation}%
Hence, the injective tensor product norm is given by%
\begin{equation*}
\left\Vert w\right\Vert _{\vee }:=\sup \left\{ \sum_{i=1}^{n}x^{\ast }\left(
x_{i}\right) y^{\ast }\left( y_{i}\right) :x^{\ast },y^{\ast }\in \mathcal{O}%
_{V^{\ast }}\text{ }\right\} .
\end{equation*}
\end{definition}

The injective tensor product norm is an example of a cross-norm and hence of
a compatible norm. It has the additional property that it is the smallest
cross norm whose dual norm $\left\Vert .\right\Vert _{\vee }^{\ast }$ on $%
V^{\ast }\otimes V^{\ast }$ -- given by%
\begin{equation*}
\left\Vert \sum_{i=1}^{n}f_{i}\otimes g_{i}\right\Vert _{\vee }^{\ast
}:=\sup \left\{ \sum_{i=1}^{n}\left\langle f_{i}\otimes g_{i},x\otimes
y\right\rangle :x\otimes y\in \mathcal{O}_{V\otimes _{\vee }V}\right\}
\end{equation*}%
-- is also cross (\cite{Paulsen91}).

\subsection{A representation result for $C\left( \left[ 0,1\right] \right)
\otimes _{\vee }C\left( \left[ 0,1\right] \right) $}

We find an explicit representation of $V\otimes _{\vee }V$ in the particular
case where 
\begin{equation*}
V=\left( C\left( \left[ 0,1\right] \right) ,\left\Vert .\right\Vert _{\infty
}\right) .
\end{equation*}

\begin{proposition}
\label{Proposition: injective tensor algebra}The injective tensor product $%
C\left( \left[ 0,1\right] \right) \otimes _{\vee }C\left( \left[ 0,1\right]
\right) $ is isomorphic to the space of continuous functions on the unit
square, i.e. 
\begin{equation*}
C\left( \left[ 0,1\right] \right) \otimes _{\vee }C\left( \left[ 0,1\right]
\right) \cong C\left( \left[ 0,1\right] \times \left[ 0,1\right] \right) .
\end{equation*}
\end{proposition}

\begin{proof}
Suppose $w\in C\left( \left[ 0,1\right] \right) \otimes _{\vee }C\left( %
\left[ 0,1\right] \right) $ has a representation%
\begin{equation*}
w=\sum_{i=1}^{n}f_{i}\otimes g_{i}.
\end{equation*}%
By the \emph{Radon-Riesz Theorem }(Theorem 6.19 in \cite{Rudin87}) the
topological dual of $\left( C\left( \left[ 0,1\right] \right) ,\left\Vert
.\right\Vert _{\infty }\right) $ is isomorphic to the space of Borel
measures of finite total variation -- denoted by $M^{1}\left( \left[ 0,1%
\right] \right) $ with the norm given by the total variation $\left\vert
.\right\vert _{M^{1}\left( \left[ 0,1\right] \right) }$. Thus, the injective
tensor norm has the following form 
\begin{equation*}
\left\Vert w\right\Vert _{\vee }=\sup \left\{
\sum_{i=1}^{n}\int_{0}^{1}\int_{0}^{1}f_{i}\left( u\right) g_{i}\left(
v\right) \mu \left( du\right) \nu \left( dv\right) :\mu ,\nu \in \mathcal{O}%
_{M^{1}\left( \left[ 0,1\right] \right) }\right\} .
\end{equation*}%
But if $\mu ,\nu \in \mathcal{O}_{M^{1}\left( \left[ 0,1\right] \right) }$,
then the product measure $\mu \times \nu \left( du,dv\right) :=\mu \left(
du\right) \nu \left( dv\right) $ is in $\mathcal{O}_{M^{1}\left( \left[ 0,1%
\right] \times \left[ 0,1\right] \right) }$. Hence,%
\begin{eqnarray*}
&&\sup \left\{ \sum_{i=1}^{n}\int_{0}^{1}\int_{0}^{1}f_{i}\left( u\right)
g_{i}\left( v\right) \mu \times \nu \left( du,dv\right) :\mu ,\nu \in 
\mathcal{O}_{M^{1}\left( \left[ 0,1\right] \right) }\left( 1\right) \right\}
\\
&\leq &\sup \left\{ \sum_{i=1}^{n}\int_{0}^{1}\int_{0}^{1}f_{i}\left(
u\right) g_{i}\left( v\right) \gamma \left( du,dv\right) :\gamma \in 
\mathcal{O}_{M^{1}\left( \left[ 0,1\right] \times \left[ 0,1\right] \right)
}\left( 1\right) \right\}
\end{eqnarray*}%
Using (\ref{eqn: Hahn-Banach}), it follows that%
\begin{equation*}
\left\vert \sum_{i=1}^{n}f_{i}\left( .\right) g_{i}\left( .\right)
\right\vert _{\infty }=\sup_{\gamma \in \mathcal{O}_{M^{1}\left( \left[ 0,1%
\right] \times \left[ 0,1\right] \right) }}\left\vert
\int_{0}^{1}\int_{0}^{1}\sum_{i=1}^{n}f_{i}\left( u\right) g_{i}\left(
v\right) \gamma \left( du,dv\right) \right\vert
\end{equation*}%
Hence,%
\begin{equation*}
\left\Vert \sum_{i=1}^{n}f_{i}\otimes g_{i}\right\Vert _{\vee }\leq
\left\vert \sum_{i=1}^{n}f_{i}\left( .\right) g_{i}\left( .\right)
\right\vert _{\infty }
\end{equation*}%
To prove the other direction, we note that since $f_{i},g_{i}\in C\left( %
\left[ 0,1\right] \right) $ for $1\leq i\leq n$, 
\begin{equation*}
w\left( .,.\right) =\sum_{i=1}^{n}f_{i}\left( .\right) g_{i}\left( .\right)
\end{equation*}
-- and hence $\left\vert w\left( .,.\right) \right\vert $ -- is in $C\left( %
\left[ 0,1\right] \times \left[ 0,1\right] \right) $. But $\left[ 0,1\right]
\times \left[ 0,1\right] $ is compact, so that by the Bolzano-Weierstrass
theorem\ there exists a pair $\left( x_{\left\vert w\right\vert
},y_{\left\vert w\right\vert }\right) \in \left[ 0,1\right] \times \left[ 0,1%
\right] $ such that%
\begin{equation*}
\left\vert w\left( x_{_{\left\vert w\right\vert }},y_{_{\left\vert
w\right\vert }}\right) \right\vert =\sup \left\{ \left\vert w\left(
x,y\right) \right\vert :x\in \left[ 0,1\right] ,y\in \left[ 0,1\right]
\right\} .
\end{equation*}%
But $w\left( x,y\right) =\int_{0}^{1}\int_{0}^{1}w\left( u,v\right) \delta
_{x}\left( u\right) \delta _{y}\left( v\right) dudv$ and $\delta _{x},\delta
_{y}$ are in $\mathcal{O}_{M^{1}\left( \left[ 0,1\right] \right) }$ so that 
\begin{eqnarray*}
&&\sup \left\{ \left\vert w\left( x,y\right) \right\vert :x\in \left[ 0,1%
\right] ,y\in \left[ 0,1\right] \right\} \\
&=&\sup \left\{ \left\vert \int_{0}^{1}\int_{0}^{1}w\left( u,v\right) \delta
_{x}\left( u\right) \delta _{y}\left( v\right) dudv\right\vert :x\in \left[
0,1\right] ,y\in \left[ 0,1\right] \right\} \\
&\leq &\left\Vert w\right\Vert _{\vee },
\end{eqnarray*}

Hence,%
\begin{equation*}
\left\vert \sum_{i=1}^{n}f_{i}\left( .\right) g_{i}\left( .\right)
\right\vert _{\infty }\leq \left\Vert \sum_{i=1}^{n}f_{i}\otimes
g_{i}\right\Vert _{\vee }.
\end{equation*}%
This concludes the proof that the imbedding map 
\begin{equation*}
j:\left( C\left( \left[ 0,1\right] \right) \otimes C\left( \left[ 0,1\right]
\right) ,\left\Vert .\right\Vert _{\vee }\right) \hookrightarrow \left(
C\left( \left[ 0,1\right] \times \left[ 0,1\right] \right) ,\left\Vert
.\right\Vert _{\infty }\right)
\end{equation*}%
is continuous and norm-preserving.

To prove that 
\begin{equation*}
\overline{j\left( C\left( \left[ 0,1\right] \right) \otimes C\left( \left[
0,1\right] \right) \right) }=C\left( \left[ 0,1\right] \times \left[ 0,1%
\right] \right) ,
\end{equation*}%
we note that $j\left( C\left( \left[ 0,1\right] \right) \otimes C\left( %
\left[ 0,1\right] \right) \right) $ contains the algebra $A$ of bivariate
polynomials, i.e.%
\begin{equation*}
A:=\left\{ \left. \sum_{i=1}^{n}\sum_{j=1}^{m}\alpha
_{ij}x^{i}y^{j}\right\vert n,m\in \mathbb{N\cup }\left\{ 0\right\}
;a_{ij}\in \mathbb{R}\right\} .
\end{equation*}%
This algebra separates points on $\left[ 0,1\right] \times \left[ 0,1\right] 
$. Hence, we may invoke the Stone-Weierstrass Theorem to deduce that%
\begin{equation*}
\overline{A}=C\left( \left[ 0,1\right] \times \left[ 0,1\right] \right)
\end{equation*}
and since $A$ is contained in $j\left( C\left( \left[ 0,1\right] \right)
\otimes C\left( \left[ 0,1\right] \right) \right) $, that%
\begin{equation*}
C\left( \left[ 0,1\right] \right) \otimes _{\vee }C\left( \left[ 0,1\right]
\right) =\overline{j\left( C\left( \left[ 0,1\right] \right) \otimes C\left( %
\left[ 0,1\right] \right) \right) }=C\left( \left[ 0,1\right] \times \left[
0,1\right] \right) .
\end{equation*}
\end{proof}

\section{Multiplicative level $2$ path lifts: The relevance of an area
process}

\subsection{Multiplicative level $2$ path lifts}

Suppose $x$ is a Banach space $V$-valued path. From Example \ref{Ex: level 1
multiplicative lift} in Chapter 1, we know that $x$ has a canonical
multiplicative lift to $T^{1}\left( V\right) $ given by $\left(
1,x_{t}-x_{s}\right) $. Our aim here is to construct a multiplicative lift
of $x=\mathcal{T}^{B}$ to $T^{2}\left( V\right) $: From (\ref{symmetrisation}%
) we know that every two-tensor can be written as the sum of a symmetric and
an antisymmetric two-tensor. Our method of construction is to choose a
symmetric two-tensor process, $\mathbf{s}^{x^{2}}$, and contingent on the
choice of $\mathbf{s}^{x^{2}}$, an antisymmetric two-tensor process $\mathbf{%
a}^{x^{2}}$, in such a way that 
\begin{equation*}
x^{2}:=\mathbf{s}^{x^{2}}+\mathbf{a}^{x^{2}}
\end{equation*}%
satisfies the multiplicative two-tensor condition, i.e. that for any $s\leq
t\leq u$%
\begin{equation}
x_{s,u}^{2}=x_{s,t}^{2}+x_{t,u}^{2}+x_{s,t}^{1}\otimes x_{t,u}^{1}.
\label{multiplicative two-tensor condition}
\end{equation}%
Any two-tensor process $x^{2}$ that satisfies (\ref{multiplicative
two-tensor condition}) determines a multiplicative functional $\left(
1,x_{t}-x_{s},x_{s,t}^{2}\right) $ in $T^{2}\left( V\right) $. If we choose 
\begin{equation*}
\mathbf{s}_{s,t}^{x^{2}}:=\frac{1}{2}x_{s,t}^{1}\otimes x_{s,t}^{1},
\end{equation*}%
$\mathbf{a}^{x^{2}}$ has to satisfy 
\begin{equation}
\mathbf{a}_{s,u}^{x^{2}}=\mathbf{a}_{s,t}^{x^{2}}+\mathbf{a}_{t,u}^{x^{2}}+%
\frac{1}{2}\left( x_{s,t}^{1}\otimes x_{t,u}^{1}-x_{t,u}^{1}\otimes
x_{s,t}^{1}\right)  \label{area process}
\end{equation}%
for $\mathbf{s}^{x^{2}}+\mathbf{a}^{x^{2}}$ to obey (\ref{multiplicative
two-tensor condition}).

\begin{definition}
Given a $V$-valued path $x$, an anti-symmetric two-tensor process $\mathbf{a}
$ that satisfies (\ref{area process}) for any $s\leq t\leq u$ is called an
area process of the path $x$.
\end{definition}

If we choose $\ \frac{1}{2}x_{s,t}^{1}\otimes x_{s,t}^{1}$ \ as the
symmetric component, then the problem of finding a multiplicative lift of $x$
to $T^{2}\left( V\right) $ reduces to finding an area process $\mathbf{a}$:
Such a process is not unique; different choices of $\mathbf{a}^{x^{2}}$%
define different lifts of $x$.

\subsection{The L\'{e}vy Area}

A particular example of an area is the L\'{e}vy area of a path $x$:

\begin{definition}[L\'{e}vy area random variable ]
For a fixed pair $s\leq t$%
\begin{equation*}
\mathcal{A}_{s,t}\left( x\right) :=\frac{1}{2}\lim_{n\rightarrow \infty
}\sum_{2^{n}s\leq k}^{k\leq 2^{n}t}x_{s,\frac{k}{2^{n}}}^{1}\otimes x_{\frac{%
k}{2^{n}},\frac{k+1}{2^{n}}}^{1}-x_{\frac{k}{2^{n}},\frac{k+1}{2^{n}}%
}^{1}\otimes x_{s,\frac{k}{2^{n}}}^{1}
\end{equation*}%
is called the \emph{L\'{e}vy area of }$x$\emph{\ run from }$s$\emph{\ to }$t$
-- provided the limit exists in $V\otimes _{c}V$.
\end{definition}

\begin{example}
\label{Ex: R^2}Suppose $x=\left( f,g\right) $ is a path in $\mathbb{R}^{2}$
and that $f$ and $g$ are both of bounded variation. In this case, we have 
\begin{equation*}
\mathcal{A}_{s,t}\left( x\right) =\frac{1}{2}\underset{s\leq u\leq v\leq t}{%
\int \int }df_{u}\otimes dg_{v}-dg_{u}\otimes df_{v}
\end{equation*}%
and so if $e_{1}$ and $e_{2}$ represent an orthonormal basis of $\mathbb{R}%
^{2}$ such that 
\begin{equation*}
\left\langle f\otimes g,e_{1}\otimes e_{2}\right\rangle =fg,
\end{equation*}%
then 
\begin{equation*}
\left\langle \mathcal{A}_{s,t}\left( x\right) ,e_{1}\otimes
e_{2}\right\rangle =\frac{1}{2}\underset{s\leq u\leq v\leq t}{\int \int }%
df_{u}dg_{v}-dg_{u}df_{v}.
\end{equation*}%
$\left\langle \mathcal{A}_{s,t}\left( x\right) ,e_{1}\otimes
e_{2}\right\rangle $ has a geometric meaning:\ It is the oriented area
enclosed by the path $x$ run from $s$ to $t$ and the chord between the
points $x\left( s\right) $ and $x\left( t\right) $.
\end{example}

\begin{remark}
Given a path $x$, for each fixed $n$, the dyadic polygonal approximation $%
x\left( n\right) $ of $x$ (c.f. Definition \ref{Defn: dyadic polygonal
approximation} in Chapter 2) has bounded variation and so defines a sequence 
$\mathbf{x}\left( n\right) $ of smooth rough paths (Definition\ \ref{Defn:
smooth rough path}\ , Chapter 1) given by Chen's Theorem (Theorem \ref%
{Theorem: Chen's theorem}, Chapter 1): The two-tensor coordintate $x\left(
n\right) ^{2}$ of $\mathbf{x}\left( n\right) $ is given by 
\begin{equation*}
x\left( n\right) _{s,t}^{2}:=\int_{s}^{t}\int_{s}^{v}x\left( n\right)
_{u}\otimes dx\left( n\right) _{v}
\end{equation*}%
the antisymmetric component of which is 
\begin{eqnarray*}
&&\frac{1}{2}\int_{s}^{t}\left( x\left( n\right) _{u}-x\left( n\right)
_{s}\right) \otimes dx\left( n\right) _{u}-dx\left( n\right) _{u}\otimes
\left( x\left( n\right) _{u}-x\left( n\right) _{s}\right) \\
&=&\frac{1}{2}\sum_{2^{n}s\leq k}^{k\leq 2^{n}t-1}x_{s,\frac{k}{2^{n}}%
}^{1}\otimes x_{\frac{k}{2^{n}},\frac{k+1}{2^{n}}}^{1}-x_{\frac{k}{2^{n}},%
\frac{k+1}{2^{n}}}^{1}\otimes x_{s,\frac{k}{2^{n}}}^{1}.
\end{eqnarray*}%
Thus, if $x$ has a L\'{e}vy area, we get a second representation of the L%
\'{e}vy area as 
\begin{equation*}
\mathcal{A}_{s,t}\left( x\right) =\frac{1}{2}\lim_{n\rightarrow \infty }%
\underset{s\leq u\leq v\leq t}{\int \int }dx\left( n\right) _{u}\otimes
dx\left( n\right) _{v}-dx\left( n\right) _{v}\otimes dx\left( n\right) _{u}.
\end{equation*}
\end{remark}

\section{Non-Existence of the L\'{e}vy Area of $\mathcal{T}^{B}$ in the
injective tensor algebra}

\subsection{Introduction\label{Section: area introduction}}

Let us assume for the moment that the L\'{e}vy Area $\mathcal{A}_{0,1}\left( 
\mathcal{T}^{B}\right) $ of the frame process $\mathcal{T}^{B}$ exists in
the injective tensor algebra. By Proposition \ref{Proposition: injective
tensor algebra}, this is equivalent to the existence of a continuous process 
$\mathcal{A}_{0,1}\left( \mathcal{T}^{B}\right) \left( .,.\right) $ on the
unit square, i.e. the antisymmetric two-tensor $\mathcal{A}_{0,1}\left( 
\mathcal{T}^{B}\right) $ has to be a $C\left( \left[ 0,1\right] \times \left[
0,1\right] \right) $-valued random-variable. By Definition \ref{Defn:
symmetric + antisymmetric two-tensor}, for any $\left( s,t\right) \in \left[
0,1\right] \times \left[ 0,1\right] $ and $\mathbb{P}$-a.e. sample path $%
f\in C_{0}\left( \left[ -1,1\right] \right) $, we have that%
\begin{eqnarray}
&&\mathcal{A}_{0,1}\left( \mathcal{T}^{f}\right) \left( s,t\right)  \notag \\
&=&\left\langle \mathcal{A}_{0,1}\left( \mathcal{T}^{f}\right) ,\delta
_{s}\otimes \delta _{t}\right\rangle  \notag \\
&=&-\left\langle \mathcal{A}_{0,1}\left( \mathcal{T}^{f}\right) ,\delta
_{t}\otimes \delta _{s}\right\rangle  \notag \\
&=&-\mathcal{A}_{0,1}\left( \mathcal{T}^{f}\right) \left( t,s\right) .
\label{Eqn: antisymmetry}
\end{eqnarray}%
Hence, since we assumed that $\mathcal{A}_{0,1}\left( \mathcal{T}^{f}\right) 
$ is continuous on $\left[ 0,1\right] \times \left[ 0,1\right] $, we require
that for any sequence $s_{n}\rightarrow s$ through $\left[ 0,1\right] ,$ 
\begin{equation}
\mathcal{A}_{0,1}\left( \mathcal{T}^{f}\right) \left( s_{n},s\right)
\rightarrow 0.  \label{continuity condition}
\end{equation}%
However, we show that for any fixed $s\in \left[ 0,1\right] $ and any
sequence $s_{n}\nearrow s$, for $\mathbb{P}$-a.e. sample path $f\in
C_{0}\left( \left[ -1,1\right] \right) $, 
\begin{equation*}
\mathcal{A}_{0,1}\left( \mathcal{T}^{f}\right) \left( s_{n},s\right)
\rightarrow -\frac{1}{2},
\end{equation*}%
thus contradicting (\ref{continuity condition}) and hence disproving our
initial assumption that $\mathcal{A}_{0,1}\left( \mathcal{T}^{B}\right) $
exists in the injective tensor product.

\subsection{Existence of the $\mathbb{R}$-valued random variable $\mathcal{A}%
_{0,1}\left( \mathcal{T}^{B}\right) \left( s,t\right) $ for dyadic $s$ and $%
t $}

We show the existence of $\mathcal{A}_{0,1}\left( \mathcal{T}^{B}\right) $
on the dyadic unit square $D\left( \left[ 0,1\right] ^{2}\right) $. By (\ref%
{Eqn: antisymmetry}), if we have defined $\mathcal{A}_{0,1}\left( \mathcal{T}%
^{B}\right) $ on $\left\{ \left( s,t\right) :0\leq s<t\leq 1\right\} $, we
have defined it on the entire unit square ($\mathcal{A}_{0,1}\left( \mathcal{%
T}^{B}\right) \left( s,s\right) :=0$ for $s\in \left[ 0,1\right] $).

\begin{proposition}
\label{Propn: Representation of Levy Area random variable}

For $\mathbb{P}$-a.e. $f\in C_{0}\left( \left[ -1,1\right] \right) $ and $%
\left( s,t\right) \in D\left( \left[ 0,1\right] ^{2}\right) $, the limit 
\begin{equation*}
\lim_{n\rightarrow \infty }\mathcal{A}_{0,1}\left( \mathcal{T}_{.}^{B}\left(
n\right) \right) \left( s,t\right)
\end{equation*}%
exists. Moreover, on the open upper dyadic triangle 
\begin{equation*}
\left\{ \left( s,t\right) \in D\left( \left[ 0,1\right] ^{2}\right) :0\leq
s<t\leq 1\right\}
\end{equation*}%
(which \emph{excludes }the diagonal), 
\begin{equation*}
\mathcal{A}_{0,1}\left( \mathcal{T}^{B}\right) :=\lim_{n\rightarrow \infty }%
\mathcal{A}_{0,1}\left( \mathcal{T}_{.}^{B}\left( n\right) \right)
\end{equation*}%
has the following representation:%
\begin{eqnarray*}
&&\mathcal{A}_{0,1}\left( \mathcal{T}^{B}\right) \left( s,t\right) \\
&=&\frac{1}{2}\int_{t-1}^{t}B_{v-\left( t-s\right) }dB_{v}-\frac{1}{2}%
\int_{t-1}^{s}B_{v}dB_{v}+\frac{1}{2}B_{t-1}\left( B_{s}-B_{t-1}\right) -%
\frac{1}{2}B_{s-1}\left( B_{t}-B_{t-1}\right) \\
&&-\frac{1}{2}\left( 1-t+s\right) \\
&&+\frac{1}{2}\int_{1-s}^{2-s}\hat{B}_{v-\left( t-s\right) }d\hat{B}_{v}-%
\frac{1}{2}\int_{1-s}^{2-t}\hat{B}_{v}d\hat{B}_{v}-\frac{1}{2}\hat{B}%
_{2-t}\left( \hat{B}_{2-s}-\hat{B}_{2-t}\right) .
\end{eqnarray*}%
The integrals are understood in an It\^{o} sense and the Brownian motion $%
\hat{B}$ is defined as 
\begin{equation*}
\hat{B}_{u}:=B_{1}-B_{1-u},\ \ u\in \left[ 0,2\right] .
\end{equation*}%
By (\ref{Eqn: antisymmetry}), this defines $\mathcal{A}_{0,1}\left( \mathcal{%
T}^{B}\right) $ on the open lower dyadic triangle, 
\begin{equation*}
\left\{ \left( s,t\right) \in D\left( \left[ 0,1\right] ^{2}\right) :0\leq
t<s\leq 1\right\} ,
\end{equation*}%
via 
\begin{equation*}
\mathcal{A}_{0,1}\left( \mathcal{T}^{B}\right) \left( s,t\right) :=-\mathcal{%
A}_{0,1}\left( \mathcal{T}^{B}\right) \left( t,s\right) .
\end{equation*}%
On the dyadic diagonal $\left\{ \left( s,s\right) :s\in D\left( \left[ 0,1%
\right] \right) \right\} $,%
\begin{equation*}
\mathcal{A}_{0,1}\left( \mathcal{T}^{B}\right) \left( s,s\right) :=0.
\end{equation*}
\end{proposition}

\begin{proof}
We fix $n$ so large that $\frac{1}{2^{n}}\leq t-s$. Then 
\begin{eqnarray}
&&\mathcal{A}_{0,1}\left( \mathcal{T}^{B}\left( n\right) \right) \left(
s,t\right)  \notag \\
&=&\frac{1}{2}\sum_{v=1}^{2^{n}-1}\left( \mathcal{T}_{\frac{v}{2^{n}}%
}^{B}\left( s\right) -\mathcal{T}_{0}^{B}\left( s\right) \right) \left( 
\mathcal{T}_{\frac{v+1}{2^{n}}}^{B}\left( t\right) -\mathcal{T}_{\frac{v}{%
2^{n}}}^{B}\left( t\right) \right)  \notag \\
&&-\frac{1}{2}\sum_{v=1}^{2^{n}-1}\left( \mathcal{T}_{\frac{v}{2^{n}}%
}^{B}\left( t\right) -\mathcal{T}_{0}^{B}\left( t\right) \right) \left( 
\mathcal{T}_{\frac{v+1}{2^{n}}}^{B}\left( s\right) -\mathcal{T}_{\frac{v}{%
2^{n}}}^{B}\left( s\right) \right)  \notag \\
&=&\frac{1}{2}\sum_{v=1}^{2^{n}-1}\sum_{u=0}^{v-1}\left( B_{\frac{u+1}{2^{n}}%
-1+s}-B_{\frac{u}{2^{n}}-1+s}\right) \left( B_{\frac{v+1}{2^{n}}-1+t}-B_{%
\frac{v}{2^{n}}-1+t}\right)  \notag \\
&&-\frac{1}{2}\sum_{v=1}^{2^{n}-1}\sum_{u=0}^{v-1}\left( B_{\frac{u+1}{2^{n}}%
-1+t}-B_{\frac{u}{2^{n}}-1+t}\right) \left( B_{\frac{v+1}{2^{n}}-1+s}-B_{%
\frac{v}{2^{n}}-1+s}\right) .  \label{Eqn: first area representation}
\end{eqnarray}

\begin{figure}
\centering
\includegraphics[scale=0.5]{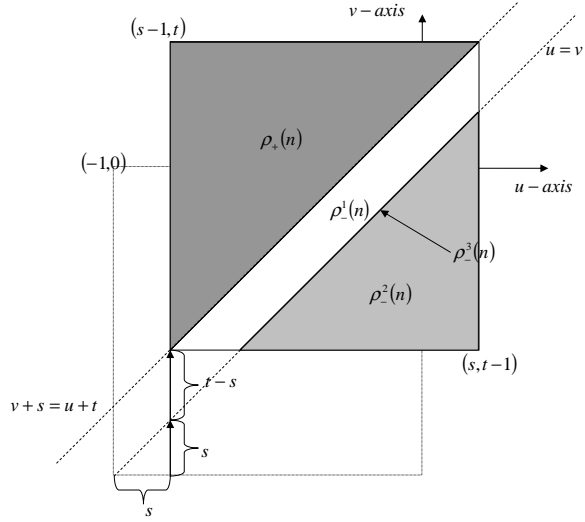}
\caption{Region of integration for $s<t$} \label{myfig}
\end{figure}

Figure (\ref{myfig}) illustrates that summing terms of the form

\begin{eqnarray*}
&\left( B_{\frac{u+1}{2^{n}}-1+s}-B_{\frac{u}{2^{n}}-1+s}\right) \left( B_{%
\frac{v+1}{2^{n}}-1+t}-B_{\frac{v}{2^{n}}-1+t}\right) \\
&-\left( B_{\frac{u+1}{2^{n}}-1+t}-B_{\frac{u}{2^{n}}-1+t}\right) \left( B_{%
\frac{v+1}{2^{n}}-1+s}-B_{\frac{v}{2^{n}}-1+s}\right)
\end{eqnarray*}%
over points $(\frac{u}{2^{n}},\frac{v}{2^{n}})$ in the upper dyadic triangle
of the unit square as in (\ref{Eqn: first area representation}), is
equivalent to summing terms of the form 
\begin{equation*}
\left( B_{\frac{u+1}{2^{n}}}-B_{\frac{u}{2^{n}}}\right) \left( B_{\frac{v+1}{%
2^{n}}}-B_{\frac{v}{2^{n}}}\right)
\end{equation*}%
over points $(\frac{u}{2^{n}},\frac{v}{2^{n}})$ in the upper dyadic triangle
of the unit square shifted by $\left( s-1,t-1\right) $. This set of shifted
points can be written as the (disjoint) union of 
\begin{eqnarray*}
\rho _{+}\left( n\right) &=&\left\{ \left( \frac{u}{2^{n}},\frac{v}{2^{n}}%
\right) :u,v\in \left\{ 0,1,2,...,2^{n}-1\right\} \text{ and}\right. \\
&&\left. \frac{u+1}{2^{n}}+t\leq \frac{v}{2^{n}}+s\text{ and }\frac{u}{2^{n}}%
\geq s-1\text{ and }\right. \\
&&\left. t-1\leq \frac{v}{2^{n}}\leq t-\frac{1}{2^{n}}\right\} ,
\end{eqnarray*}%
\begin{eqnarray*}
\rho _{-}^{1}\left( n\right) &=&\left\{ \left( \frac{u}{2^{n}},\frac{v}{2^{n}%
}\right) :u,v\in \left\{ 0,1,2,...,2^{n}-1\right\} \text{ and}\right. \\
&&\left. \frac{u}{2^{n}}+t\geq \frac{v}{2^{n}}+s\text{ and }\frac{v}{2^{n}}%
\geq \left( t-1\right) \vee \frac{u+1}{2^{n}}\right. \\
&&\left. \text{and }s-1\leq \frac{u}{2^{n}}\leq s-\frac{1}{2^{n}}\right\} ,
\end{eqnarray*}%
\begin{eqnarray*}
\rho _{-}^{2}\left( n\right) &=&\left\{ \left( \frac{u}{2^{n}},\frac{v}{2^{n}%
}\right) :u,v\in \left\{ 0,1,2,...,2^{n}-1\right\} \text{ and}\right. \\
&&\left. t-1\leq \frac{v}{2^{n}}\leq \frac{u-1}{2^{n}}\text{ and }t-1\leq 
\frac{u}{2^{n}}\leq s-\frac{1}{2^{n}}\right\}
\end{eqnarray*}%
and%
\begin{eqnarray*}
\rho _{-}^{3}\left( n\right) &=&\left\{ \left( \frac{u}{2^{n}},\frac{v}{2^{n}%
}\right) :u,v\in \left\{ 0,1,2,...,2^{n}-1\right\} \text{ and}\right. \\
&&\left. u=v\text{ and }t-1\leq \frac{u}{2^{n}}\leq s-\frac{1}{2^{n}}%
\right\} .
\end{eqnarray*}

\ In this way, we can rewrite $\mathcal{A}_{0,1}\left( \mathcal{T}^{B}\left(
n\right) \right) \left( s,t\right) $ as 
\begin{eqnarray*}
&&\mathcal{A}_{0,1}\left( \mathcal{T}^{B}\left( n\right) \right) \left(
s,t\right) \\
&=&\frac{1}{2}\sum_{\left( \frac{u}{2^{n}},\frac{v}{2^{n}}\right) \in \rho
_{+}\left( n\right) }\left( B_{\frac{u+1}{2^{n}}}-B_{\frac{u}{2^{n}}}\right)
\left( B_{\frac{v+1}{2^{n}}}-B_{\frac{v}{2^{n}}}\right) \\
&&-\frac{1}{2}\sum_{\left( \frac{u}{2^{n}},\frac{v}{2^{n}}\right) \in \rho
_{-}^{1}\left( n\right) }\left( B_{\frac{u+1}{2^{n}}}-B_{\frac{u}{2^{n}}%
}\right) \left( B_{\frac{v+1}{2^{n}}}-B_{\frac{v}{2^{n}}}\right) \\
&&-\frac{1}{2}\sum_{\left( \frac{u}{2^{n}},\frac{v}{2^{n}}\right) \in \rho
_{-}^{2}\left( n\right) }\left( B_{\frac{u+1}{2^{n}}}-B_{\frac{u}{2^{n}}%
}\right) \left( B_{\frac{v+1}{2^{n}}}-B_{\frac{v}{2^{n}}}\right) \\
&&-\frac{1}{2}\sum_{\left( \frac{u}{2^{n}},\frac{v}{2^{n}}\right) \in \rho
_{-}^{3}\left( n\right) }\left( B_{\frac{u+1}{2^{n}}}-B_{\frac{u}{2^{n}}%
}\right) \left( B_{\frac{v+1}{2^{n}}}-B_{\frac{v}{2^{n}}}\right) .
\end{eqnarray*}%
Since $t>s$, the process $B_{.-\left( t-s\right) }$ is adapted to the
Brownian filtration 
\begin{equation*}
\left\{ \sigma \left( B_{u}:-1\leq u\leq t\right) :t\in \left[ -1,1\right]
\right\}
\end{equation*}%
and so 
\begin{eqnarray*}
&&\sum_{\left( \frac{u}{2^{n}},\frac{v}{2^{n}}\right) \in \rho _{+}\left(
n\right) }\left( B_{\frac{u+1}{2^{n}}}-B_{\frac{u}{2^{n}}}\right) \left( B_{%
\frac{v+1}{2^{n}}}-B_{\frac{v}{2^{n}}}\right) \\
&=&\sum_{v=2^{n}\left( t-1\right) }^{2^{n}t-1}\left( B_{\frac{v}{2^{n}}%
-\left( t-s\right) }-B_{s-1}\right) \left( B_{\frac{v+1}{2^{n}}}-B_{\frac{v}{%
2^{n}}}\right)
\end{eqnarray*}%
converges to the It\^{o} integral $\int_{t-1}^{t}\left( B_{v-\left(
t-s\right) }-B_{s-1}\right) dB_{v}$ as $n\rightarrow \infty $, $\mathbb{P}$%
-a.s.

Similarly,%
\begin{eqnarray*}
&&\sum_{\left( \frac{u}{2^{n}},\frac{v}{2^{n}}\right) \in \rho
_{-}^{2}\left( n\right) }\left( B_{\frac{u+1}{2^{n}}}-B_{\frac{u}{2^{n}}%
}\right) \left( B_{\frac{v+1}{2^{n}}}-B_{\frac{v}{2^{n}}}\right) \\
&=&\sum_{v=2^{n}\left( t-1\right) }^{2^{n}s-1}\left( B_{\frac{v}{2^{n}}%
}-B_{t-1}\right) \left( B_{\frac{v+1}{2^{n}}}-B_{\frac{v}{2^{n}}}\right)
\end{eqnarray*}%
converges to the It\^{o} integral $\int_{t-1}^{s}\left( B_{v}-B_{t-1}\right)
dB_{v}$, $\mathbb{P}$-a.s.

Furthermore, for \ the sum over $\rho _{-}^{3}\left( n\right) $, we get that 
\begin{eqnarray*}
&&\sum_{\left( \frac{u}{2^{n}},\frac{v}{2^{n}}\right) \in \rho
_{-}^{3}\left( n\right) }\left( B_{\frac{u+1}{2^{n}}}-B_{\frac{u}{2^{n}}%
}\right) \left( B_{\frac{v+1}{2^{n}}}-B_{\frac{v}{2^{n}}}\right) \\
&=&\sum_{u=\left( t-1\right) 2^{n}}^{2^{n}s-1}\left( B_{\frac{u+1}{2^{n}}%
}-B_{\frac{u}{2^{n}}}\right) ^{2}
\end{eqnarray*}%
converges to the quadratic variation $\left\langle B,B\right\rangle
_{t-1}^{s}=1-t+s$ , $\mathbb{P}$-a.s.

Finally, for the sum over $\rho _{-}^{1}\left( n\right) ,$ we find that

\begin{eqnarray*}
&&\sum_{\left( \frac{u}{2^{n}},\frac{v}{2^{n}}\right) \in \rho
_{-}^{1}\left( n\right) }\left( B_{\frac{u+1}{2^{n}}}-B_{\frac{u}{2^{n}}%
}\right) \left( B_{\frac{v+1}{2^{n}}}-B_{\frac{v}{2^{n}}}\right) \\
&=&\sum_{v=2^{n}\left( s-1\right) }^{2^{n}s-1}\left( B_{\frac{v+1}{2^{n}}%
+t-s}-B_{\left( t-1\right) \vee \frac{v+1}{2^{n}}}\right) \left( B_{\frac{v+1%
}{2^{n}}}-B_{\frac{v}{2^{n}}}\right) \\
&=&\sum_{k=2^{n}\left( 1-s\right) +1}^{2^{n}\left( 2-s\right) }\left( B_{1-%
\frac{k-1}{2^{n}}+\left( t-s\right) }-B_{\left( t-1\right) \vee 1-\frac{k-1}{%
2^{n}}}\right) \left( B_{1-\frac{k-1}{2^{n}}}-B_{1-\frac{k}{2^{n}}}\right) \\
&=&\sum_{k=2^{n}\left( 1-s\right) +1}^{2^{n}\left( 2-s\right) }\left[ \left(
B_{1}-B_{\left( t-1\right) \vee 1-\frac{k-1}{2^{n}}}\right) -\left(
B_{1}-B_{1-\frac{k-1}{2^{n}}+\left( t-s\right) }\right) \right] \\
&&\times \left[ \left( B_{1}-B_{1-\frac{k}{2^{n}}}\right) -\left( B_{1}-B_{1-%
\frac{k-1}{2^{n}}}\right) \right] \\
&=&\sum_{k=2^{n}\left( 1-s\right) +1}^{2^{n}\left( 2-s\right) }\left( \hat{B}%
_{2-t\wedge \frac{k-1}{2^{n}}}-\hat{B}_{\frac{k-1}{2^{n}}-\left( t-s\right)
}\right) \left( \hat{B}_{\frac{k}{2^{n}}}-\hat{B}_{\frac{k-1}{2^{n}}}\right)
\end{eqnarray*}%
where $\hat{B}_{u}:=B_{1}-B_{1-u}$ is a Brownian motion with respect to the
filtration 
\begin{equation*}
\left\{ \sigma \left( \hat{B}_{u}:0\leq u\leq t\right) :t\in \left[ 0,2%
\right] \right\} .
\end{equation*}%
Taking the limit as $n\rightarrow \infty $, we find that%
\begin{eqnarray*}
&&\sum_{k=2^{n}\left( 1-s\right) +1}^{2^{n}\left( 2-s\right) }\left( \hat{B}%
_{2-t\wedge \frac{k-1}{2^{n}}}-\hat{B}_{\frac{k-1}{2^{n}}-\left( t-s\right)
}\right) \left( \hat{B}_{\frac{k}{2^{n}}}-\hat{B}_{\frac{k-1}{2^{n}}}\right)
\\
&\rightarrow &\int_{1-s}^{2-t}\hat{B}_{u}d\hat{B}_{u}+\hat{B}_{2-t}\left( 
\hat{B}_{2-s}-\hat{B}_{2-t}\right) -\int_{1-s}^{2-s}\hat{B}_{u-\left(
t-s\right) }d\hat{B}_{u}
\end{eqnarray*}

After some tidying up, we find that for dyadic times $s<t$, 
\begin{eqnarray*}
&&\mathcal{A}_{0,1}\left( \mathcal{T}^{B}\right) \left( s,t\right) \\
&=&\frac{1}{2}\int_{t-1}^{t}B_{v-\left( t-s\right) }dB_{v}-\frac{1}{2}%
B_{s-1}\left( B_{t}-B_{t-1}\right) -\frac{1}{2}\int_{t-1}^{s}B_{v}dB_{v}+%
\frac{1}{2}B_{t-1}\left( B_{s}-B_{t-1}\right) \\
&&-\frac{1}{2}\left( 1-t+s\right) \\
&&+\frac{1}{2}\int_{1-s}^{2-s}\hat{B}_{v-\left( t-s\right) }d\hat{B}_{v}-%
\frac{1}{2}\int_{1-s}^{2-t}\hat{B}_{v}d\hat{B}_{v}-\frac{1}{2}\hat{B}%
_{2-t}\left( \hat{B}_{2-s}-\hat{B}_{2-t}\right) .
\end{eqnarray*}
\end{proof}

\subsection{Non-existence of the L\'{e}vy area random variable $\mathcal{A}%
_{0,1}\left( \mathcal{T}^{B}\right) $ in the injective tensor product}

In Proposition \ref{Propn: Representation of Levy Area random variable}, we
showed that $\mathcal{A}_{0,1}\left( \mathcal{T}^{B}\right) $ exists on the
dyadic unit square $D\left( \left[ 0,1\right] ^{2}\right) $. In this
section, we prove that $\mathcal{A}_{0,1}\left( \mathcal{T}^{B}\right) $ has
a unique continuous extension to the entire unit square \emph{excluding the
diagonal}. We show that this extension has a jump discontinuity on the
diagonal and hence, that the L\'{e}vy area does not exist in the injective\
tensor product.

\begin{proposition}
\label{Propn: discontinuity of LA}With $\mathbb{P}$-probability $1$, the L%
\'{e}vy Area $\mathcal{A}_{0,1}\left( \mathcal{T}^{B}\right) $ as defined in
Proposition \ref{Propn: Representation of Levy Area random variable} on the
dyadic open upper triangle 
\begin{equation*}
\left\{ \left( s,t\right) \in D\left( \left[ 0,1\right] ^{2}\right) :0\leq
s<t\leq 1\right\}
\end{equation*}%
has a unique uniformly continuous extension $\mathcal{A}_{0,1}^{\ast }\left( 
\mathcal{T}^{B}\right) $ to the open upper triangle 
\begin{equation*}
\left\{ \left( s,t\right) \in \left[ 0,1\right] \times \left[ 0,1\right]
:0\leq s<t\leq 1\right\}
\end{equation*}%
and hence by (\ref{Eqn: antisymmetry}) to the open lower triangle of the
unit square 
\begin{equation*}
\left\{ \left( s,t\right) \in \left[ 0,1\right] \times \left[ 0,1\right]
:0\leq t<s\leq 1\right\} .
\end{equation*}%
$\mathcal{A}_{0,1}^{\ast }\left( \mathcal{T}^{B}\right) $ agrees with $%
\mathcal{A}_{0,1}\left( \mathcal{T}^{B}\right) $ on the dyadic open upper
triangle (dyadic open lower triangle) on a set $\Omega _{0}$ of full $%
\mathbb{P}$-measure.

$\mathcal{A}_{0,1}^{\ast }\left( \mathcal{T}^{B}\right) $ has a jump
discontinuity on the diagonal and so does not exist in the injective tensor
algebra $C\left( \left[ 0,1\right] \times \left[ 0,1\right] \right) $ (c.f.
Proposition \ref{Proposition: injective tensor algebra}): There does not
exist a continuous extension of $\mathcal{A}_{0,1}$ to the unit square.
\end{proposition}

\begin{proof}
We start by proving the existence of a continuous extension of $\mathcal{A}%
_{0,1}\left( \mathcal{T}^{B}\right) $ to the open upper triangle:\ In
Proposition \ref{Propn: Representation of Levy Area random variable}, we
showed that on the \emph{open }upper dyadic triangle, 
\begin{equation*}
\left\{ \left( s,t\right) \in D\left( \left[ 0,1\right] ^{2}\right) :0\leq
s<t\leq 1\right\} ,
\end{equation*}%
$\mathcal{A}_{0,1}\left( \mathcal{T}^{B}\right) \left( s,t\right) $ has the
following representation: 
\begin{eqnarray}
&&\mathcal{A}_{0,1}\left( \mathcal{T}^{B}\right) \left( s,t\right)  \notag \\
&=&\frac{1}{2}\int_{t-1}^{t}B_{v-\left( t-s\right) }dB_{v}-\frac{1}{2}%
\int_{t-1}^{s}B_{v}dB_{v}+\frac{1}{2}B_{t-1}\left( B_{s}-B_{t-1}\right) -%
\frac{1}{2}B_{s-1}\left( B_{t}-B_{t-1}\right)  \notag \\
&&-\frac{1}{2}\left( 1-t+s\right)  \notag \\
&&+\frac{1}{2}\int_{1-s}^{2-s}\hat{B}_{v-\left( t-s\right) }d\hat{B}_{v}-%
\frac{1}{2}\int_{1-s}^{2-t}\hat{B}_{v}d\hat{B}_{v}-\frac{1}{2}\hat{B}%
_{2-t}\left( \hat{B}_{2-s}-\hat{B}_{2-t}\right) .
\label{eqn: LA representation}
\end{eqnarray}%
For a Brownian Motion $W$, the It\^{o} integral $M_{b}:=%
\int_{const}^{b}W_{u}dW_{u}$ is a continuous martingale with respect to the
filtration $\left\{ \sigma \left( W_{u}:const\leq u\leq b\right) :b\in \left[
const,T\right] \right\} $, so that $M\left( a,b\right) :=M\left( b\right)
-M\left( a\right) $ is continuous on a set of full $\mathbb{P}$-measure.
Furthermore, the closed upper triangle 
\begin{equation*}
\left\{ \left( s,t\right) \in \left[ 0,1\right] \times \left[ 0,1\right]
:0\leq s\leq t\leq 1\right\}
\end{equation*}%
is compact. Hence, 
\begin{equation*}
F_{1}\left( s,t\right) :=\frac{1}{2}\int_{t-1}^{s}B_{v}dB_{v}
\end{equation*}%
and 
\begin{equation*}
\hat{F}_{1}\left( s,t\right) :=\frac{1}{2}\int_{1-s}^{2-t}\hat{B}_{v}d\hat{B}%
_{v}
\end{equation*}%
have uniformly continuous sample paths on the closed upper triangle. $%
B_{t-1}\left( B_{s}-B_{t-1}\right) $, $B_{s-1}\left( B_{t}-B_{t-1}\right) $
and $\hat{B}_{2-t}\left( \hat{B}_{2-s}-\hat{B}_{2-t}\right) $ are uniformly
continuous in $s$ and $t$ \ on the closed upper triangle because of the $%
\mathbb{P}$-a.s. sample path continuity of Brownian motion and compactness
of the closed upper triangle. This leaves us with having to prove that the
terms 
\begin{equation*}
F_{2}\left( s,t\right) :=\int_{t-1}^{t}B_{v-\left( t-s\right) }dB_{v}
\end{equation*}%
and%
\begin{equation*}
\hat{F}_{2}\left( s,t\right) :=\int_{1-s}^{2-s}\hat{B}_{u-\left( t-s\right)
}d\hat{B}_{u}
\end{equation*}%
have uniformly continuous modifications on the closed dyadic upper triangle.
Since 
\begin{equation*}
\left\{ \left( s,t\right) \in D\left( \left[ 0,1\right] ^{2}\right) :0\leq
s\leq t\leq 1\right\}
\end{equation*}%
is dense in 
\begin{equation*}
\left\{ \left( s,t\right) \in \left[ 0,1\right] \times \left[ 0,1\right]
:0\leq s\leq t\leq 1\right\} ,
\end{equation*}%
this is equivalent to showing that $F_{2}$ and $\hat{F}_{2}$ each have a
unique uniformly continuous extension to the closed upper triangle.

We start by establishing the existence of a continuous modification of $F_{2}
$ by using Kolmogorov's Lemma (Theorem \ref{Thm: Kolmogorov's Lemma} in
Chapter $2$) for a two-dimensional parameter set. Suppose that $\left(
s_{1},t_{1}\right) $ and $\left( s_{2},t_{2}\right) $ are in $\left\{ \left(
s,t\right) \in D\left( \left[ 0,1\right] ^{2}\right) :0\leq s<t\leq
1\right\} $ and w.l.o.g. suppose that $t_{2}\geq t_{1}$. For $m\in \mathbb{N}
$, we have that 
\begin{eqnarray*}
&&\left\vert F_{2}\left( s_{2},t_{2}\right) -F_{2}\left( s_{1},t_{1}\right)
\right\vert ^{2m} \\
&=&\left\vert F_{2}\left( s_{2},t_{2}\right) -F_{2}\left( s_{1},t_{2}\right)
+F_{2}\left( s_{1},t_{2}\right) -F_{2}\left( s_{1},t_{1}\right) \right\vert
^{2m} \\
&\leq &4^{2m-1}\left\{ \left\vert \int_{t_{2}-1}^{t_{2}}\left( B_{u-\left(
t_{2}-s_{2}\right) }-B_{u-\left( t_{2}-s_{1}\right) }\right)
dB_{u}\right\vert ^{2m}+\left\vert \int_{t_{1}}^{t_{2}}B_{u-\left(
t_{2}-s_{1}\right) }dB_{u}\right\vert ^{2m}\right.  \\
&&\left. +\left\vert \int_{t_{2}-1}^{t_{1}}\left( B_{u-\left(
t_{2}-s_{1}\right) }-B_{u-\left( t_{1}-s_{1}\right) }\right)
dB_{u}\right\vert ^{2m}+\left\vert \int_{t_{1}-1}^{t_{2}-1}B_{u-\left(
t_{1}-s_{1}\right) }dB_{u}\right\vert ^{2m}\right\} 
\end{eqnarray*}%
by Jensen's inequality. We require the following result from \cite{Shreve91}%
, equations (3.24) and (3.25): If $X$ is a measurable, adapted process
satisfying 
\begin{equation*}
\mathbb{E}\left[ \int_{0}^{T}X_{t}^{2m}dt\right] <\infty ,
\end{equation*}%
for some real number $T>0$ and $m\in \mathbb{N}$ and $W_{.}$ is a standard,
one-dimensional Brownian Motion on $\left[ 0,T\right] $, then 
\begin{equation*}
\mathbb{E}\left[ \left\vert \int_{0}^{T}X_{u}dW_{u}\right\vert ^{2m}\right]
\leq \left( m\left( 2m-1\right) \right) ^{m}T^{m-1}\mathbb{E}\left[
\int_{0}^{T}\left\vert X_{u}\right\vert ^{2m}du\right] .
\end{equation*}%
We apply this to our setup where the Brownian motion $B$ is defined on $%
\left[ -1,1\right] $ ( instead of on $\left[ 0,T\right] $ ) to find the
following moment bounds: We find that 
\begin{eqnarray*}
&&\mathbb{E}\left[ \left\vert \int_{t_{2}-1}^{t_{2}}\left( B_{u-\left(
t_{2}-s_{2}\right) }-B_{u-\left( t_{2}-s_{1}\right) }\right)
dB_{u}\right\vert ^{2m}\right]  \\
&\leq &\left( m\left( 2m-1\right) \right) ^{m}\frac{2^{m}}{\sqrt{\pi }}%
\Gamma \left( \frac{2m+1}{2}\right) \left\vert s_{2}-s_{1}\right\vert ^{m},
\end{eqnarray*}%
where we used (\ref{Gaussian moments}) in Chapter 2. Similarly, we find that 
\begin{eqnarray*}
&&\mathbb{E}\left[ \left\vert \int_{t_{1}}^{t_{2}}B_{u-\left(
t_{2}-s_{1}\right) }dB_{u}\right\vert ^{2m}\right]  \\
&\leq &\left( m\left( 2m-1\right) \right) ^{m}\frac{2^{m}}{\sqrt{\pi }}%
\Gamma \left( \frac{2m+1}{2}\right) \left\vert t_{2}-t_{1}\right\vert
^{m-1}\int_{t_{1}}^{t_{2}}\left\vert u-t_{2}+s_{1}\right\vert ^{m}du \\
&\leq &\left( m\left( 2m-1\right) \right) ^{m}\frac{2^{m}}{\sqrt{\pi }}%
\Gamma \left( \frac{2m+1}{2}\right) \left\vert t_{2}-t_{1}\right\vert ^{m},
\end{eqnarray*}%
and%
\begin{eqnarray*}
&&\mathbb{E}\left[ \left\vert \int_{t_{2}-1}^{t_{1}}\left( B_{u-\left(
t_{2}-s_{1}\right) }-B_{u-\left( t_{1}-s_{1}\right) }\right)
dB_{u}\right\vert ^{2m}\right]  \\
&\leq &\left( m\left( 2m-1\right) \right) ^{m}\frac{2^{m}}{\sqrt{\pi }}%
\Gamma \left( \frac{2m+1}{2}\right) \left\vert t_{2}-t_{1}\right\vert ^{m},
\end{eqnarray*}%
and finally that 
\begin{eqnarray*}
&&\mathbb{E}\left[ \left\vert \int_{t_{1}-1}^{t_{2}-1}B_{u-\left(
t_{1}-s_{1}\right) }dB_{u}\right\vert ^{2m}\right]  \\
&\leq &\left( m\left( 2m-1\right) \right) ^{m}\frac{2^{m}}{\sqrt{\pi }}%
\Gamma \left( \frac{2m+1}{2}\right) \left\vert t_{2}-t_{1}\right\vert ^{m}.
\end{eqnarray*}%
Therefore, for $m=3$, we have the bound

\begin{equation*}
\mathbb{E}\left[ \left\vert F_{2}\left( s_{2},t_{2}\right) -F_{2}\left(
s_{1},t_{1}\right) \right\vert ^{6}\right] \leq \frac{15^{3}.2^{15}}{\sqrt{%
\pi }}\Gamma \left( \frac{7}{2}\right) \left( \max \left( \left\vert
s_{2}-s_{1}\right\vert ,\left\vert t_{2}-t_{1}\right\vert \right) \right)
^{3}.
\end{equation*}

Hence, by Kolmogorov's Theorem (c.f. Theorem \ref{Thm: Kolmogorov's Lemma}\
in Chapter 2), $F_{2}\left( .,.\right) $ has a continuous modification $%
F_{2}^{\ast }$ on the closed upper dyadic triangle 
\begin{equation*}
\left\{ \left( s,t\right) \in D\left( \left[ 0,1\right] ^{2}\right) :0\leq
s\leq t\leq 1\right\} .
\end{equation*}

Since the closed upper dyadic triangle is dense in the closed upper triangle
\begin{equation*} 
\left\{ \left( s,t\right) :0\leq s\leq t\leq 1\right\}, 
\end{equation*}
$F_{2}^{\ast }$ has a unique uniformly continuous extension $G_{2}^{\ast }$ to the closed
upper triangle.

The argument for $\hat{F}_{2}$ is similar -- its continuous extension to the
closed upper triangle is denoted as $\hat{G}_{2}^{\ast }$.

In summary, 
\begin{eqnarray}
&&\mathcal{A}_{0,1}^{\ast }\left( \mathcal{T}^{B}\right) \left( s,t\right) :=%
\frac{1}{2}G_{2}^{\ast }\left( s,t\right) -\frac{1}{2}%
\int_{t-1}^{s}B_{v}dB_{v}  \notag \\
&&+\frac{1}{2}B_{t-1}\left( B_{s}-B_{t-1}\right) -\frac{1}{2}B_{s-1}\left(
B_{t}-B_{t-1}\right)  \notag \\
&&-\frac{1}{2}\left( 1-t+s\right)  \notag \\
&&+\frac{1}{2}\hat{G}_{2}^{\ast }\left( s,t\right) -\frac{1}{2}%
\int_{1-s}^{2-t}\hat{B}_{v}d\hat{B}_{v}-\frac{1}{2}\hat{B}_{2-t}\left( \hat{B%
}_{2-s}-\hat{B}_{2-t}\right) ,  \label{continuous area extension}
\end{eqnarray}%
is $\mathbb{P}$-a.s. uniformly continuous on the closed upper triangle 
\begin{equation*}
\left\{ \left( s,t\right) \in \left[ 0,1\right] \times \left[ 0,1\right]
:0\leq s\leq t\leq 1\right\} .
\end{equation*}%
The restriction of $\mathcal{A}_{0,1}^{\ast }\left( \mathcal{T}^{B}\right) $
to the dyadic open upper triangle%
\begin{equation*}
\left\{ \left( s,t\right) \in D^{2}\left( \left[ 0,1\right] ^{2}\right)
:0\leq s<t\leq 1\right\}
\end{equation*}%
agrees with $\mathcal{A}_{0,1}\left( \mathcal{T}^{B}\right) $ on a set $%
\Omega _{0}$ of full $\mathbb{P}$-measure.

To see that $\mathcal{A}_{0,1}^{\ast }\left( \mathcal{T}^{B}\right) $ has a
jump discontinuity on the diagonal suppose that $u_{n}\nearrow u$ such that 
\begin{equation*}
\left( u_{n},u\right) \in \left\{ \left( s,t\right) \in \left[ 0,1\right]
\times \left[ 0,1\right] :0\leq s<t\leq 1\right\}
\end{equation*}%
and 
\begin{equation*}
\left( u,u_{n}\right) \in \left\{ \left( s,t\right) \in \left[ 0,1\right]
\times \left[ 0,1\right] :0\leq t<s\leq 1\right\}
\end{equation*}%
for all $n\in \mathbb{N}$. Since%
\begin{equation*}
B_{u-1}\left( B_{u_{n}}-B_{u-1}\right) -B_{u_{n}}\left( B_{u}-B_{u-1}\right)
\rightarrow 0\text{, }\mathbb{P}\text{-a.s.}
\end{equation*}%
and 
\begin{equation*}
G_{2}^{\ast }\left( u_{n},u\right) \rightarrow \int_{u-1}^{u}B_{v}dB_{v},%
\mathbb{P}\text{-a.s.,}
\end{equation*}%
\begin{equation*}
\hat{G}_{2}^{\ast }\left( u_{n},u\right) \rightarrow \int_{1-u}^{2-u}\hat{B}%
_{v}d\hat{B}_{v},\mathbb{P}\text{-a.s.}
\end{equation*}%
and 
\begin{equation*}
\int_{1-u_{n}}^{2-u}\hat{B}_{v}d\hat{B}_{v}\rightarrow \int_{1-u}^{2-u}\hat{B%
}_{v}d\hat{B}_{v},\mathbb{P}\text{-a.s.,}
\end{equation*}%
we find that $\mathcal{A}_{0,1}^{\ast }\left( \mathcal{T}^{B}\right) \left(
u_{n},u\right) \rightarrow -\frac{1}{2},\mathbb{P}$-a.s. Hence, by (\ref%
{Eqn: antisymmetry}), $\mathcal{A}_{0,1}^{\ast }\left( \mathcal{T}%
^{B}\right) \left( u,u_{n}\right) \rightarrow \frac{1}{2}$, $\mathbb{P}$%
-a.s., which concludes the proof of the fact that $\mathcal{A}_{0,1}^{\ast
}\left( \mathcal{T}^{B}\right) $ has a jump discontinuity on the diagonal.

Now, suppose that $\mathcal{\bar{A}}_{0,1}$ is a continuous extension of $%
\mathcal{A}_{0,1}$ to the unit square $\left[ 0,1\right] \times \left[ 0,1%
\right] $. Since the unit square is compact, $\mathcal{\bar{A}}_{0,1}$ is
uniformly continuous. By the supposed continuity and (\ref{Eqn: antisymmetry}%
), we have that 
\begin{equation*}
\mathcal{\bar{A}}_{0,1}\left( s,s\right) \text{ }\equiv 0\text{ for all }%
s\in \left[ 0,1\right] .
\end{equation*}%
Since $\mathcal{\bar{A}}_{0,1}$ is an extension of $\mathcal{A}_{0,1}$ there
exists a set $\bar{\Omega}_{0}$ of full $\mathbb{P}$-measure such 
\begin{equation*}
\left. \mathcal{\bar{A}}_{0,1}\right\vert _{\left\{ \left( s,t\right) \in
D\left( \left[ 0,1\right] ^{2}\right) :0\leq s<t\leq 1\right\} }=\mathcal{A}%
_{0,1}\text{ on }\bar{\Omega}_{0}\text{.}
\end{equation*}%
But 
\begin{equation*}
\left. \mathcal{A}_{0,1}^{\ast }\right\vert _{\left\{ \left( s,t\right) \in
D\left( \left[ 0,1\right] ^{2}\right) :0\leq s<t\leq 1\right\} }=\mathcal{A}%
_{0,1}\text{ on }\Omega _{0}\text{.}
\end{equation*}%
The uniqueness of the uniformly continuous extension implies that 
\begin{equation*}
\mathcal{A}_{0,1}^{\ast }=\mathcal{\bar{A}}_{0,1}\text{ on }\Omega _{0}\cap 
\bar{\Omega}_{0}\text{,}
\end{equation*}%
which cannot be since $\mathcal{A}_{0,1}^{\ast }$ does \emph{not }vanish on
the diagonal, whereas $\mathcal{\bar{A}}_{0,1}$ does.
\end{proof}

\end{document}